\documentclass[11pt,a4paper,reqno]{amsart}

\usepackage[colorlinks,bookmarksopen,bookmarksnumbered,citecolor=red,urlcolor=red]{hyperref}
\usepackage{amsaddr}
\usepackage{epsfig,amsmath,amssymb,amsfonts,mathrsfs}
\usepackage{graphicx}
\usepackage{float} 
\usepackage{textcomp, eufrak}
\usepackage{verbatim}
\usepackage{fancyhdr}
\usepackage[numbers]{natbib}
\usepackage[dvipsnames]{xcolor}
\usepackage{soul}

\usepackage{caption, threeparttable}
\usepackage{enumitem}
\usepackage{scalerel}

\usepackage{textgreek}

\numberwithin{equation}{section}
\usepackage[mathscr]{euscript}

\newtheorem{theorem}{Theorem}[section]
\newtheorem{definition}[theorem]{Definition}
\newtheorem{lem}[theorem]{Lemma}
\newtheorem{prop}[theorem]{Proposition}
\newtheorem{rem}[theorem]{Remark}
\newtheorem{cor}[theorem]{Corollary}

\newtheorem{ass}[]{Assumption}
\numberwithin{equation}{section}

\newcommand{\cccdot}{\scaleobj{1.25}{\ccdot}}
\newcommand{\brr}{\mathring{b}}

\newcommand{\R}{\mathbb{R}}

\newcommand{\Rn}{\mathbb{R}^n}
\newcommand{\N}{\mathbb{N}}

\newcommand{\Rp}{\mathbb{R}^+}
\newcommand{\F}{\mathcal{F}}

\newcommand{\PS}{(\Omega, \mathcal{F}, \mathbb{P})}

\newcommand{\I}{\mathbb{I}}

\newcommand{\Bb}{\mathcal{B}}

\newcommand{\D}{\mathcal{D}}

\newcommand{\T}{\mathbb{T}}

\newcommand{\p}{\mathbb{P}}

\newcommand{\bn}{\begin{definition}}
\newcommand{\en}{\end{definition}} 
\newcommand{\bt}{\begin{theorem}}                
\newcommand{\et}{\end{theorem}}
 \newcommand{\bnm}{\begin{enumerate}}              
\newcommand{\enm}{\end{enumerate}}
\newcommand{\br}{\begin{rem}} 
\newcommand{\er}{\end{rem}}

\newcommand{\om}{\omega}

\newcommand{\Om}{\Omega}

\newcommand{\btm}{\begin{itemize}}
 \newcommand{\etm }{\end{itemize}}

\newcommand{\E}{\mathbb{E}}

\newcommand{\Rd}{\R^d}

\newcommand{\Ms}{{\mathfrak{P}\hspace{.01cm}}}

\newcommand{\Df}{\D_{\varphi}}

\newcommand{\Ic}{\mathcal{I}}

\newcommand{\DK}{\mathcal{D}_{\textsc{kl}}}

\newcommand{\Solm}{\phi_{t_0,t}}
\newcommand{\LG}{\mathcal{L}}

 \newcommand{\RL}{\mathcal{R}}

\newcommand{\Xt}{\mathcal{M}}

\newcommand{\ccdot}{\,\cdot\,}

\newcommand{\PP}{\hspace*{.000cm}\mathord{\raisebox{-0.095em}{\scaleobj{.86}{\includegraphics[width=1em]{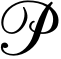}}}}\hspace*{.02cm}}

\begin{document}
\title[Path-based expansion rates and Lagrangian uncertainty in stochastic flows]{Measures of path-based nonlinear expansion rates and Lagrangian uncertainty in stochastic flows}

\author{Micha\l \;Branicki$^{\,\dagger,*}$ and Kenneth Uda$^\dagger$}

\address{\vspace*{-0.2cm} $^\dagger$ Department of Mathematics, University of Edinburgh, Scotland, UK \\
\hspace*{-1.9cm}$^*$ The Alan Turing Institute for Data Science, London, UK}

\email{M.Branicki@ed.ac.uk, K.Uda@ed.ac.uk}
\thanks{This work was supported by the Office of Naval Research grant ONR N00014-15-1-2351.}

\begin{abstract} 
We develop a probabilistic characterisation of trajectorial expansion rates in non-autonomous stochastic dynamical systems that can be  defined over a finite time interval and used for the subsequent  uncertainty quantification in Lagrangian (trajectory-based) predictions. These expansion rates are quantified via certain divergences (pre-metrics) between probability measures induced by the laws of the stochastic flow associated with the underlying dynamics.   We  construct scalar fields of finite-time divergence/expansion rates,  show their existence and space-time continuity for general stochastic flows. Combining these divergence rate fields  with information inequalities derived in \cite{branuda_luq} allows for quantification  and mitigation of  the uncertainty in path-based observables estimated from  simplified models in a way that is amenable to algorithmic implementations, and it can be  utilised in information-geometric  analysis of statistical estimation and inference, as well as in a data-driven machine/deep learning of coarse-grained models. We also  derive a  link between the divergence rates  and  finite-time Lyapunov exponents for probability measures and for path-based observables.

 
\end{abstract}

\maketitle

\newpage

\tableofcontents

\section{Introduction}
\newcommand{\tint}{[t_0,\,T\hspace{.02cm})}

Consider a  dynamical system on a smooth finite-dimensional manifold $\Xt$ generated by the map, $\phi_{t_0,t}\,{:}\;\Xt\,{\times}\, \Omega \rightarrow \Xt$, $t\in \Ic:=\tint$, such that 
\begin{equation}\label{dynsys}
\phi_{t_0,t_0}(\ccdot,\om)=\textrm{id}_\Xt, \quad \phi_{t_0,t}(\ccdot,\om) = \phi_{u,t}(\phi_{t_0,u}(\ccdot,\om),\om)  \quad \textrm{a.a.} \;\,\om\in \Om \quad \forall\, t,u\in \Ic,
\end{equation}
where $(\Om,\mathcal{F},\p)$ is a complete probability space. It is a basic fact that  for any fixed $y$, $x+y\in \Xt$, and the paths $t\mapsto \phi_{t_0,t}(x,\om)$, $t\mapsto \phi_{t_0,t}(x+y,\om)$ the distance  $\|\phi_{t_0,t}(x,\om)-\phi_{t_0,t}(x+y,\om)\|$ can be large (on the scale of $\Xt$) after some finite time even if $|y|\ll1$. More importantly, it is also well known that for two distinct measurable families of maps $\{\phi_{t_0,t}^\mu\}_{t\in \Ic}$, $\{\phi^\nu_{t_0,t}\}_{t\in \Ic}$ generated by difference or differential equations the paths $t\mapsto \phi^\mu_{t_0,t}(x,\om)$, $t\mapsto \phi^\nu_{t_0,t}(x,\om)$ can be very different so that  $\|\phi_{t_0,t}^\mu(x,\om)-\phi_{t_0,t}^\nu(x,\om)\|\gg 1$ $\p$\,-\,a.s.~for a non-negligible set of initial conditions $x\in \Xt$ even if the vector fields generating these maps are arbitrarily $L^p$-close (but not identical);  see, e.g., \cite{kolm64,arnld64,macKay87,wiggins92,katok97,ottino89,bourgain00})  amongst a vast literature concerned with dynamical systems theory. In other words, for a generic dynamical system even small uncertainties in the initial condition or in the underlying dynamics are  likely to result in a drastically different fate of the corresponding trajectories/paths and thus lead  to significant errors in predictions based on the trajectories of the underlying dynamical system.

Given the paths, $ t\mapsto \phi_{t_0,t}(x,\om)$, $\om\in \Om$, labelled by the (potentially uncertain) initial conditions $x\in \Xt$, and a measurable functional $f$ on these paths,  we refer to the estimation of {\it observables} $x\mapsto \E\big[f\big(\phi_{t_0,\scaleobj{1.25}{\ccdot}}(x)\big)\big]$ as a {\it Lagrangian prediction}.  In contrast,  estimation of  vector fields generating  the family of maps $(\phi_{t_0,t})_{t\in \Ic}$  is referred to as an {\it Eulerian prediction}; see~\cite{branuda_luq}.  
This work is motivated by the desire to devise a systematic approach to  improving path-based/Lagrangian predictions in complex dynamical systems based on simplified, data-driven models. Thus, our main objective, following on from \cite{branuda_luq}, is to develop a general  framework for quantification and mitigation of error in  Lagrangian predictions which are obtained from (Eulerian) vector fields generating the underlying dynamical system in a way which naturally applies to both  the deterministic and stochastic settings.  Importantly, as highlighted above, the effects  of uncertainties on the skill of Lagrangian predictions are distinctly different from uncertainties affecting approximations of the Eulerian fields (generating the maps $\{\phi^\nu_{t_0,t}\}_{t\in \Ic}$) and different `metrics' have to be developed for estimating and mitigating  the error in path-based observables. These issues  are naturally addressed  in the probabilistic setting, where  information-theoretic/geometric  tools can be used  to bound and minimise the lack of information in Lagrangian estimates obtained from uncertain Eulerian fields. The following major steps are needed to achieve this objective:
\begin{itemize}[leftmargin=0.8cm] 

\vspace{-0.0cm}\item[\bf (i)] Determination of an appropriate measure of discrepancy between two Lagrangian predictions which is 
 analytically and computationally tractable, and such that it can be  utilised in information-geometric  analysis of statistical estimation and inference on families of models.

\smallskip
\item[\bf (ii)] Derivation of bounds on the error in  Lagrangian observables estimated from approximate models of the original dynamics, and bounding errors in the underlying probability measures. 

\smallskip
\item[\bf (iii)] Identification of the most important Lagrangian structures, i.e., subsets of true and approximate trajectories  which need to be  systematically  `tuned'  in order to minimise the loss of information in Lagrangian predictions based on  approximate Eulerian fields.
\end{itemize}
A framework addressing (i)-(ii)  was developed in \cite{branuda_luq}.  Here, we address step (iii) above and consider fields of trajectorial expansion/divergence rates (around every $x\in \Xt$)  which utilise the evolution of time-marginal probability associated with the (path space) law of the stochastic flow. It turns out minimising the discrepancy between such fields generated by  the true dynamics and its approximations  can be used to systematically bound the uncertainty the Lagrangian predictions and minimise the error in path space observables. Importantly, this procedure can be carried out within the statistical estimation/inference framework or within a purely data-driven setting. As illustrated later on, the features which models have to reproduce for best Lagrangian prediction skill are generally very different from the Eulerian features one would optimise for when predicting, say, the evolution of velocity fields in atmosphere-ocean science applications. Apart from purely theoretical aspects relevant for finite-dimensional dynamical systems (e.g., \cite{macKay87,wiggins92,katok97,ottino89,bourgain00,Arnold86,Baxendale1,Baxendale2,carverhill85,Crauel-89}), 
Hamiltonian PDEs, or some parabolic PDEs (e.g., \cite{constantin01a,constantin01b}), a major motivation for these efforts arose from the desire to provide Lagrangian/path-based predictions in geophysical flows or in molecular dynamics.

\smallskip
In what follows, we first construct scalar fields of finite-time expansion/divergence rates based on the evolution of time-marginal probability measures associated with projections of path space measures induced by the underlying dynamics.  
We then focus on the analytical issues related to the existence and properties of such  fields. At this stage we do not consider the Lagrangian `tuning' of the divergence rate fields, nor do we focus on approaches to estimating the fields generated by the true dynamics which is, in  applications, likely to be  only sampled through partial and noisy observations. These practically important but also numerically involved issues are postponed to a subsequent publication devoted to applications. The use of the bound on the discrepancy between the divergence rate fields as a loss function in machine/deep learning of coarse-grained models is of special interest and is a subject on an ongoing work.

The contents of this article are as follows. Section \ref{mreslt} outlines the general setup, the hierarchy of main results and notation. In section \ref{s_setup} we introduce the main tools and background results which are needed for the subsequent derivations. Section \ref{s_exp_rates} is devoted to the main results. The finite-time divergence rates ($\varphi$-DR) are defined in \S\ref{ftdr_sec} and their general properties are outlined in \S\ref{ss_exrates_prop}. Fields of finite-time divergence rates ($\varphi$-DRF) are discussed in \S\ref{uni_exp_rates}. A tangential but interesting issue concerned with the links between $\varphi$-DRF  and the fields of finite-time Lyapunov exponents is discussed in \S\ref{s_lyp_fun_stoch}.  Section~\ref{s_comp}  is devoted to numerical examples 
illustrating the finite-time divergence rate fields for several dynamical systems. 
We close with some remarks on future work in section~\ref{s_concl}. Several technical proofs are delegated to the Appendix to improve the read.


\bigskip
\section{Main results}\label{mreslt}

The setup and the main results of this work are outlined below.  Frequently used notation is summarised in the Glossary at the end of the paper. 

\subsection{Basic setup and notation}\label{setp}

 For the most part, either for convenience or by necessity, we will assume that the dynamics is  generated  by
 the Stratonovich\footnote{\,We start from the Stratonovich form of the SDE rather than the It\^o form, since the former is consistent with the physical limit that  leads to idealised, i.e. stochastic, perturbations in the deterministic dynamics (e.g., \cite{Gardiner10}).} stochastic differential equation (SDE)  on~$\mathcal{M}$, that is 
\begin{equation}\label{gen_sde}
d X_t^{t_0,x} = \brr\big(t,X_t^{ t_0,x}\big)d t+\sigma\big(t,X_t^{t_0,x}\big)\circ d W_{t-t_0},\qquad X_{t_0}^{t_0,x}  \sim \mu_{t_0}\in\PP(\Xt),\quad t\in\Ic:=\tint,
\end{equation}    
where $W_t$ is an m-dimensional  Wiener process and  the uncertainty in the initial conditions, or a distribution of initial conditions of interest, is prescribed by $\mu_{t_0}\in \PP(\Xt)$, where $\PP(\Xt)$ is the space of Borel probability measures on $\Xt$. Throughout,   $\Xt = \R^d$ or a flat torus $\Xt\,{=}\,\bar{\T}^d$ and all probability measures defined on the appropriate Borel $\mathfrak{S}$-algebra. In this case $\PS$ associated with (\ref{dynsys}) is the Wiener probability space   and  the law of $X_t^{t_0,x}$ is given by $\Ms_{t_0} \simeq \mu_{t_0}\otimes \p$. Furthermore, $\Om\simeq\mathcal{C}_{0}(\R;\R^m)$ is identified with a subspace of continuous functions $\mathcal{C}(\R;\R^m)$ which vanish at zero (e.g., \cite{Arnold1}), 
and $\p$ is the Wiener measure on $\mathcal{F}$, where $\F$ is the Borel  $\mathfrak{S}$-algebra on~$\Om$ generated by the Wiener process. In what follows, we will consider   $(\brr, \sigma)$ such that (\ref{gen_sde}) has unique global solutions  which can be represented by stochastic flows (\S\ref{s_setup}); namely 
 $$X^{t_0,x}_t(\om) = \phi_{t_0,t}(x,\om)\quad \p\,\textrm{-\,a.s.}\quad t\in \Ic, \;x\in \Xt. $$
 Under some weak assumptions specified   in~\S\ref{s_setup} the family of maps $\big\{\phi_{t_0,t}(\ccdot,\ccdot)\big\}_{t\in \Ic}$ represents  a {\it stochastic flow} of homeomorphisms on $\Xt$; we will  exploit this fact extensively. Many results derived in the sequel apply to a broader class of flows than those induced by solutions of SDEs in the  Markovian setting but the outlined case serves as a useful reference setup. 
 
Given the family of time-marginal probability measures $(\mu_t)_{t\in \Ic}$, $\mu_t\in \PP(\Xt)$ generated, for example,  by the law of $\phi_{t_0,\cccdot}$ on $\Xt$, a unified probabilistic approach to studying trajectorial expansion rates and uncertainty associated with the paths $x\mapsto\phi_{t_0,t}(x,\omega)$ is based on utilising  the so-called $\varphi$-divergencies $D_\varphi(\mu_t\|\mu_{t_0})$ which are defined as 
\begin{align}\label{phidiv}
\Df(\mu_t\|\mu_{t_0}) = \scaleobj{.9}{\int_\Xt}\,\varphi\big(d\mu_t/d\mu_{t_0}\big)d\mu_{t_0}, 
\end{align}
where $d\mu_t/d\mu_{t_0}$ is the Radon-Nikodym derivative, and $\varphi\in \mathcal{C}^2(\R^+)$ is a strictly convex scalar function (see \S\ref{s_exp_rates}).  This probabilistic notion of trajectorial expansion rates has an information-theoretic interpretation as the loss of information in the measure $\mu_{t_0}$ on the initial conditions  for approximating~$\mu_{t}$ and it is used in  \cite{branuda_luq} for Lagrangian uncertainty quantification in path-based predictions. Here, we further develop the link between Lagrangian uncertainty  and certain fields of $\varphi$-divergences which can be utilised in  data-driven model learning and optimization.

\subsection{Outline of main results}\label{out_main}
\noindent  Our results  can be summarised as follows:
\begin{itemize}[leftmargin = .85cm]
\item[{\bf (I)}] We define and study nonlinear expansion rates, termed  the {\it $\varphi$-divergence rates} ($\varphi$-DR),
 which are based on the evolution of time-marginal probability measures $(\mu_t)_{t\in \Ic}$, ${\mu_t\in\PP(\Xt)}$ in terms of  $\D_\varphi(\mu_t\|\mu_{t_0})$.  Then, we construct scalar $\varphi$-DR fields, $x\mapsto \D_\varphi^{t_0,t}(\mu^x_t\|\mu^x_{t_0})$, where $\mu^x_{t_0}$ is a  probability measure localised around  $x\in \Xt$, and $\D^{t_0,t}_\varphi:=|t-t_0|^{-1}\D_\varphi$. These fields are defined for any family $(\mu_t)_{t\in \Ic}$ such that $\mu_t\ll\mu_{t_0}$ but we focus on time-marginal measures of the (path space) law $\Ms_{t_0}$ of the stochastic flow $\big\{\phi_{t_0,t}\big\}_{t\in \Ic}$ in which case $\Ms_{t_0}\circ\phi^{-1}_{t_0,t} = \mu_t$, thus providing a trajectorial interpretation. In section \ref{s_exp_rates} (Theorem \ref{exp_rate}) we show  that  the $\varphi$-DR fields exist and are space-time continuous under fairly general conditions, and that they provide a computable diagnostic of nonlinear expansion induced by the stochastic flow.

Importantly, $\varphi$-DR fields  can be viewed as a general measure-based diagnostic of the  growth of uncertainty in stochastic 
flows that operates  in  both the deterministic and the stochastic settings. 
The information-theoretic loss of information in $\mu^x_t\in \PP(\Xt)$ for describing the initial probability measure $\mu^x_{t_0}$ is used  in  \cite{branuda_luq} for uncertainty quantification in path-based predictions of statistical observables.
These results  rely on the hierarchy of bounds on statistical observables based on two different flows, $\{\Solm^\mu\}_{t\in \Ic}$ and $\{\Solm^\nu\}_{t\in \Ic}$, which employ  general {\it information bounds/inequalities} in the form\footnote{\, In fact, the information inequalities also hold for path space measures $\mathfrak{P}^x_{t_0}$, $\mathfrak{Y}^x_{t_0}$ s.t.~$\mathfrak{P}^x_{t_0}\ll \mathfrak{Y}^x_{t_0}$ and $\mathcal{D}_\varphi^{t-t_0}(\Ms^x_{t_0}\|\mathfrak{Y}^x_{t_0})$ can be bounded by a functional of the rhs of (\ref{bnd2}) but here we focus on time-resolved observables.} 
  \begin{align}\label{bnd1}
\hat{\mathcal{K}}_{\varphi,f}^{\nu}\big(-\mathcal{D}_\varphi(\mu_t\|\nu_t)\big)\leqslant  \E^{\mu_t}[f] - \E^{\nu_t}[f] \leqslant \mathcal{K}_{\varphi,f}^{\nu}\big(\mathcal{D}_\varphi(\mu_t\|\nu_t)\big), \quad t\in \Ic, 
\end{align}
where  $\E^{\mu_t}[f] = \E^{\mu_{t_0}\otimes\p}\big[f\big(\Solm^\mu\big)\big]$, $\E^{\nu_t}[f] = \E^{\mu_{t_0}\otimes\p}\big[f\big(\Solm^\nu\big)\big]$ (see (\ref{k_lift})), $\mu_t,\nu_t\in \PP(\Xt)$, and  $\mathcal{K}_{\varphi,f}^{\nu}$, $\hat{\mathcal{K}}_{\varphi,f}^{\nu}$ are bounded and such that $\mathcal{K}_{\varphi,f}^{\nu}(u)\rightarrow 0$, $\hat{\mathcal{K}}_{\varphi,f}^{\nu}(-u)\rightarrow 0$ as $u\,{\downarrow}\, 0$. In particular, setting $f=\textrm{Id}_\mathscr{A}$
corresponds to estimating the fate of trajectories with initial conditions supported on the set $\mathscr{A}$.  Importantly, it turns out \cite[\S5.2]{branuda_luq} that for two families of probability measures $(\mu^x_t)_{t\in \Ic}$, $(\nu^x_t)_{t\in \Ic}$ such that $\mu^x_{t_0} = \nu^x_{t_0}$,  $\mu^x_t\ll\mu^x_{t_0}$, $\nu^x_t\ll\mu^x_{t_0}$
  one has 
\begin{align}\label{bnd2}
\mathcal{D}_\varphi^{t-t_0}(\mu^x_t\|\nu^x_t)\leqslant  \big\vert \mathcal{D}_\varphi^{t-t_0}(\mu^x_t\|\mu^x_{t_0})-\mathcal{D}_\varphi^{t-t_0}(\nu^x_t\|\mu^x_{t_0})\big\vert, \quad  \forall \;x\in \Xt, \,t\in \Ic,
\end{align}
so that the combination of (\ref{bnd1}) and (\ref{bnd2}) allows for quantification of uncertainty in approximations of the observables $\E^{\mu_t}[f]$, 
 by considering discrepancies between the $\varphi$-DR fields  $\mathcal{D}_\varphi^{t-t_0}(\mu^x_t\|\mu^x_{t_0})$ and $\mathcal{D}_\varphi^{t-t_0}(\nu^x_t\|\mu^x_{t_0})$ defined over $\Xt$. 
Thus, given a parametric family of models, one can choose (in the simplest case) the model generating $\big(\nu^{x}_t(\eta^\dagger;dy)\big)_{t\in \Ic}$ where
\begin{equation}\label{ml_bnd}
\eta^\dagger = \min_\eta \int_{\Ic}\int_\Xt\big\vert \mathcal{D}_\varphi^{t_0,t}(\mu^x_t\|\mu^x_{t_0})-\mathcal{D}_\varphi^{t_0,t}(\nu^x_t(\eta)\|\mu^x_{t_0})\big\vert^p\kappa(dx)dt, \qquad \kappa\in \PP(\Xt).
\end{equation}
As mentioned earlier, the practically important but also numerically involved issues related to model optimisation are postponed to a separate publication devoted to applications.
In particular, the use of the bound (\ref{ml_bnd}) as a loss function in machine/deep learning of coarse-grained models is of special interest for the subsequent work.

\item[{\bf (II)}]  In addition to the above,  the $\varphi$-DR framework elucidates trajectorial connections between the evolution of time-marginal probability measures on $\Xt$ and average local stretching rates 
obtained from general Lyapunov exponents for probability measures\footnote{\,See Corollary \ref{localbound} and Reemark \ref{lapsum} for introductory results in this direction. The (infinite-time) Lyapunov exponents for probability measures in the ergodic  setting  were considered by, e.g.,  Kunita \cite{Kunitabook}, Arnold \cite{Arnold1}, and  Baxendale \cite{Baxendale1,Baxendale2}; here, we focus on the finite-time case and a different, simpler setting used in applications. 
}.
Moreover,  it turns out that for a specific $\varphi$-divergence, namely the KL-vivergence $\D_\textsc{kl}$, and a probability measure $\mu_t^x$ evolving from  $\mu^x_{t_0}$ concentrated on a neighbourhood of $x\in \Xt$  the map $x\mapsto\DK^{t_0,t}\big({\mu}^x_t\|{\mu}^x_{t_0}\big)$ can be related  to fields of the so-called  finite-time Lyapunov exponents which for deterministic dynamics are defined as  (see \S\ref{s_lyp_fun_stoch})
$$\Lambda^{t-t_0}_{t_0}(x,y) = \scaleobj{.8}{\frac{1}{|t-t_0|}}\log \frac{|Y_{t}^{t_0,y}(x)|}{|y|}, \qquad y\ne0,\;y\in \Xt,$$
and  are commonly used in applications to assess  the linearised  growth of a perturbation about solutions with the initial condition $x\in \Xt$; i.e., for $\mathcal{Y}^{t_0,y}_t(x) = \phi_{t_0,t}(x)-\phi_{t_0,t}(x+y)$ we have ${Y}^{t_0,y}_t(x) = D\phi_{t_0,t}(x)y$. The bounds on the average finite-time Lyapunov exponent turn out to have the same form as those in (\ref{bnd1}); namely
\begin{align}\label{KK_1_summ}
{\mathcal{K}}^{x}_{\scaleobj{1}{-}}\big(-\D_\textsc{kl}^{t_0,t}\big(\tilde{\mu}^{x}_t\|\tilde{\mu}^{x}_{t_0}\big)\big)\leqslant \E^{\tilde{\mu}^{x}_{t_0}}\big[{\Lambda}^{t-t_0}_{t_0}(x)\big] \leqslant \mathcal{K}^{x}_+\big(\D_\textsc{kl}^{t_0,t}\big(\tilde{\mu}^{x}_t\|\tilde{\mu}^{x}_{t_0}\big)\big), \hspace{.8cm} \forall \,t\in \Ic,
\end{align}
where  ${\mathcal{K}}^{x}_{\scaleobj{1}{-}}$, ${\mathcal{K}}^{x}_{\scaleobj{1}{+}}$ are bounded and such that ${\mathcal{K}}^{x}_{\scaleobj{1}{-}}(-u)\rightarrow 0$, ${\mathcal{K}}^{x}_{\scaleobj{1}{+}}(u)\rightarrow 0$ as $u\,{\downarrow}\, 0$, and the locally averaged Lyapunov exponent is given by ${\E^{\tilde{\mu}^{x}_{t_0}}\big[{\Lambda}^{t-t_0}_{t_0}(x)\big] = \int_\Xt {\Lambda}^{t-t_0}_{t_0}(x,y)\mu^x_{t_0}(dy)}$. A more general version for in the fully stochastic case is given in \S\ref{s_lyp_fun_stoch}. These result further highlights the utility and universality of the general information bounds (\ref{bnd1})-(\ref{bnd2}) derived~in~\cite{branuda_luq}.
\end{itemize}
 
 It is worth stressing that this study is aimed systematising the link (\ref{bnd1})-(\ref{bnd2}) between the expected error in path-based predictions and $\varphi$-DR fields in the context of Lagrangian uncertainty quantification \cite{branuda_luq} and our results are predominantly concerned with the existence and properties of $\varphi$-DR fields as measures of trajectorial expansion rates and uncertainty in stochastic dynamics.  We are not concerned with the so-called Lagrangian transport analysis. Thus, we neither discuss nor prove any relationship between $\varphi$-DR fields and (approximately) flow-invariant  structures in stochastic or deterministic flows, which is a topic for another study.

\section{General framework and background results}\label{s_setup}

In order to make the presentation relatively self-contained for the target audience, we gather below several  general results (some well known, some new) and assumptions that are needed in the subsequent sections. 

\begin{definition}[\textbf{\textit{Stochastic flow}} \cite{Kunitabook, Kunitanote}]\label{Flow} \rm 
For any $s,t\in \Ic\subseteq\R$, $ x \in \Xt$, let $\phi_{s,t}(x,\om)\in \Xt$  be a random field on some  probability space $(\Om, \F, \p)$. 
The two-parameter family $\{\phi_{s,t}\!:\, s,t\in \Ic\subseteq\R \}$ is called a {\it stochastic flow of homeomorphisms} if there exists a null set $\mathfrak{N}\subset \Om$ such that for any $\om\,{\notin}\, \mathfrak{N},$ there exists a family of continuous maps $\{\phi_{s,t}(\ccdot,\om)\!:\, s,t\in\Ic\}$ on $\Xt$ satisfying 
\begin{itemize}[leftmargin=0.9cm]
\item[(i)] $\phi_{s,t}(\ccdot,\om) = \phi_{u,t}(\phi_{s,u}(\ccdot,\om),\om)$ holds for any $s,t,u\in \Ic,$
\item[(ii)] $\phi_{s,s}(\ccdot,\om) = \textup{id}_{\Xt}$, for all $s\in\Ic,$
\item[(iii)] $\phi_{s,t}(\ccdot,\om): \Xt\rightarrow \Xt$ is a homeomorphism for any $t,s\in \Ic$.
\end{itemize}
The map $\phi_{s,t}(\ccdot,\om)$ is a {\it stochastic flow of $\mathcal{C}^l$-diffeormorphisms}on $\Xt$, if it is a homeomorphism and $\phi_{s,t}(x,\om)$ is $l$-times continuously differentiable with respect to $x\in \Xt$ for all $s,t\in \Ic\subseteq\R$. 
 The {\it two-parameter filtration} $\{\F_s^{\hspace{.02cm}t}: s\leqslant t\}$  is the smallest complete sub $\mathfrak{S}$-algebra of $\mathcal{F}$ s.t.~
 $\cap_{\varepsilon>0}\,\mathfrak{S}\big(\phi_{u,v}\!:\; s-\varepsilon\leqslant u, v\leqslant t+\varepsilon\big)$.
The flow is referred to as {\it `forward'} for $s\leqslant t$, and as {\it `backward'} for $t\leqslant s$.  
\end{definition}

\begin{definition}[\textbf{\textit{Transition kernel induced by a stochastic flow}}]\label{Mkernel}\rm
Consider the stochastic flow $\big\{\phi_{s,t}(\ccdot,\ccdot)\!: \,s,t\in \Ic\big\}$.  
 Given the Borel-measurable space $\big(\Xt,\mathcal{B}(\Xt\hspace{.03cm})\big)$, the {\it transition probability kernel} $P(s,x; t, \ccdot)$ induced by the stochastic flow is given~by
 \begin{equation}\label{phiP}
 P(s,x; t, A) = \p\big(\{\om\in\Om: \phi_{s,t}(x,\om )\in A\}\big), \quad \forall\;s,t\in\Ic,\; A\in\Bb\big(\Xt\hspace{.03cm}\big).
 \end{equation}
 The transition kernel generated by a forward stochastic flow  satisfies the Chapman-Kolmogorov equation (e.g., \cite{Kunitabook,bogachev10,Roc-Kry}); note that $\{\phi_{s,t}\}_{s,t\in \Ic}$ is not required to be Markov w.r.t~$(\F_s^{\hspace{.02cm}t})_{s\leqslant t}$. 
 \end{definition}
 \begin{definition}[\textbf{\textit{Transition evolution and its dual}}]\label{trans_ev}\rm
 Given the stochastic flow $\big\{\phi_{s,t}\}_{s,t\in \Ic}$ and the transition kernel (\ref{phiP}), 
 the operator  $\mathcal{P}_{s,t}: \mathbb{M}_\infty(\Xt)  \rightarrow\mathbb{M}_\infty(\Xt)$, $s,t\in \Ic$, called the {\it transition evolution} is defined by 
 \begin{equation}\label{calP}
 \big(\mathcal{P}_{s,t}f\big)(x) = \int_{\Xt}f(y)P(s,x; t, dy) = \E\big[f\big(\phi_{s,t}(x,\ccdot)\big)\big], \quad \quad \forall\,s,t\in\Ic,\;x\in \Xt.
 \end{equation}
  For any Borel probability measure $\mu_s \in \PP(\Xt)$, $s\in \Ic$, on $\big(\Xt,\Bb(\Xt)\big)$, the (formal) $L^1(\Xt,\mu_s)$ 
  dual $\mathcal{P}_{s,t}^{*}$ of the transition evolution $\mathcal{P}_{s,t}$ is  defined as
 \begin{equation}\label{P*mu}
 \big(\mathcal{P}_{s,t}^{*}\mu_s\big)(A) = \int_{\Xt}P(s,x; t, A)\mu_s(dx),  \qquad  \quad \forall\;s,t\in\Ic,\;A\in \Bb\big(\Xt\big).
 \end{equation}
 Consequently,  for any $s,u,t\,{\in}\, \Ic$ 
  and for all $A\in \Bb(\Xt)$ we can formally write 
 \begin{equation}\label{mut}
 \mu_t(A) = \big(\mathcal{P}_{s,t}^{*}\mu_s\big)(A) = \big(\mathcal{P}_{u,t}^{*}\mathcal{P}_{s,u}^{*}\hspace{.03cm}\mu_s\big)(A) = \big(\mathcal{P}_{u,t}^{*}\hspace{.03cm}\mu_u\big)(A).
 \end{equation}
 \end{definition}

\begin{definition}[\textbf{\textit{Random probability measures}}]\rm
Given the probability measure  $\mu_s\in \PP(\Xt)$ and  the stochastic flow $\big\{\phi_{s,t}(\ccdot,\ccdot)\!: \,s,t\in \Ic\big\}$, the pushforward of $\mu_s$ by, respectively,   $\phi_{s,t}(\ccdot,\om)$ and its inverse $\phi_{s,t}^{-1}(\ccdot,\om):=\phi_{s,t}(\ccdot,\om)^{-1}$ given by    
\begin{align}\label{Random measure1}
\Pi_{s,t}(A,\omega) := (\phi_{s,t}\mu_s)(A) = \mu_s\big(\phi_{s,t}^{-1}(A,\omega)\big), \qquad A\in \Bb(\Xt),
\end{align}
\begin{equation}\label{imm}
\Check{\Pi}_{s,t}(A,\omega) := (\phi_{s,t}^{-1}\mu_s)(A)= \mu_s\big(\phi_{s,t}(A,\omega)\big), \qquad A\in \Bb(\Xt),
\end{equation}
are referred to as {\it random probability measures} induced by the stochastic flow $\{\phi_{s,t}\}_{s,t\in \Ic}$.  
\end{definition}

\begin{prop} \label{al_pi} Assume that $\mu_s\ll m_d$, $ d\mu_s/dm_d>0$.\footnote{\,Throughout, $\mu\ll m_d$ means that $\mu\in \PP(\Xt)$ is absolutely continuous with respect to $m_d$ - the Lebesgue measure on $\Rd$; $d\mu/dm_d$ is the Radon-Nikodym derivative.} If  
$\big\{\phi_{s,t}\}_{s,t\in \Ic}$ is a stochastic flow of diffeomorphisms on $\Xt$, then the densities  $\pi_{s,t} = d\Pi_{s,t}/d\mu_s$ and $\alpha_{s,t} = d\Check{\Pi}_{s,t}/d\mu_s$ 
are related~by 
\begin{align}\label{rel11}
 \pi_{s,t}(x,\omega)^{-1} = \alpha_{s,t}(\phi_{s,t}^{-1}(x,\omega),\omega) \qquad \p\,\textrm{\rm -\,a.s.} 
 \end{align}
\end{prop}
\noindent {\it Proof.} See Appendix  \ref{app_al_pi}. 

\smallskip

 \vspace{-0.1cm}
 \begin{theorem}\label{mut_rnd} 
  Assume that $d\mu_s/dm_d>0$  and $\big\{\phi_{s,t}\}_{s,t\in \Ic}$ is a stochastic flow of diffeomorphisms on $\Xt$.
Then, the time-marginals $ \mu_t=\mathcal{P}_{s,t}^{*}\mu_s $ of the law of $\phi_{s,\scaleobj{.9}{\Ic}}$   are such that  $\mu_{t}\ll m_d$ $\forall t\in \Ic$ and the density $\rho_t =  d\mu_t/dm_d>0$  is given~by
  \begin{align}\label{rho_pi}
 \rho_{t}(x)  = \rho_{s}(x)\E\big[\pi_{s,t}( x)\big]. 
\end{align}  
 \end{theorem}
 
 \vspace{-.2cm}
 \noindent {\it Proof.} See Appendix \ref{app_mutrnd}. 
 
 \begin{rem}\rm
The random, path-dependent measures (\ref{Random measure1}) and (\ref{imm}) are needed in the subsequent uncertainty quantification in Lagrangian considerations, since they  enable a unified treatment of both the deterministic and stochastic dynamics.  The relationship in (\ref{rel11}) is important, since  
 it is  generally not possible to write the stochastic integral governing $\pi_{s,t}$, as opposed to $\alpha_{s,t}$ (see (\ref{alph_Ga})), and both are needed later on.
\end{rem}
Consider the dynamics (\ref{gen_sde}) in the It\^o form either on $\Xt = \Rd$ or $\Xt = \bar{\mathbb{T}}^d$ given by 
\begin{equation}\label{gen_sde_ito}
d X_t^{t_0,x} = b\big(t,X_t^{ t_0,x}\big)d t+\sigma\big(t,X_t^{t_0,x}\big)d W_{t-t_0},\qquad X_{t_0}^{t_0,x}  \sim \mu_{t_0}\in\PP(\Xt),\quad t\in\Ic,
\end{equation} 
where $ b_i(t,x) \,{:=}\, \brr_i(t,x)\,{+}\,c_i(t,x)$, $c_i(t,x) := \frac{1}{2}\sum_{k,j=1}^{m,d} \sigma_{jk}(t,x)\partial_{x_j} \sigma_{ik}(t,x)$ is the Stratonovich correction, and $W_t$ is an $m$-dimensional Wiener process. 

\begin{theorem}[\textbf{\textit{Representation of solutions to SDE's via stochastic flows}}]\label{SDE_flow}\mbox{}

\noindent 
Consider the dynamics (\ref{gen_sde_ito}) with $b(t,\ccdot):\Xt\rightarrow\Xt$, and $\sigma(t,\ccdot): \Xt\rightarrow \Xt^{\otimes m}$  measurable functions for all $t\in \Ic =[t_0,\,T)$, and $t\mapsto b(t,\ccdot)$, $t\mapsto \sigma(t,\ccdot)$ continuous on $\Ic$ such that\,\footnote{\,Throughout $|\cdot|$ is the Euclidean norm on $\Xt$, and  $\|\cdot\|_\textsc{hs}$ is the Hilbert-Schmidt  norm on $\Xt^{\otimes m}$.} 
\begin{alignat}{2}
&|\langle b(t,x),x\rangle|+\|\sigma(t,x)\|^2_\textsc{hs}\leqslant C(1+|x|^2),& \hspace{.8cm} &x\in \Xt, \,t\in \Ic, \label{k_growth}\\[.1cm]
&\vert b(t, x) - b(t,y)\vert + \Vert \sigma(t,x)-\sigma(t,y)\Vert_{\textsc{hs}}\leqslant L_K\vert x-y\vert, && x,y\in \Xt, \,t\in \Ic, \;L_K>0, \label{Lipsch}
\end{alignat}
where  $K$ is any compact subset of $\Xt$.
Then, the  solutions to (\ref{gen_sde_ito})  can be represented (have a continuous modification)  as a stochastic flow of homeomorphisms, i.e.,  
\begin{equation}\label{SDE_sol}
X^{t_0,x}_t(x,\om) = \phi_{t_0,t}(x,\om) \quad \p\,\textrm{\rm -\,a.s.}, 
\end{equation}
see, e.g., \cite[Theorem 4.7.1 combined with Theorem 3.4.6]{Kunitabook}. Moreover,   $\phi_{t_0,t}(\ccdot,\om)$ in  (\ref{SDE_sol}) is a $C^l$-diffeomorphism on $\Xt$ if, in addition to (\ref{k_growth})-(\ref{Lipsch}), the $l$-th derivatives of $b(t,\ccdot)$ and $\sigma(t,\ccdot)$
 are  $\delta$-H\"older continuous for all $t\in \Ic$ and $0<\delta\leqslant 1$ (e.g., \cite[Theorem 4.7.2 combined with Theorem 3.4.6 giving global solutions on $\Ic$]{Kunitabook}). The same holds for solutions of (\ref{gen_sde}) with $l$~derivatives of $\brr(t,\ccdot)$ and $l+1$ derivatives of $\sigma(t,\ccdot)$ $\delta$-H\"older continuous for all $t\in \Ic$, $\delta >0$.

\end{theorem}

\begin{ass}\label{cf_remk}\rm
Specific function spaces of interest that  contain the coefficients of (\ref{gen_sde}) or (\ref{gen_sde_ito}) that generate flows of homeomorphisms/diffeomorphisms are denoted by $\tilde{\mathcal{C}}^{l,\delta}(\Xt)$, $\bar{\mathcal{C}}^{l,\delta}(\Xt)$,
 $l\in \mathbb{N}_0$, $0<\delta\leqslant 1$ (see the Glossary). The solutions of (\ref{gen_sde})  are represented by a flow of $C^l$-diffeomorphisms~for  $t$-continuous $(\brr,\sigma)$ such that
  \begin{equation}\label{reg_growth_cond}
\brr(t,\ccdot)\in \tilde{\mathcal{C}}^{l,\delta}(\Xt), \qquad \sigma_k(t,\ccdot)\in \bar{\mathcal{C}}^{l+1,\delta}(\Xt), \qquad l\geqslant 2, \; k=1,\dots,m, \; t\in \Ic, 
\end{equation}
where  $\sigma_k$, $1\leqslant k \leqslant m$, are the columns of $\sigma$.
For the solutions of (\ref{gen_sde_ito}) one needs  
\begin{equation}\label{reg_growth_cond_ito}
b(t,\ccdot)\in \tilde{\mathcal{C}}^{l,\delta}(\Xt), \qquad \sigma_k(t,\ccdot)\in \bar{\mathcal{C}}^{l,\delta}(\Xt), \qquad l\geqslant 2, \; k=1,\dots,m, \; t\in \Ic.
\end{equation}
We will assume throughout that (\ref{reg_growth_cond}) holds since it implies (\ref{reg_growth_cond_ito}).
 \end{ass}
\begin{rem}\rm
Largely analogous results to those in Theorem \ref{SDE_flow} hold for SDE's driven by more general semimartingales than $F(t,x) = \int_{t_0}^t b(s,x)dt+\int_{t_0}^t  \sigma(s,x)dW_s$ as long as their mean is in $ \tilde{\mathcal{C}}^{l,\delta}(\Xt)$ and the covariance is in $\bar{{\mathcal{C}}}^{l,\delta}(\Xt)$ (see \cite[\S4.4]{Kunitabook}).  
\end{rem}

 If the flow $\big\{\phi_{t_0,t}\}_{t\in \Ic}$ is generated by the SDE  (\ref{gen_sde}) or (\ref{gen_sde_ito}) with sufficiently regular coefficients (e.g., (\ref{reg_growth_cond_ito}) or (\ref{reg_growth_cond})),  the time-marginal probability measures  $(\mu_t)_{t\in \Ic}$ of the law of $\phi_{t_0,\scaleobj{.9}{\Ic}}$  satisfy uniquely (in the distributional sense) the {\it forward Kolmogorov equation} (e.g., \cite{Roc-Kry, Figali, Stroock79})
\begin{align}\label{fPDE}
\partial_t\mu_{t}=\LG_t^*\mu_{t},    \qquad \mu_{t_0} \in \PP(\Xt), \qquad t \in \Ic, 
\end{align}
where $\LG^*_t$ is the formal $L^1(\Xt,\mu_{t})$ dual of
\begin{equation}\label{gen}
 \LG_t f(x) = \sum_{i=1}^\ell b_i(t,x)\partial_{x_i}f(x)+\frac{1}{2}\sum_{i,j=1}^\ell a_{ij}(t,x)\partial^2_{x_i x_j}f(x), \qquad f\in \mathcal{C}^{2}(\Xt),
\end{equation}
with $ b_i(t,x) \,{:=}\, \mathring{b}_i(t,x)\,{+}\,c_i(t,x)$, $c_i(t,x) := \frac{1}{2}\sum_{k,j=1}^{m,d} \sigma_{jk}(t,x)\partial_{x_j} \sigma_{ik}(t,x)$, and $a_{ij}:=  \sum_{k=1}^m \sigma_{ik}\sigma_{jk}$. When  $\big\{\mathcal{P}_{t_0,t}\big\}_{t\in \Ic}$  are generated by the flow of solutions of  (\ref{gen_sde}) or (\ref{gen_sde_ito}),  $\mu_t = \mathcal{P}^*_{t_0,t}\mu_{t_0}$  solves~(\ref{fPDE}) and it can be represented as (e.g., \cite[Theorem~2.6]{Figali} for finite measures or~\cite{Roc-Kry})
\begin{align}\label{k_lift}
 \int_{\Xt}f(x)\mu_t(dx) = \int_{\Xt}\int_\Om f\big(\phi_{t_0,t}(x,\om)\big)\Ms_{t_0,x}(d\om)\mu_{t_0}(dx), \qquad \forall\;f\in \mathcal{C}_{\infty}^2(\Xt), 
 \end{align}
 where $\Ms_{t_0,x}$ is the martingale solution for $\mathcal{L}_t$ starting at $x\in \Xt$ (e.g., \cite{Stroock79}); here, given the existence of the flow generated by (\ref{gen_sde}), $\Ms_{t_0,x}\circ \phi^{-1}_{t_0,t}(x,\ccdot) = \p\circ \phi^{-1}_{t_0,t}(x,\ccdot)$.

 In many cases the forward Kolmogorov equation~(\ref{fPDE}) can be written in terms of the  density of $\mu_{t}$ w.r.t.~the Lebesgue measure $m_d$ on $\mathcal{M}$,  
 and one seeks the (distributional) solutions to (e.g., \cite{Roc-Kry, Figali, Stroock79}) 
\begin{align}\label{rho_FP}
\partial_t\rho_{t}(x) = \LG_t^*\rho_{t}(x), \qquad  \rho_{t_0}\in L^1_+(\Xt;m_d)\cap L^\infty(\Xt;m_d). 
\end{align}

\addtocontents{toc}{\protect\setcounter{tocdepth}{2}}

\smallskip
In the sequel we will need the following result which links the evolution of random probability measures $\Pi_{t_0,t}$ in (\ref{Random measure1}) transported by the stochastic flow $\big\{\phi_{t_0,t}\big\}_{t\in \Ic}$ of  solutions of the SDE~(\ref{gen_sde}) to  the solutions  of the forward Kolmogorov~equation (\ref{fPDE}): 

 \vspace{-0.1cm}
 \begin{theorem}[\textbf{\textit{Solutions of PDE}} (\ref{fPDE})]\label{FKPt1} Suppose the coefficients in  (\ref{gen_sde}) satisfy Assumption~\ref{cf_remk}.
Then, for any $\mu_{t_0}\in\PP(\Xt)$  
the PDE (\ref{fPDE}) has a unique $($weak$-^*$$)$ solution ${\mu_{t}\in \PP(\Xt)}$ on $\Ic$. 
 If $\mu_{t_0}\ll m_d$ and the  (Lebesgue) density $\rho_{t_0} \in \mathcal{C}^2_\infty(\Xt; \Rp)\cap L^1_+(\Xt;dx)$, the solution to~(\ref{fPDE}) is absolutely continuous with respect to $m_d$ with $\rho_{t} \in  \mathcal{C}^2_\infty(\Xt; \Rp)\cap L^1_+(\Xt;dx)$, $t\in \Ic$, and it  coincides with ~(\ref{rho_pi}). If, in addition, $\sigma\sigma^*$  is uniformly elliptic, then (\ref{rho_FP}) has a unique classical solution  $\rho_{t}=d\mu_t/dm_d\in \mathcal{C}_\infty^{2}(\Xt;\R^+)\cap L^1_+(\Xt;dx)$, $t>t_0$,  even if the density of the initial measure $\mu_{t_0}$ is a Dirac mass (in the sense of Schwarz distributions). 

\vspace{.1cm}
\noindent {\it Proof.} {\rm These are  well known results (e.g.,~\cite{Roc-Kry,Stroock79,Figali,leBris08,Kunitabook}) but additional comments on this compressed statement are listed  in Appendix \ref{app_FKPt1}.} 
 \end{theorem}

\begin{theorem} {\rm (\textbf{\textit{Derivative flow}}  \cite[Theorem 3.3.4 and  Corollary 4.6.7]{Kunitabook}}\label{der_fl_thm}
If the coefficients of the SDE~(\ref{gen_sde}) satisfy Assupmtion \ref{cf_remk}, there exists a (stochastic) {\it derivative flow}  of diffeomorphisms $\big\{D\Solm\big\}_{t\in \Ic}$~on~$\Xt$\footnote{\,Note that $D\phi_{t_0,t}(x,\om):\; T_x\Xt \rightarrow T_{\phi_{t_0,t}(x,\om)}\Xt$ but we utilise the isomprphism between  $T_{x}\Xt$ and $\Xt$ given the assumed `flat' geometry of $\Xt$.}, associated with the flow $\big\{\Solm\big\}_{t\in \Ic}$ of solutions of (\ref{gen_sde}) 
 that satisfies 
 \begin{align}\label{var_dphi}
 D\Solm(x,\om) &= {\rm id}_{\mathcal{M}}\;{+}\int_{t_0}^t\big(\nabla_{\!x} \brr\,\big)(\xi,\phi_{t_0,\xi}(x,\om))D\phi_{t_0,\xi}(x,\om)d\xi \notag\\
 &\hspace{2cm}+\sum_{k=1}^d\int_{t_0}^t\big(\nabla_{\!x}\sigma_k\big)(\xi,\phi_{t_0,\xi}(x,\om))D\phi_{t_0,\xi}(x,\om)\circ dW_\xi^k \quad \p\,\textrm{\rm -\,a.s.}
 \end{align}
 \end{theorem}

\begin{definition}[\textit{Centred two-point motion}]\label{cnt_def}\rm
Consider the stochastic flow $\big\{\Solm\big\}_{t\in \Ic}$ on $\Xt$. 
The family $\big\{\Phi_{t_0,t}^x\big\}_{t\in \Ic}$, where $\Phi^x_{t_0,t}:\, \Xt\times\Omega\,{\rightarrow}\,\Xt$ is given~by 
\begin{equation}\label{cnt_fl}
\Phi_{t_0,t}^{x}(y,\om): = \Solm(x+y,\om) - \Solm(x,\om),  \qquad x,y\in \Xt,
\end{equation}
is called   the {\it centred two-point motion} associated with the stochastic flow $\big\{\Solm\big\}_{t\in \Ic}$ on $\Xt$. 
\end{definition}

\vspace{-0.2cm}
\begin{prop}\label{prp_centfl}
Consider the stochastic flow $\big\{\Solm)\big\}_{t\in \Ic}$ of diffeomorphisms. 
Then,  the centred two-point motion $\big\{\Phi_{t_0,t}^x\big\}_{t\in \Ic}$ is such that  
\begin{align}\label{Phiv_phixv}
 D\Phi^{x}_{t_0,t}(\om) = D\Solm(x,\om) \qquad  \p\,\textrm{\rm-\,a.s.} \;\;\forall \; x\in \Xt, t\in \Ic,
\end{align}
and is nonsingular. 
Moreover, $\big\{D\Phi^x_{t_0,t}\big\}_{t\in \Ic}$ is a stochastic flow satisfying~(\ref{var_dphi}), 
 and  it  induces the  transition  evolution   
\begin{align}\label{2pme}
\big(\mathcal{P}^{\Phi^x}_{t_0,t}g\big)(y) = \E\big[g\big(\varPhi^x_{t_0,t}(y)\big)\big], \qquad g\in \mathbb{M}(\Xt),  \quad t\in \Ic. 
\end{align}
If $\big\{\Solm\big\}_{t\in \Ic}$ is induced by the  SDE~(\ref{gen_sde}) with coefficients satisfying Assumption \ref{cf_remk},
then  the generator $\mathcal{L}^x_t$ of $\varPhi^{x}_{t_0,t}$ is given~by 
\begin{align}\label{2pp1}
\mathcal{L}^x_tg(y) &= {b}_i^x(t,y)\partial_{y_i}\,g(y)+{\textstyle \frac{1}{2}}\text{tr}\Big( (\sigma^x(t,y))^T\mathcal{H}_g(y)\sigma^x(t,y)\Big),
\end{align}
where  ${b}^x(t,y) = {b}(t,x+y)- {b}(t,x)$, $\sigma^x(t,y) = \sigma(t,x+y)- \sigma(t,x)$, and $\mathcal{H}_g(v)$ is the Hessian of~$g$.   
\end{prop}

\vspace{-0.1cm}
\noindent {\it Proof}. See Appendix~\ref{app_prp_centfl}.

\smallskip


\section{Finite-time expansion rates in stochastic flows}\label{s_exp_rates}

The new probabilistic measure of nonlinear expansion in stochastic flows is based on  what we term  $\varphi$-divergence rates ($\varphi$-DR).  This approach allows to consider finite-time expansion rates for both deterministic and stochastic flows in a unified framework without the need for linearisation of the underlying dynamics. Moreover, this framework is amenable to systematic extensions to flows on non-Euclidean manifolds which is a subject of ongoing work.  This probabilistic/information-theoretic approach is well-suited for deriving bounds in  
 the  Lagrangian Uncertainty Quantification (LUQ) \cite{branuda_luq}.  In addition, some specific connections of the $\varphi$-divergence rates  to some other descriptors of Lagrangian expansion are discussed in \S\ref{s_lyp_fun_stoch}.

\smallskip
Let $\varphi\!: \,\mathcal{J}\,{\subset}\,\R\,{\rightarrow}\, \R$  be a strictly convex, locally bounded function satisfying 
\begin{align}\label{Normality}
\varphi(1) = 0, \quad \varphi^{\prime}(1) = 0, \quad \inf_{a>0}\varphi(a)>-\infty,
\end{align}
and  let $\mu,\nu\in\PP(\Xt)$ be two probability measures on a measurable space $(\mathcal{M}, \Bb(\Xt))$. Then, the  $\varphi$-divergence $\Df\!: \,\PP(\Xt)\times\PP(\Xt) \rightarrow\Rp$ between $\mu$ and $\nu$  is defined by\footnote{\,Definition of $\Df$ in (\ref{dphi}) is related to that of $f$-divergence due to Csisz\'ar \cite{Csiszar72,csiszar91,Ciszar08} when $\mu\ll\nu$.  
However, the conditions~(\ref{Normality}) are often not imposed and $f$-divergences might not even be premetrics. 
Here, the constraint imposed on $\varphi$  removes the symmetries $\mathcal{D}_{f+c(u-1)} \,{=}\, \mathcal{D}_f$, $\mathcal{D}_{cf} \,{=}\, \mathcal{D}_f$, $c\ne0$, present in $f$-divergences. }

\begin{align}\label{dphi}
\Df(\mu\|\nu) = \begin{cases} \int_{\Xt}\varphi\left( \frac{d\mu}{d\gamma}\big/\frac{d\nu}{d\gamma}\right)d\nu, &\;\;\ \textrm{if }\; \mu,\nu\ll\gamma,\;\varphi\left( \frac{d\mu}{d\gamma}\big/\frac{d\nu}{d\gamma}\right)\in L^1(\Xt, \nu),\\
+\infty, &\;\;\; \textrm{otherwise},
\end{cases}
\end{align}
where $\gamma$ is any dominating  measure on $(\mathcal{M}, \Bb(\Xt))$.  Note that the definition in (\ref{dphi}) is independent of the choice of the dominating measure due to the properties of the Radon-Nikodym derivative.  
Clearly, $\Df$ in (\ref{dphi}) is generally not symmetric and it does not satisfy the triangle inequality. However,   $\Df$ is {\it information monotone} in the sense that 
  $\Df(\mu\|\nu)\geqslant \Df(\mu_\mathcal{A}\|\nu_\mathcal{A}),$
for $\mu_{\mathcal{A}}(B)= \mu(A_n\cap B)$, $\nu_{\mathcal{A}}(B) =\nu(A_n\cap B)$ for all $B\in \Bb(\Xt)$ and for any measurable partition $\mathcal{A} = \{A_n: n\in \N\}$ of $\Xt$. Information monotonicity is naturally imposed by physical constraints when {\it coarse-graining} the underlying dynamics and it  implies that $\Df$ is a {\it premetric}; i.e.,  $\Df(\mu\|\nu)\geqslant 0$ and $\Df(\mu\|\nu) =0$ iff $\mu = \nu$ almost everywhere. Moreover, $\Df$ is jointly convex and lower semi-continuous in its  arguments.
These properties follow readily  from  the following variational  representation (e.g.,~\cite{Gmeas})
\begin{align}\label{Dvar}
\Df(\mu\|\nu) &= \sup_{f\in \mathcal{C}_\infty(\Xt)}\bigg\{\int_{\Xt} f(x)\mu(dx) - \int_{\mathcal{M}}\varphi^*(f(x))\nu(dx)\bigg\}=:\sup_{f\in \mathcal{C}_\infty(\Xt)}\Big\{\big\langle f, \mu\big\rangle - \big\langle \varphi^*(f), \nu\big\rangle\Big\},
\end{align}
where $\varphi^*$ is the Legendre-Fenchel dual of $\varphi$.  Various well-known divergencies (some of them proper metrics) used in information theory, probability theory  and statistics are obtained from~(\ref{dphi}) with an appropriate choice of the convex function $\varphi$. 
In particular, the Kullback-Leibler divergence, $\mathcal{D}_{\textsc{kl}}$,  is obtained by setting $\varphi(u) = u\log u-u+1$ for $u>0$ in (\ref{dphi}) so that 
\begin{equation}\label{dkl}
\D_\textsc{kl}(\mu\|\nu) =\int_{\mathcal{M}}\log\left(\frac{d\mu}{d\lambda}\Big/\frac{d\nu}{d\lambda}\right)d\mu.
\end{equation}
 The variational representation of the   KL-divergence is given by  
\begin{equation}\label{KL_duality}
\DK(\mu\|\nu)= \sup_{f\in\mathcal{C}_\infty(\mathcal{M})}\Big(\big\langle f, \mu\big\rangle-\big\langle e^{f}, \nu\big\rangle+1\Big).
\end{equation}
Information-monotonicity of  $\varphi$-divergencies allows one to construct  (e.g., \cite{chentsov72}) a special Riemannian geometry on the manifold of probability measures in which a Pythagorean-like decomposition and (non-metric) geodesic projections are crucial in applications of {\it information-geometric} framework to statistical estimation and model selection (e.g., \cite{chentsov72,amari00,amari09,Amari10, Amari16}). The suitability of  the geometry induced by a given $\varphi$-divergence in applications depends on  the considered submanifold of probability measures (e.g., \cite{amari00,amari09,Amari16,Csiszar72,csiszar91,Ciszar08}).  Given that we aim to exploit these geometric properties in future work on uncertainty quantification in reduced-order models, we consider the whole family of $\varphi$-divergencies  in the framework developed in the subsequent sections\,\footnote{\,A number of other divergences and/or metrics, including Chernoff \cite{cher52}, Renyi \cite{ren61}, Bregman \cite{breg67} divergencies or  Wasserstein distance (e.g., \cite{Gmeas}), have been extensively used  in various contexts  including information theory, statistical inference, optimisation, image processing, neural networks; e.g.,~\cite{burnh02,niels88,Amari16, Vajda06, Amari10, Csiszar72, Ciszar08}). However, these divergencies and metrics are generally not information monotone and are thus not suitable for our future~purposes.}.

\subsection{Divergence rates ($\pmb{\varphi}\hspace{0.04cm}$-DR)}\label{ftdr_sec}
The Stroock-Varadhan martingale solutions and the support theorem \cite{Stroock79, Stroock72} provide a link between path-space interpretation of SDE's, stochastic flows  and time-marginal probability of their laws  (see  (\ref{k_lift}) and, e.g., \cite{Figali,Roc-Kry,leBris08}).  Thus, it is natural to quantify average trajectorial expansion rates  in stochastic flows via the associated time-marginal probability measures. Here, we define expansion rates in stochastic flows  as follows: 
\begin{definition}[\textbf{\textit{$\pmb{\varphi}$-divergence rate}}; $\pmb{\varphi}\hspace{0.04cm}$-{\bf DR}]\label{ftdr1}
\rm
Let $(\mu_t)_{t\in \Ic}$, $\Ic = [t_0,\,T)$, be a measurable family of probability measures in $\PP(\Xt)$ generated by  $\big(\mathcal{P}_{t_0,t}^{*}\big)_{t\in \Ic}$  so that $\mu_t = \mathcal{P}_{t_0,t}^*\mu_{t_0}$. The  {\it divergence rate}  between $\mu_t$ and $\mu_{t_0}$ for $t\in \Ic$ is given by
\begin{align}\label{Dphi}
\D^{t_0,t}_{\varphi}(\mu_t\|\mu_{t_0}):= \scaleobj{.85}{\frac{1}{ t-t_0}}\,\mathcal{D}_{\varphi}(\mu_t\|\mu_{t_0}), \qquad t\in \Ic.
\end{align} 
\end{definition}

\noindent We will show in the subsequent sections (see  \S\ref{s_lyp_fun_stoch}-\ref{s_comp}) that divergence rates can be used to  generalise the standard notions of expansion rates based (finite-time) Lyapunov exponents.
  For KL-divergence, the rate $\D^{t_0,t}_\textsc{kl}$ can be interpreted as the rate of loss of information in the measure $\mu_t$ relative to the initial measure $\mu_{t_0}$, thus providing a direct  information-theoretic characterisation of dynamic uncertainty in Lagrangian/path-based predictions due to errors in the initial conditions characterised by $\mu_{t_0}$.

\begin{rem}\label{rem_div_rts}\rm \mbox{}
\begin{itemize}[leftmargin=.75cm]
\item[(i)] Note that $\D^{t_0,t}_{\varphi}$ in (\ref{Dphi}) is well defined  any measurable family $(\mu_t)_{t\in \Ic}$ in $\PP(\Xt)$  as long as the absolute continuity conditions in (\ref{dphi}) are satisfied and $\varphi\in L^1(\Xt;\mu_{t_0})$. In particular, this includes the case of more general semimatringale flows and martingale solutions of SDE's with less regularity than required in Assumption \ref{cf_remk}. Additional conditions on the moments of $\mu_t$ are typically needed for $\varphi(d\mu_t/d\mu_{t_0}) \in L^1(\Xt;\mu_{t_0})$ as $T\rightarrow \infty$ (see, e.g., \cite{MBUda20}).

\newcommand{\Avg}{\mathcal{A}_{\varepsilon}}
 
\item[(ii)] Analogously to (\ref{Dphi}), one can define a (time) backward divergence rate between as 
\begin{align}\label{Dphi_bck}
\D^{T,T-u}_{\varphi}(\mu_{T-u}\|\mu_{T}):= \scaleobj{.85}{\frac{1}{u}}\,\mathcal{D}_{\varphi}(\mu_{T-u}\|\mu_{T}), \qquad u\in [0,\,T-t_0].
\end{align}
Most properties of the forward rates $\D^{t_0,t}_{\varphi}(\mu_t\|\mu_{t_0})$ discussed in the sequel can be shown to hold for the backward rates $\D^{T,T-u}_{\varphi}(\mu_{T-u}\|\mu_{T})$. However, the backward divergence rates do not seem to be important when considering the growth of uncertainty in Lagrangian prediction problems and we do not delve into the analytical treatment of this setup. We consider  some numerical examples in \S\ref{pics} as the backward rates might be of use in the analysis of the so-called Lagrangian transport which is tangential to our considerations.

\item[(iii)] 
  A measure of expansion/stretching  rates  known as the finite-time entropy (FTE) was considered and applied to transport in deterministic problems in \cite{froy12}; this measure  is defined~as 
 \begin{equation*}
\text{FTE}(\mu_{t_0}, \mu_t) := \lim_{\varepsilon\rightarrow 0}\frac{1}{\scriptstyle\vert t-t_0\vert}\Big[\mathcal{H}(\mathcal{A}_{\varepsilon}\mu_{t}) - \mathcal{H}(\mu_{t_0})\Big],\quad \mathcal{H}(\mu) := - \int_{\Rd}\log\left(\frac{d\mu}{dm_d}\right)\mu(dx),
\end{equation*}
where $\mathcal{H}$ is the  differential entropy  and  $\Avg$ is some averaging operator with a symmetric kernel.
FTE can be naturally extended to apply to  stochastic flows.  
However, it   requires an ad-hoc introduction of the operator $\mathcal{A}_\varepsilon$ in order to avoid a trivial behaviour for incompressible flows (i.e., flows that preserve the Lebesque measure) in which case $\mathcal{H}(\mu_{t})=\mathcal{H}(\mu_{t_0})$. 
The class of expansion rates defined in~(\ref{Dphi}) is devoid of such problems, and its properties are naturally suited for Lagrangian uncertainty quantification considered in \cite{branuda_luq}.  
 \end{itemize} 
 \end{rem}

\subsection{General properties of $\pmb{\varphi}$-divergence rates}\label{ss_exrates_prop}
The expansion rate  $\Df^{t_0,t}(\mu_t\|\mu_{t_0})$ introduced in (\ref{Dphi}) is well defined when $\mu_t, \mu_{t_0}\ll \gamma$  for all  $t\in \Ic$ and for some strictly positive dominating measure  $\gamma$ on $(\Xt, \Bb(\Xt))$; however, even if $\mu_{t_0}\ll \gamma$ the absolute continuity $\mu_t\ll \gamma$ is not guaranteed and it depends on  the underlying dynamics. 
If the family  $(\mu_t)_{t\in\Ic}$, $\mu_t\in \PP(\Xt)$, is  associated with a stochastic flow  generated by the SDE ~(\ref{gen_sde}) or (\ref{gen_sde_ito}), there is a wide range of  configurations  where the absolute continuity w.r.t.~the Lebesgue measure  $m_d$ on $\Xt,$ is automatically satisfied.

 \smallskip
 We begin the analysis with some general bounds on $\Df^{t_0,t}(\mu_t\| \mu_{t_0})$ and proceed to certain localised measures which are suitable for numerical approximations of $\varphi$-FTDR  in \S\ref{uni_exp_rates}. 

\begin{prop}\label{Dphi_posit}
 Assume that the coefficients of the SDE 
(\ref{gen_sde}) 
are sufficiently regular (e.g., they satisfy (\ref{reg_growth_cond})) so that the dynamics generates a global flow of $\mathcal{C}^{2}$-diffeomorphisms on $\Xt$, and $\rho_{t_0}{=}d\mu_{t_0}/dm_d\in \mathcal{C}^2_\infty(\Xt; \Rp)\cap L^1_+(\Xt;dx)$.
  Then    
\begin{align*}
\Df^{t_0,t}(\mu_t\| \mu_{t_0}) = \scaleobj{.9}{\frac{1}{|t-t_0|}}\int_{\Xt}\varphi\left(\frac{\rho_t(x)}{\rho_{t_0}(x)}\right)\rho_{t_0}(x)m_d(dx)<\infty. 
\end{align*}
\end{prop}
\noindent {\it Proof.} Given that $ \rho_{t_0}\in \mathcal{C}^2_\infty(\Xt; \Rp)\cap L^1_+(\Xt;dx)$, the same holds for $\rho_{t}$, $t\in \Ic$, by Theorem \ref{FKPt1}, and 
$\varphi(\rho_t/\rho_{t_0})$ is bounded  by Theorems \ref{mut_rnd}, Proposition \ref{al_pi} and Lemma \ref{alph_lem}.

 \newcommand{\BLf}{\Check{\LG}}
 \newcommand{\DP}{\alpha_{t_0}}
 
 \medskip

The following sequence of results provides  bounds on  $\Df^{t_0,t}(\mu_t\|\mu_{t_0})$ in terms of a Lyapunov exponent of the random density $\alpha_{t_0,t}$ in Proposition \ref{al_pi}, and then in terms of the coefficients of~(\ref{gen_sde}).  These bounds are reminiscent of bounds derived in \cite[\S5]{branuda_luq} and are utilised later.   

\begin{lem}[\cite{Kunitabook}, Lemma 4.3.4]\label{alph_lem} 
Assume that the conditions of Proposition \ref{Dphi_posit} hold and ${\rho_{t_0}\in \mathcal{C}^3_\infty(\Xt; \Rp)\cap L^1_+(\Xt;dx)}$. Then, the density $\alpha_{t_0,t} = d\Check\Pi_{t_0,t}/d\mu_{t_0}$, $t\in \Ic$, of the random  probability measure  (\ref{imm}) on $(\mathcal{M}, \Bb(\Xt))$ by the flow of solutions of (\ref{gen_sde}) can be  represented as  
\begin{align}\label{alph_Ga}
\alpha_{t_0,t}(x,\omega) &= \exp\bigg\{\int_{t_0}^tG\big(s,\phi_{t_0,s}(x,\omega)\big)dW_s-\frac{1}{2}\int_{t_0}^t a^{\scriptscriptstyle G}_s\big(\phi_{t_0,s}(x,\om), \phi_{t_0,s}(x,\om)\big)ds\bigg\}\notag\\[.3cm]
& \hspace{3cm} \times\exp\bigg\{\int_{t_0}^t   \frac{1}{\rho_{t_0}\big(\phi_{t_0,s}(x,\omega)\big)}\BLf_s^*\rho_{t_0}\big(\phi_{t_0,s}(x,\omega)\big)ds\bigg\},
 \end{align}
where $\rho_{t_0} = d\mu_{t_0}/dm_d >0$, and $\check{\mathcal{L}}_t^*$ is the $L^2(\Xt,\mu_{t})$ dual of  the generator
 \begin{equation}\label{bsdegen}
\check{\LG}_t f(x) =   -\sum_{i=1}^\ell \brr_i(t,x)\partial_{x_i} f(x)+\frac{1}{2} \sum_{i,j=1}^\ell a_{ij}(t,x)\partial_{x_ix_j} f(x),\quad f\in \mathcal{C}^2(\Xt),
\end{equation}
 of the backward flow of solutions to (\ref{gen_sde}), and  $a^{\scriptscriptstyle G}_t(x,x) {=} \langle G(\ccdot,x)\rangle_t$ 
is the quadratic variation
of  $G(t,x) \,{=}\, \rho_{t_0}(x)^{-1}\textrm{\rm div}\big(\rho_{t_0}(x)M(t,x)\big)$, where $M(t,x) =\int_{t_0}^t\sigma(s,x)dW_s.$ 
 \end{lem}
\begin{prop}\label{alpha_bnd1} Given the assumptions as in Proposition \ref{Dphi_posit}, the following holds  
\begin{align}\label{N121}
\Df^{t_0,t}(\mu_t\|\mu_{t_0}) \leqslant \scaleobj{.9}{\frac{1}{|t-t_0|}}\int_{\Xt}\E\big[\varphi_{\ddagger}(\alpha_{t_0,t}(x))\big]\mu_{t_0}(dx), \quad t\in \Ic,
\end{align}
where $\varphi_{\ddagger}(z) = z\varphi(z^{-1})\geqslant 0$ for all $z>0$, and $\E\big[\varphi_{\ddagger}(\alpha_{t_0,t}(x))\big]:=\int_\Om \varphi_{\ddagger}(\alpha_{t_0,t}(x,\om))\p(d\om)$.
\end{prop}
\noindent{\it Proof.} See Appendix~\ref{app_alpha_bnd1}.

\begin{theorem}\label{local-bound}
Suppose the drift and diffusion coefficients of (\ref{gen_sde}) satisfy Assumption \ref{cf_remk}.
If 
$\mu_{t_0}\ll m_d$ and  $\rho_{t_0} \in \mathcal{C}^3(\Xt, \Rp)\cap L^1(\Xt,\mathbb{R}^+)$, then
\begin{align*}
\Df(\mu_t\|\mu_{t_0})&\leqslant  \sum_{t_0\leqslant s\leqslant t}\E\big[\Delta\varphi_{\ddagger}(\alpha_{t_0,s})\big]+\int_{\Xt}\E\left[\int_{t_0}^tD^{-}\varphi_{\ddagger}(\alpha_{t_0,s})\alpha_{t_0,s}\frac{\check{\LG}_s^*\rho_{t_0}(\phi_{t_0,s})}{\rho_{t_0}(\phi_{t_0,s})}ds\right]d\mu_{t_0}\\[.2cm]
&\qquad +\frac{1}{2}\int_{\Xt}\E\bigg[\int_{(0,\infty)}L^\ell_t(\alpha)\varphi^{\prime\prime}_{\ddagger}(d\ell)\bigg] d\mu_{t_0}, \quad t,t_0\in \Ic.
\end{align*}
where $\mathscr{L}^{\ell}_t(\alpha)$ is the local time of $\alpha_{t_0,t}(x,\ccdot)$ at level $\ell\geqslant 0$. 
Moreover, if $\varphi_{\ddagger}\in \mathcal{C}^2\big((0,\infty),\R\big)$, then  
\begin{align*}
\Df(\mu_t\|\mu_{t_0})\leqslant \int_{\Xt}\E\left[\int_{t_0}^t\Big(\varphi^\prime_{\ddagger}(\alpha_{t_0,s})\alpha_{t_0,s}\frac{\check{\LG}_s^*\rho_{t_0}(\phi_{t_0,s})}{\rho_{t_0}(\phi_{t_0,s})}ds
 +\frac{1}{2}\varphi_{\ddagger}^{\prime\prime}(\alpha_{t_0,s})\alpha^2_{t_0,s}G^2(s,\phi_{t_0,s})ds\Big)\right]\mu_{t_0}(dx). \end{align*}
\end{theorem}

\vspace{-0.cm}
\noindent {\it Proof.} See Appendix \ref{app_local-bound}. 

\begin{cor}\label{localbound}
Let $\varphi(u) = u\log u -u+1,$ $u>0$, so that $\varphi_{\ddagger}(u) = -\log u +u-1$ 
 and  $\varphi^{\prime\prime}_{\ddagger}(u) = u^{-2}$. Then, for $t\in \Ic$
 \begin{align}\label{Local_FORM2}
\DK^{t_0,t}(\mu_t\|\mu_{t_0}) \leqslant  \int_{\Xt}\left\{ \scaleobj{.9}{\frac{1}{|t-t_0|}}\E\big[\alpha_{t_0,t}-1\big] -\E\big[\Lambda^\alpha_{t_0,t}\big]\right\} \mu_{t_0}(dx)ds,
\end{align}
where 
\begin{equation*}
\Lambda^\alpha_{t_0,t}(x,\om) =  \scaleobj{.9}{\frac{1}{|t-t_0|}} \log \big(\alpha_{t_0,t}(x,\om)\big)
\end{equation*}
 is the finite-time Lyapunov exponent for the random density $\alpha_{t_0,t}$ (see \cite{Kunitabook}). 
\end{cor}
\noindent {\it Proof.} This follows by direct computation utilising the specific form of $\varphi$ in $\varphi_\ddagger$, and using the fact that $\E\big[\alpha_{t_0,t}\big] = 1+\int_{t_0}^t\alpha_{t_0,s}\,\rho_{t_0}(\phi_{t_0,s})^{-1}\check{\LG}_s^*\rho_{t_0}(\phi_{t_0,s})ds$ (see (\ref{alph_dw})).

\begin{rem}\label{lapsum}\rm\mbox{}
{If the dynamics in (\ref{gen_sde}) is autonomous and  the transition semigroup $(\mathcal{P}_{t-t_0})_{t\in \Ic}$  admits a stationary ergodic probability measure $\bar \mu\in \PP(\Xt)$, a related bound on $\E\big[\D^{t_0,t}_\textsc{kl}(\Pi_{t_0,t}\|\bar\mu)\big]$  in terms of the sum  $\lambda_\Sigma$ of (infinite-time)  Lyapunov exponents was obtained in \cite[Theorems~4.2,~4.3]{Baxendale1}. We will return to it in \S\ref{Almost_sure}.  See, also  \cite{carverhill85,Crauel-89}} for an extensive discussion of Lyapunov exponents in stochastic flows in the ergodic setting.  
\end{rem}

\subsection{Fields of finite-time divergence rates ($\pmb{\varphi}$-DRF)}\label{uni_exp_rates}
Here, we establish the existence and continuity of $\varphi\hspace{0.04cm}$-DR fields  of locally averaged trajectorial expansion rates from a measure on the initial conditions 
 that is localised around some $x\in \Xt$. These results are important for the subsequent use of $\varphi\hspace{0.04cm}$-DR fields in the bounds of the type (\ref{bnd2}) in \cite{branuda_luq}, and for a systematic justification of computational approximations discussed in \S\ref{s_comp}. 
  In addition, 
 we employ these results are to identify connections of $\varphi\hspace{0.04cm}$-DRD's  to other well-known descriptors of Lagrangian expansion in deterministic flows, which are subsequently generalised  to stochastic~flows.

\begin{definition}[\textbf{\textit{Regularised uniform measure}}]\label{gmolf} \rm Given the uniform probability measure 
\begin{equation}
\mu_{B_{\varepsilon}(x)}(A) = \int_A\widehat{\I}_{B_{\varepsilon(x)}}(y)m_d(dy), \quad \widehat{\I}_{B_{\varepsilon(x)}}(y) = \frac{1}{m_d(B_\varepsilon(x))}\I_{B_{\varepsilon(x)}}(y)\qquad  A\in \mathcal{B}(\Xt), 
\end{equation}
with the probability density $\widehat{\I}_{B_{\varepsilon(x)}}$ supported on the ball $B_\varepsilon(x)$, $0\leqslant \varepsilon< \infty$, $x\in \Xt$, the measure 
\begin{equation}\label{reg_u_meas}
\tilde \mu^{x,\varepsilon}(A) = \int_A\tilde \rho^{x,\varepsilon}(y)m_d(dy), \qquad \tilde \rho^{x,\varepsilon}(y)=\eta^\varkappa\star\widehat{\I}_{B_{\varepsilon(x)}}(y) = \int_{\Rd} \eta^\varkappa(y-z)\widehat{\I}_{B_{\varepsilon(x)}}(z)dz,
\end{equation}
where 
\begin{equation}\label{mol}
\eta^\varkappa(\xi)= \big(2\pi\varkappa^2\big)^{-d/2}\exp\Big(-| \xi|^2/2\varkappa^2\Big), \quad 0<\varkappa<\infty,
\end{equation}
is referred to as the {\it Gaussian-regularised uniform measure} on $\Xt = \Rd$. Analogously,  for $\Xt= \bar{\mathbb{T}}^d$ (a flat $d$-dimensional torus) and $\bar B_\varepsilon(x)$, $0\leqslant \varepsilon$, $x\in \bar{\mathbb{T}}$, a  ball such that $\bar B_\varepsilon(x)\subseteq \bar{\mathbb{T}}$ the measure 
\begin{equation}
\tilde \mu^{x,\varepsilon}(A) = \int_A\tilde \rho^{x,\varepsilon}(y)m_d |_{\bar{\mathbb{T}}}(dy), \qquad \tilde \rho^{x,\varepsilon}(y)=\eta^\varkappa\star\widehat{\I}_{\bar B_{\varepsilon(x)}}(y) = \int_{\bar{\mathbb{T}}} \eta^\varkappa(y-z)\widehat{\I}_{\bar B_{\varepsilon(x)}}(z)dz,
\end{equation}
where   $0<\varkappa< \infty$ such that $\textrm{supp}(\,\widehat{\I}_{\bar B_\varkappa})\in \bar{\mathbb{T}} $ and 
\begin{equation}\label{molT}
\eta^\varkappa(\xi)= 
\mathcal{Z}_\varkappa^{-1}\Big(\varkappa+\I_{|\xi|\leqslant \varkappa}\exp\big(-1/(1- |\xi/\varkappa|^2)\big)\Big), \qquad \mathcal{Z}_\varkappa = \int_{\bar{\mathbb{T}}}\eta^\varkappa(\xi)d\xi, 
 \end{equation}
is referred to as the {\it regularised uniform measure} on $\Xt = \bar{\mathbb{T}}^d$.
\end{definition}

\begin{rem} \rm We skip the explicit dependence of $\tilde \mu^{x,\varepsilon}$ on $\varkappa$ in the subsequent derivations since all our results  hold for any 
  $0<\varkappa<\infty$.  The choice of the regulariser $\eta^\varkappa$ 
is not restricted to  (\ref{mol}) or (\ref{molT}) as long as $\eta^\varkappa>0$ is continuous,  $\int_\Xt\eta^\varkappa(\xi)d\xi = 1$ and $\lim_{\varkappa\rightarrow 0^+} \eta^\varkappa(x) = \delta(x)$ (in the distributional sense).  We focus on  $\varkappa\ll 1$ so that the  measure $\tilde\mu^{x,\varepsilon}$ is localised on the ball $B_\varepsilon(x)$ with the (Lebesgue) density $\tilde \rho^{x,\varepsilon}>0$ but such that it is approximately uniform on $B_\varepsilon(x)$ to any computational accuracy. 
\end{rem}

Now, consider the  measurable family of time-marginal probability measures $(\tilde \mu^{x,\varepsilon}_t)_{t\in \Ic}$, 
$x\in \Xt$, such that  $\tilde \mu^{x,\varepsilon}_t = \mathcal{P}^{ \Phi^x*}_{t_0,t}\tilde \mu^{x,\varepsilon}_{t_0}$, $\tilde \mu^{x,\varepsilon}_{t_0} \ll m_d$, and  \begin{equation}\label{Mollify}
\tilde \mu_t^{x,\varepsilon}(dy) =\tilde \rho^{x,\varepsilon}_{t}(y)m_d(dy),\quad  \tilde \rho_t^{x,\varepsilon}(y) = \mathcal{P}^{\Phi^{x*}}_{t_0,t}\int_{\Xt}\widehat\I_{B_\varepsilon(x)}(y-z)\eta^{\varkappa}(z)m_d(dz), \;\; 0<\varkappa<\infty,
\end{equation}
 where $\mathcal{P}^{ \Phi^x*}_{t_0,t}$ is the  $L^1(\Xt,\mu_{t_0})$ dual of the transition evolution $\mathcal{P}^{\Phi^x}_{t_0,t}$ in~(\ref{2pme}) 
  and $\big\{\Phi^x_{t_0,t}\big\}_{t\in\Ic}$  is the centred two-point motion  (\ref{cnt_fl}) induced by the stochastic flow $\big\{\Solm\big\}_{t\in \Ic}$ on $\Xt$.
 In particular, the above holds if the family of measures $(\tilde \mu^{x,\varepsilon}_t)_{t\in \Ic}$ is associated with  a stochastic flow generated by an ODE/SDE with sufficiently regular coefficients (e.g., such as those in (\ref{reg_growth_cond})); see~Theorem~\ref{FKPt1}.

\begin{definition}[\textbf{\textit{$\pmb{\varphi}$-divergence rate field; $\pmb{\varphi}$-DRF}}]\label{l_ftdr}\rm Let $(\tilde \mu^{x,\varepsilon}_t)_{t\in \Ic}$ be a measurable family of Borel probability measures on $\Xt$. Denote the $\varphi$-divergence rate from the initial measure $\tilde\mu_{t_0}^{x,\varepsilon}\in \PP(\Xt)$ localised at $x\in \Xt$ 
 by 
\begin{equation}\label{xlocrt}
 \mathcal{R}^{\Phi}_{\varphi,\varepsilon}(x,t_0,t) :=\D_\varphi^{t_0,t}\big(\mathcal{P}^{\Phi^x*}_{t_0,t}\tilde\mu_{t_0}^{x,\varepsilon}\|\tilde\mu_{t_0}^{x,\varepsilon}\big), \qquad t\in \Ic=[t_0,\,T).
\end{equation}
The map $x\mapsto  \mathcal{R}^{\Phi}_{\varphi,\varepsilon}(x,t_0,t)$ is referred to as the (time) forward $\varphi$-divergence rate field ($\varphi$-DRF) over the interval $[t_0,\,t]\subseteq\Ic$. 
\end{definition}

\begin{theorem}\label{exp_rate}  Let $\big\{\Solm\big\}_{t\in \Ic}$ be a  stochastic flow of diffeomorphisms on $\Xt$ with the associated  time-marginal probability measures such that $\D^{t_0,t}_\varphi\big(\mathcal{P}^{\Phi^x*}_{t_0,t}\tilde\mu_{t_0}^{x,\varepsilon}\|\tilde\mu_{t_0}^{x,\varepsilon}\big)<\infty$ for all $t\in \Ic$,  $x\in \Xt$.
Then, for any $0\leqslant \varepsilon< \infty$ and all $x\in \Xt$ such that $B_\varepsilon(x)\subset \Xt$ the following hold:

\begin{itemize}[leftmargin=0.9cm]
\item[{\rm (i)}] $ t\mapsto \RL_{\varphi,\varepsilon}^{\Phi}(x,t_0,t)$ is continuous  for all $x\in \Xt$,
\item[{\rm (ii)}] $ x\mapsto \RL_{\varphi,\varepsilon}^{\Phi}(x,t_0,t)$ is continuous for all $t\in \Ic$,

\vspace{.1cm}\item[{\rm (iii)}] 
$\lim_{\varepsilon\rightarrow 0^+}\RL_{\varphi,\varepsilon}^{\Phi}(x,t_0,t) = \RL^{D\Phi}_\varphi(x,t_0,t)$.
\end{itemize}
 \end{theorem}
\begin{rem}\rm The following comments are in order:
\begin{itemize}[leftmargin=0.5cm]
\item[(i)]  The $\varphi$-divergence rate fields  $\mathcal{R}^{\Phi}_{\varphi,\varepsilon}$ are defined for a general  stochastic flow.  
 Thus, we  study  $\varphi$-DRF's  under  weaker assumptions that are 
automatically satisfied for the dynamics in~(\ref{gen_sde}) when $(\brr,\,\sigma)$ 
 satisfy~(\ref{reg_growth_cond}); see Proposition~\ref{Dphi_posit}. 
In what follows,  we consider initial probability measures with smooth Lebesgue densities but  the results can be generalised to  weak solutions of (\ref{rho_FP}) as long as $\tilde \rho_t^{x,\varepsilon}>0$ for all $t\in \Ic$.
\item[(ii)] It can be shown that analogous statements to those in Theorem \ref{exp_rate} apply to time backward $\varphi$-divergence fields defined as 
\begin{equation}\label{xlocrt_bck}
 \mathcal{R}^{\Phi}_{\varphi,\varepsilon}(x,T,T-t) :=\D_\varphi^{t_0,t}\big(\mathcal{P}^{\Phi^x*}_{t_0,T-t}\tilde\mu_{t_0}^{x,\varepsilon}\|\mathcal{P}^{\Phi^x*}_{t_0,T}\tilde\mu_{t_0}^{x,\varepsilon}\big), \qquad t\in [0,\,T-t_0].
\end{equation}
We do not elaborate  on the backward $\varphi$-divergence rate fields since they do not seem important for uncertainty quantification in (forward) Lagrangian predictions (see also Remark \ref{rem_div_rts}(ii)). Some numerical examples are shown in \S\ref{pics}, since these fields are related to fields of time-backward Lyapunov exponents which are used in Lagrangian transport analysis. 
\end{itemize}

\end{rem}

 The proof of Theorem \ref{exp_rate} 
is on the tedious side except  for part (i) which  follows  from the strict convexity of $\varphi$ and time continuity of the centred two-point motion $\{\Phi^x_{t_0,t}\}_{t\in \Ic}$, which implies  time continuity of the transition evolutions $(\mathcal{P}_{t_0,t}^{\Phi^x})_{t\in\Ic}$ in (\ref{2pme}).
The rest of the proof 
is split into two parts  given in Propositions~\ref{Con} and {\ref{prt_iii}}, corresponding to parts (ii)--(iii) of Theorem \ref{exp_rate}.

\begin{prop}[Property (ii) in Theorem \ref{exp_rate}]\label{Con}
Let $\big\{\Solm\big\}_{t\in \Ic}$ be a  stochastic flow on $\Xt$ such that  $\D^{t_0,t}_\varphi\big(\mathcal{P}^{\Phi^x*}_{t_0,t}\tilde\mu_{t_0}^{x,\varepsilon}\|\tilde\mu_{t_0}^{x,\varepsilon}\big)<\infty$ for all $t\in \Ic$, $x\in \Xt$. Then,  the map $x\mapsto \RL_{\varphi,\varepsilon}^{\Phi}(x,t_0,t) $  is continuous for any $0\leqslant \varepsilon<\infty$  such that $B_\varepsilon(x)\subset \Xt$.
\end{prop}
\newcommand{\Mdp}{\rho_{t_0}^{\varepsilon, \alpha}}
\noindent {\it Proof.} 
Given the definition of $\RL_{\varphi,\varepsilon}^{\Phi}$ 
it is sufficient to prove the claim for the $\varphi$-divergence $\D_\varphi$. 
 
Consider $\tilde \mu_{t_0}^{x,\varepsilon}\,{\in}\, \PP(\Xt)$  in (\ref{Mollify}) 
  and let $(x_n)_{n\ge 1}$ be a sequence in $\Xt$ such that $x_n\rightarrow x\,{\in}\, \Xt$ as $n\rightarrow\infty.$ Given the assumptions 
$\D_\varphi\big(\mathcal{P}_{t_0,t}^{\Phi^x*}\tilde\mu_{t_0}^{x,\varepsilon}\|\tilde \mu_{t_0}^{x,\varepsilon}\big)$ and $\D_\varphi\big(\mathcal{P}^{\Phi^{x_n}*}_{t_0,t} \tilde\mu_{t_0}^{x,\varepsilon}\|\tilde\mu_{t_0}^{x,\varepsilon}\big)$ 
are well defined and finite for all $t_0,t\in \Ic.$ First, we  show that $x\mapsto \D_\varphi(\mathcal{P}^{\Phi^x*}_{t_0,t}\tilde\mu_{t_0}^{x,\varepsilon}\|\tilde \mu_{t_0}^{x,\varepsilon})$ is lower semicontinuous for all $t_0,t\in \Ic$. This property can be derived from the variational  representation of $\D_\varphi$ which yields 
 \begin{align*}
 \D_\varphi(\mathcal{P}^{\Phi^x*}_{t_0,t}\tilde\mu_{t_0}^{x,\varepsilon}\|\tilde\mu_{t_0}^{x,\varepsilon}) &= \sup_{h\in \mathcal{C}_\infty(\Xt)}\bigg\{\langle \mathcal{P}^{\Phi^x}_{t_0,t}h, \tilde\mu_{t_0}^{x,\varepsilon}\rangle -\langle \varphi^*(h), \tilde\mu_{t_0}^{x,\varepsilon}\rangle\bigg\}\\[.2cm]
 &\hspace{1cm}\leqslant \D_\varphi\big(\mathcal{P}^{\Phi^{x_n}*}_{t_0,t}\tilde\mu_{t_0}^{x,\varepsilon}\|\tilde\mu_{t_0}^{x,\varepsilon}\big)+\sup_{h\in\mathcal{C}_\infty(\Xt)}\bigg\{\big\langle\mathcal{P}^{\Phi^x}_{t_0,t}h, \tilde\mu_{t_0}^{x,\varepsilon}\big\rangle - \big\langle \mathcal{P}_{t_0,t}^{\Phi^{x_n}}h, \tilde\mu_{t_0}^{x,\varepsilon}\big\rangle\bigg\}.
 \end{align*}
 Note that $ \Phi_{t_0,t}^{x}(v,\om)= \Solm(x+v,\om) - \Solm(x,\om)$ and $\mathcal{P}^{\Phi^x}_{t_0,t}h(v) = \E\big[h\big(\Solm(x+v)-\Solm(x)\big)\big]$ for $h\in \mathcal{C}_\infty(\Xt)$ so that  
\begin{align*}
&\big\langle\mathcal{P}^{\Phi^x}_{t_0,t}h, \tilde\mu_{t_0}^{x,\varepsilon}\big\rangle - \big\langle \mathcal{P}_{t_0,t}^{\Phi^{x_n}}h, \tilde\mu_{t_0}^{x,\varepsilon}\big\rangle= \int_{\Xt}\E\big[h(\Solm(x+v)-\Solm(x))\big]\tilde\mu_{t_0}^{x,\varepsilon}(dv)\\[.2cm]
&\hspace{6cm}-\int_{\Xt}\E\big[h(\Solm(x_n+v)-\Solm(x_n))\big]\tilde\mu_{t_0}^{x,\varepsilon}(dv),
\end{align*}
which converges to zero wp1 (with probability 1) as $n\rightarrow \infty$  due to continuity of  $x\mapsto\Solm(x,\ccdot)$  for all $t_0,t\in \Ic$. Then,  by Fatou's lemma, we have the required lower semicontinuity
\begin{align}\label{.521}
\D_\varphi\big(\mathcal{P}^{\Phi^x*}_{t_0,t}\tilde\mu_{t_0}^{x,\varepsilon}\|\tilde\mu_{t_0}^{x,\varepsilon}\big)\leqslant \liminf_{n\rightarrow \infty}\D_\varphi\big(\mathcal{P}^{\Phi^{x_n}*}_{t_0,t}\tilde\mu_{t_0}^{x,\varepsilon}\|\tilde\mu_{t_0}^{x,\varepsilon}\big) \qquad \forall \;t,t_0\in \Ic. 
\end{align}
\smallskip
In order to  show that $x\mapsto \D_\varphi\big(\mathcal{P}^{\Phi^x*}_{t_0,t}\tilde\mu_{t_0}^{x,\varepsilon}\|\tilde\mu_{t_0}^{x,\varepsilon}\big)$ is also upper semicontinuous, consider the function 
\begin{align*}
{H}_{\varphi}: \Xt\rightarrow [0,\infty), \quad x\mapsto H_{\varphi}(x):= \D_\varphi\big(\mathcal{P}^{\Phi^x*}_{t_0,t}\tilde\mu_{t_0}^{x,\varepsilon}\|\tilde\mu_{t_0}^{x,\varepsilon}\big).
\end{align*}
By the properties of the $\varphi$-divergence and the assumptions of the proposition,  $x\mapsto H_{\varphi}(x)$ is bounded from below and from above (i.e., $0 \leqslant H_{\varphi}(x)<\infty$). 
 Now, 
  take $k\in\N$ such that  $k|x-y|>H_{\varphi}(y) - H_{\varphi}(x)$. Then,  we have 
 \begin{equation}\label{.54}
H_{\varphi}(x)+k\vert x-y\vert >H_{\varphi}(y).
\end{equation}
   Thus, for any sequence $(x_n)_{n\geqslant 1}\subset \Xt$  we have  from (\ref{.54}) that 
\begin{align*}
H_{\varphi}(x)\geqslant  H_{\varphi}(x_n)-k\vert x-x_n\vert \quad  \Longrightarrow \quad H_{\varphi}(x_n)\leqslant H_{\varphi}(x) \quad \forall \;n\in\N,
\end{align*}
 and consequently
\begin{align}\label{.55}
\limsup_{k\rightarrow\infty} H_{\varphi}(x_k)\leqslant H_{\varphi}(x).
\end{align}
Finally, comparing (\ref{.521}) and (\ref{.55}), we conclude that 
\begin{align*}
\lim_{k\rightarrow\infty}H_{\varphi}(x_k) = H_{\varphi}(x),
\end{align*}
which implies that $x\mapsto \D_\varphi\big(\mathcal{P}^{\Phi^x*}\tilde\mu_{t_0}^{x,\varepsilon}\|\tilde\mu_{t_0}^{x,\varepsilon}\big)$ is continuous. \qed

\newcommand{\BTb}{\hat{B}_1^x(0)}
\begin{prop}[Property (iii) in Theorem \ref{exp_rate}]\label{prt_iii} 
Let $\big\{\Solm\big\}_{t\in \Ic}$ be a  stochastic flow of diffeomorphisms on $\Xt$  such that $\D_\varphi\big(\mathcal{P}^{\Phi^x*}_{t_0,t}\tilde\mu_{t_0}^{x,\varepsilon}\|\tilde\mu_{t_0}^{x,\varepsilon}\big)<\infty$ for all $t\in \Ic$,  $x\in \Xt$.
Then, the limit of $\mathcal{R}^{\Phi}_{\varphi,\varepsilon}(x;t_0,t)$ as $\varepsilon\downarrow 0$ exists and it coincides with  $\mathcal{R}^{D\phi}_{\varphi}(x,t_0,t)$.
\end{prop}
\noindent {\it Proof.} 
   Given the variational representation of $\Df$ (\ref{Dvar})  and Proposition~\ref{prp_centfl}, we have
\begin{align}\label{T_cp1}
\notag  \D_\varphi\big(\mathcal{P}^{\Phi^x*}_{t_0,t}\tilde\mu_{t_0}^{x,\varepsilon}\|\tilde\mu_{t_0}^{x,\varepsilon}\big) & =  \sup_{h\in\mathcal{C}_\infty(\Xt)}\bigg\{\int_{\Xt}(\mathcal{P}^{\Phi^x}_{t_0,t}h)(y)\tilde\mu_{t_0}^{x,\varepsilon}(dy)-\int_{\Xt}\varphi^*(h(y))\tilde\mu_{t_0}^{x,\varepsilon}(dy)\bigg\}\\[.2cm]
\notag &\leqslant \sup_{h\in \mathcal{C}_\infty(\Xt)}\bigg\{\int_{\Xt}(\mathcal{P}^{\Phi^x}_{t_0,t}h)(y)\tilde\mu_{t_0}^{x,\varepsilon}(dy)-\int_{\Xt}(\mathcal{P}^{D\phi}_{t_0,t}h)(y)\tilde\mu_{t_0}^{x,\varepsilon}(dy)\bigg\}\\[.2cm]
&\hspace{1.5cm} + \Df\big(\mathcal{P}_{t_0,t}^{D\phi*}\tilde\mu_{t_0}^{x,\varepsilon}\|\tilde\mu_{t_0}^{x,\varepsilon}\big),
\end{align}
where $\mathcal{P}_{t_0,t}^{D\phi}$ is the transition evolution of the derivative flow $D\Phi^x_{t_0,t} = D\Solm(x)$ and $\mathcal{P}_{t_0,t}^{D\phi*}$ is its dual. Given that we consider $\Xt = \Rd$ or $\Xt = \bar{\mathbb{T}}$ (so that the tangent bundle $T\Xt \simeq\Xt\times\Xt$), we do not explicitly operate on tangent spaces and tangent measures.
On the other hand, 
\begin{align}\label{T_cp2}
\notag \Df\big(\mathcal{P}_{t_0,t}^{D\phi*}\tilde\mu_{t_0}^{x,\varepsilon}\|\tilde\mu_{t_0}^{x,\varepsilon}\big) &\leqslant\sup_{h\in \mathcal{C}_\infty(\Xt)}\bigg\{ \int_{\Xt}(\mathcal{P}_{t_0,t}^{D\phi}h)(y)\tilde\mu_{t_0}^{x,\varepsilon}(dy)-\int_{\Xt}(\mathcal{P}_{t_0,t}^{\Phi^{x}}h)(y)\tilde\mu_{t_0}^{x,\varepsilon}(dy)\bigg\}\\[.2cm]
&\hspace{1cm} + \Df\big(\mathcal{P}^{\Phi^{x}*}_{t_0,t}\tilde\mu_{t_0}^{x,\varepsilon}\|\tilde\mu_{t_0}^{x,\varepsilon}\big).
\end{align}
Comparing (\ref{T_cp1}) and (\ref{T_cp2}), it is sufficient to  show that for all $ t\in \Ic$
\begin{align*}
\left\vert\int_{\Xt}(\mathcal{P}^{\Phi^x}_{t_0,t}h)(y)\tilde\mu_{t_0}^{x,\varepsilon}(dy)-\int_{\Xt}(\mathcal{P}^{D\phi}_{t_0,t}h)(y)\tilde\mu_{t_0}^{x,\varepsilon}(dy)\right\vert \underset{\varepsilon\rightarrow 0^+}\longrightarrow 0 \qquad \forall\,h\in \mathcal{C}_\infty(\Xt).
\end{align*}
First, note that $\int_{\Xt}(\mathcal{P}^{\Phi^x}_{t_0,t}h)(y)\tilde\mu_{t_0}^{x,\varepsilon}(dy)<\infty$ (by the properties of the transition evolution), and given the assumptions of the proposition, we have 
\begin{align}\label{prt1}
\int_{\Xt}(\mathcal{P}^{\Phi^x}_{t_0,t}h)(y)\tilde\mu_{t_0}^{x,\varepsilon}(dy) &= \int_{\Xt}(\mathcal{P}^{\Phi^x}_{t_0,t}h)(y)\tilde\rho_{t_0}^{x,\varepsilon}(y)m_d(dy) \notag\\[.1cm]
=& \int_{\Xt}\int_{\Xt}(\mathcal{P}^{\Phi^x}_{t_0,t}h)(y)\widehat\I_{B_\varepsilon(x)}(y-z)\eta^{\varkappa}(z)m_d(dz)m_d(dy)\notag\\[.2cm]
\overset{y=\varepsilon \zeta}{=}& \int_{\Xt}\int_{\Xt}(\mathcal{P}^{\Phi^x}_{t_0,t}h)(\varepsilon\zeta)\I_{B_\varepsilon(x)}(\varepsilon\zeta-z)\eta^{\varkappa}(z)m_d(dz)m_d(d\zeta).
\end{align}
Analogously, 
\begin{align}\label{prt2}
\int_{\Xt}(\mathcal{P}^{D\phi}_{t_0,t}h)(y)\tilde\mu_{t_0}^{x,\varepsilon}(dy)=  \int_{\Xt}\int_{\Xt}(\mathcal{P}^{D\phi}_{t_0,t}h)(\varepsilon\zeta)\I_{B_\varepsilon(x)}(\varepsilon\zeta-z)\,\eta^{\varkappa}(z)m_d(d\zeta)m_d(dz)<\infty.
\end{align} 
Combining (\ref{prt1}) and (\ref{prt2})
  leads to 
\begin{align}\label{ebnd}
&\left\vert\int_{\Xt}(\mathcal{P}^{\Phi^x}_{t_0,t}h)(y)\tilde\mu_{t_0}^{x,\varepsilon}(dy)-\int_{\Xt}(\mathcal{P}^{D\phi}_{t_0,t}h)(y)\tilde\mu_{t_0}^{x,\varepsilon}(dy)\right\vert \notag\\[.2cm]
&\hspace{2cm}= \left\vert \int_{\Xt}\int_{\Xt} \left( (\mathcal{P}^{\Phi^x}_{t_0,t}h)(\varepsilon\zeta) -(\mathcal{P}^{D\phi}_{t_0,t}h)(\varepsilon\zeta)  \right)\I_{B_\varepsilon(x)}(\varepsilon\zeta-z)\,\eta^{\varkappa}(z)m_d(d\zeta)m_d(dz)\right\vert \notag\\[.2cm]
&\hspace{2cm}\leqslant   \int_{\Xt}\int_{\Xt} \left|(\mathcal{P}^{\Phi^x}_{t_0,t}h)(\varepsilon\zeta) -(\mathcal{P}^{D\phi}_{t_0,t}h)(\varepsilon\zeta)  \right|\I_{B_\varepsilon(x)}(\varepsilon\zeta-z)\,\eta^{\varkappa}(z)m_d(d\xi)m_d(dz)\notag\\[.2cm]
&\hspace{2cm}\leqslant   \underset{\zeta\in \Xt}{\textrm{sup}} \left|(\mathcal{P}^{\Phi^x}_{t_0,t}h)(\varepsilon\zeta) -(\mathcal{P}^{D\phi}_{t_0,t}h)(\varepsilon\zeta)  \right|.
\end{align}
Next, we have  for any $t\in \Ic = [t_0,\,t_0+T)$
\begin{align}\label{T_cp3}
\notag \left\vert (\mathcal{P}_{t_0,t}^{\Phi^{x}}h)(\varepsilon \zeta)-(\mathcal{P}_{t_0,t}^{D\phi}h)(\varepsilon \zeta)\right\vert &\leqslant \|h\|_\infty\,\E\left \vert \Phi^{x}_{t_0,t}(\varepsilon \zeta) - D\phi_{t_0,t}(x)\varepsilon \zeta\right\vert \\[.2cm]
& \leqslant \|h\|_\infty\,\E\left \vert \phi_{t_0,t}(x+\varepsilon \zeta) -\phi_{t_0,t}(x)- D\phi_{t_0,t}(x)\varepsilon \zeta\right\vert\notag  \\
& \leqslant \mathfrak{C}\|h\|_\infty |\mathcal{O}(\varepsilon^2)|,
\end{align}
where $\mathfrak{C} = \mathfrak{C}(x,\xi,\Ic)$. Finally, combining (\ref{T_cp3}) and (\ref{ebnd}), 
 we have based on the dominated convergence theorem that for $t\in \Ic$
\begin{align*}
\left\vert\int_{\Xt}(\mathcal{P}^{\Phi^x}_{t_0,t}h)(y)\tilde\mu_{t_0}^{x,\varepsilon}(dy)-\int_{\Xt}(\mathcal{P}^{D\phi}_{t_0,t}h)(y)\tilde\mu_{t_0}^{x,\varepsilon}(dy)\right\vert \underset{\varepsilon\rightarrow 0^+}\longrightarrow 0 \qquad \forall \,h\in \mathcal{C}_\infty(\Xt),
\end{align*}
which completes the proof in light of (\ref{T_cp1})--(\ref{T_cp2}).\qed


\section{Divergence based expansion rates and Lyapunov functionals on a finite-time horizon for stochastic flows}\label{s_lyp_fun_stoch}

In this section, we elucidate connections between the $\varphi$-divergence rates ($\varphi$-DR) introduced in~\S\ref{s_exp_rates} and some Lyapunov functionals for stochastic flows. The primary focus  of this work lies in considering the properties of $\varphi$-DR fields and their  subsequent applications in  Lagrangian uncertainty quantification through the bounds  (\ref{bnd2})  and (\ref{bnd1}) derived in  \cite{branuda_luq}. However, it turns out that a restriction of the general bound (\ref{bnd1}) provides interesting links between $\varphi$-DR fields (\ref{xlocrt}) based on the KL-divergence (\ref{dkl}) and various  Lyapunov functionals for both  observables and  probability measures.
{Such  functionals play an important role in studies ranging from stochastic stability, multiplicative ergodic theory, large deviations, etc.~(see, e.g., \cite{Arnold1,Arnold86,Kunitabook,Baxendale1,Baxendale2,carverhill85,Crauel-89}). 
Here,  as an prelude to future work devoted to the general treatment of such functionals, we focus on finite-time Lyapunov exponents (FTLE) which are frequently used in applications (predominantly in the deterministic setting) to estimate expansion rates in fluid flows and, to some extent,  serve as a proxy for flow-invariant manifolds which represent   barriers to trajectory-based/Lagrangian transport (see, e.g., \cite{shadden,froy09, Haller_BV12,haller01,hadji17} subject to several caveats). 
 In many deterministic flows the time-forward/time-backward FTLE fields tend to align (though not provably) with flow-invariant stable/unstable manifolds of hyperbolic trajectories
 (see, e.g., \cite{branwig10}}). In stochastic dynamical systems with small-amplitude noise a different approach for identifying (diffusive) barriers to Lagrangian transport was recently developed in \cite{haller18, haller20}. 
Although such considerations are tangential to our study (we do not look for transport barriers),  the abundance and popularity of approaches to characterisation of expansion rates and transport barriers via Lyapunov exponents merits a systematic  outline of  general links between FTLE and $\varphi$-DR fields. This link is particularly relevant since the $\varphi$-DR fields are based on the general nonlinear flow (as opposed to the inherent linearisation in FTLE fields), and they  have a clear probabilistic/information-theoretic interpretation in both the deterministic and a fully stochastic setting with no restriction on non-degeneracy of the diffusion term or its norm/amplitude. 
 
\subsection{KL-divergence rates and the largest Lyapunov exponent}\label{Almost_sure}\mbox{}
 Given a stochastic flow $\{\Solm\}_{t\in \Ic}$ of $\mathcal{C}^1$-diffeomorphisms on $\Xt$, we first establish  two bounds on the empirical approximation of the largest finite-time Lyapunov exponent 
$\bar{\Lambda}^{t-t_0}_{t_0}(x)$ at $x\in \Xt$  over $[t_0, \, t]\subseteq \Ic$ in terms of $\varphi$-DR fields (\ref{xlocrt}) based on the KL-divergence (\ref{dkl}) in the form (see Proposition \ref{KL_max_Lyap_stoch})
\vspace{.2cm}
\begin{align}\label{}
-\mathfrak{C}_-^{t_0,t}- \mathcal{R}^{D\Phi}_{\textsc{kl}}(x,t_0,t)\leqslant  \E^{\tilde{\mu}^x_{t_0}}\big[\bar{\Lambda}^{t-t_0}_{t_0}(x)\big]\leqslant \mathcal{R}^{D\Phi}_{\textsc{kl}}(x,t_0,t) +\mathfrak{C}_+^{t-t_0}\quad \forall \;t\in \Ic,
\end{align}
where the measure $\tilde\mu^x_{t_0}\in \PP(\Xt)$ on the initial perturbation is localised at $x\in \Xt$, and the nonnegative constants $\mathfrak{C}_{\pm}^{t_0,t}(\tilde\mu^{x}_{t_0})$ are such that $\mathfrak{C}_\pm^{t_0,t}\rightarrow 0$ as $\tilde\mu^{x}_{t_0}\overset{\ast}{\rightharpoonup} \delta_x$. In addition, we show that  the KL-divergence rate fields $\mathcal{R}^{D\Phi}_{\textsc{kl}}(x,t_0,t)$ are bounded from above by expected divergence rates $\E\!\big[\DK^{t_0,t}\big({\Pi}^x_{t_0,t}\|\tilde{\mu}^x_{t_0}\big)\big]$ based on random measures ${\Pi}^x_{t_0,t}$ (\ref{Random measure1}) carried by the derivative flow $\{D\phi_{t_0,t}\}_{t\in \Ic}$.  This bound is useful when considering links between FTLE's and empirically averaged divergence rate fields estimated from multiple experiments. Finally, following on the results established in \cite{branuda_luq}, we derive a tight bound  
\begin{align}\label{bnd11}
\hat{\mathcal{K}}_{\varphi}^{\mu}\big(-\mathcal{D}_\textsc{kl}(\tilde\mu^x_t\|\tilde\mu^x_{t_0})\big)\leqslant  \E^{\tilde{\mu}^x_{t_0}}\big[\bar{\Lambda}^{t-t_0}_{t_0}\big] \leqslant \mathcal{K}_{\varphi}^{\mu}\big(\mathcal{D}_\textsc{kl}(\tilde\mu^x_t\|\tilde\mu^x_{t_0})\big), \quad t\in \Ic, 
\end{align}
  where $\mathcal{K}_{\varphi}^{\mu}(u)\rightarrow 0$, $\hat{\mathcal{K}}_{\varphi}^{\mu}(-u)\rightarrow 0$ as $u\,{\downarrow}\, 0$.

\medskip
We start by considering the evolution of a perturbation $y\in \Xt$ given by 
\begin{align*}
\mathcal{Y}_t^{t_0,y}(x,\om) &:= 
\Phi^x_{t_0,t}(y,\om)=\phi_{t_0,t}(x+y,\om)-\phi_{t_0,t}(x,\om),\quad \, t\in \Ic, 
\end{align*} 
and its  tangent approximation (see Proposition \ref{prp_centfl})
\begin{equation}
Y_t^{t_0,y}(x,\om) = D\Phi_{t_0,t}^x(\om)y,\quad \, t\in \Ic,
\end{equation}
 so that 
\begin{align}\label{yt}
\big| Y^{t_0,y}_t(x,\om)\big|& = \left(y^TM(t_0,t,x,\om)\,y\right)^{1/2}, \quad M(t_0,t,x,\om)= \big(D\Phi^x_{t_0,t}(\om)\big)^*D\Phi_{t_0,t}^x(\om).
\end{align}

\begin{definition}[\textbf{\textit{Finite-time Lyapunov exponents}}]\label{stoch_lyap} \rm\; 
Consider the centred two-point motion $\{\Phi^x_{t_0,t}\}_{t\in\Ic}$ on $\Xt$  given by (\ref{cnt_fl}) and generated by a stochastic flow of $\mathcal{C}^1$-diffeomorphisms. 
\begin{itemize}[leftmargin = 0.7cm] 
\item[(a)]  
 The {\it maximal stochastic time-forward Lyapunov exponent} at $x\in\Xt$ over $[t_0,\,t]\subseteq\Ic$ associated with a realisation of the stochastic flow  $x\mapsto\phi_{t_0,t}(x,\om)$  is defined for almost all $\om\in \Om$ as 
\begin{align}\label{s_lyap}
{\Lambda}^{t-t_0}_{t_0}(x,\om):=\scaleobj{.9}{\frac{1}{t-t_0}}\log\Vert D\Phi_{t_0,t}^x(\om)\Vert_{2}, 
\end{align}
where $D\Phi_{t_0,t}^x(\om) = D\phi_{t_0,t}(x,\om)$  and $\Vert \ccdot\Vert_{2}$ is the operator spectral~norm. The {\it expected maximal time-forward Lyapunov exponent} is defined as 
\begin{equation}
\bar{\Lambda}^{t-t_0}_{t_0}(x) \,{:=}\, \scaleobj{.9}{\frac{1}{t-t_0}}\E\left[ \log\Vert D\Phi_{t_0,t}^x\Vert_2\right], 
\end{equation}
where the expectation is w.r.t.~the law of $D\Phi_{t_0,\scaleobj{.9}{\Ic}}^x$.

\item[(b)] The {\it empirical stochastic time-forward  Lyapunov exponent} at $x\in\Xt$ over  $[t_0,\,t]\subseteq\Ic$ associated with the realisation of the stochastic flow  $x\mapsto \phi_{t_0,t}(x,\om)$ and the initial  perturbation $y\in \Xt$ is defined for almost all $\om\in \Om$~by 
\vspace{.1cm}
\begin{align}\label{app_s_lyap}
{\Lambda}^{t-t_0}_{t_0}(x,y,\om)& :=\scaleobj{.9}{\frac{1}{t-t_0}}\log\frac{\vert D\Phi_{t_0,t}^x(\om)y\vert+\gamma}{\vert y\vert+\gamma}= \scaleobj{.9}{\frac{1}{t-t_0}}\log\frac{\vert Y^{t_0,y}_{t}(x,\om)\vert+\gamma}{\vert y\vert+\gamma}, \quad 0<\gamma\ll1,
\end{align}
where the integrand is regularised in order to relax the computationally cumbersome constraint $y\ne 0$ when carrying out spatial averaging.  The expectation of ${\Lambda}^{t-t_0}_{t_0}(x,y,\om)$ is denoted as  $\bar{\Lambda}^{t-t_0}_{t_0}(x,y) := \E[{\Lambda}^{t-t_0}_{t_0}(x,y,\ccdot)]$.

\smallskip
\item[(c)] The {\it average empirical time-forward Lyapunov exponent} at  $x\in\Xt$ over $[t_0,\,t]\subseteq\Ic$ associated with the stochastic flow  $\{\phi_{t_0,t}(x,\ccdot)\}_{t\in \Ic}$  is  defined  as 
 \vspace{.1cm}
\begin{align}\label{muE_app_lyap}
\qquad \E^{\tilde\mu^x_{t_0}}\big[\bar{\Lambda}^{t-t_0}_{t_0}(x)\big]:= \scaleobj{.9}{\frac{1}{t-t_0}}\int_{\Xt}\E\left[\log \frac{\vert D\Phi_{t_0,t}^x\,y\vert+\gamma}{\vert y\vert+\gamma}\right]\tilde{\mu}^x_{t_0}(dy), \quad  \;{\tilde\mu}^x_{t_0}\in \PP(\Xt),
\end{align}

\smallskip
\noindent where $\tilde\mu^x_{t_0}$ is a probability measure on the initial perturbation localised at some $x\in \Xt$.
 
\end{itemize}

\begin{rem}\rm Several remarks are in order:
\begin{itemize}[leftmargin = .7cm]
\item[(i)] The finite-time Lyapunov exponents $\Lambda_{t_0}^{t-t_0}, \tilde\Lambda_{t_0}^{t-t_0}$ are often defined at $(t_0,x)$. However, such a definition is not unique and we thus define them for  the interval $[t_0,\,t]\subseteq \Ic$.

\item[(ii)] {In contrast to the linearisation of autonomous systems the finite-time maximal Lyapunov exponents are not guaranteed to be continuous in $(t_0,t)$ \cite{lyap56,sergeev90}. 
Moreover, these  Lyapunov exponents are not guaranteed to exist for $t\rightarrow \infty$. They can exist, for example, if the conditions of the Osedelets theorem are satisfied (see, e.g., \cite{osed68, Arnold1, ruelle79}).  }

\end{itemize}
\end{rem}
\end{definition}

\begin{prop}\label{KL_max_Lyap_stoch}
Consider a stochastic flow  of diffeomorphisms $\{\Solm\}_{t\in \Ic}$  on $\Xt$. Given the measurable family of time-marginal probability measures $(\tilde \mu^{x,\varepsilon}_t)_{t\in \Ic}$  such that  $\tilde \mu^{x,\varepsilon}_t = \mathcal{P}^{ D\Phi^x*}_{t_0,t}\tilde \mu^{x,\varepsilon}_{t_0}$ with $\tilde \mu^{x,\varepsilon}_{t_0}$ localised around $x\in \Xt$ as in (\ref{Mollify}),
 and $\D_\textsc{kl}^{t_0,t}\big(\tilde{\mu}^{x,\varepsilon}_t\|\tilde{\mu}^{x,\varepsilon}_{t_0}\big)<\infty$, the following holds: 
\begin{align}\label{ftle_bnd_1}
-\mathfrak{C}_{-}^{t_0,t}- \D_\textsc{kl}^{t_0,t}\big(\tilde{\mu}^{x,\varepsilon}_t\|\tilde{\mu}^{x,\varepsilon}_{t_0}\big)\leqslant  \E^{\tilde{\mu}^{x,\varepsilon}_{t_0}}\big[\bar{\Lambda}^{t-t_0}_{t_0}(x)\big]\leqslant \D_\textsc{kl}^{t_0,t}\big(\tilde{\mu}^{x,\varepsilon}_t\|\tilde{\mu}^{x,\varepsilon}_{t_0}\big) +\mathfrak{C}_+^{t_0,t}\quad \forall \;t\in \Ic,
\end{align}
where  $0\leqslant\mathfrak{C}_\pm^{t_0,t}(\tilde{\mu}^{x,\varepsilon}_{t_0})<\infty$
 and $\mathfrak{C}_\pm^{t_0,t}\rightarrow 0$ as $\tilde\mu^{x,\varepsilon}_{t_0}\overset{\ast}{\rightharpoonup} \delta_x$. 
Moreover, 
\begin{align}\label{ftle_bnd_2}
-\mathfrak{C}_-^{t_0,t}- \mathcal{R}^{D\Phi}_{\textsc{kl}}(x,t_0,t)\leqslant  \E^{\tilde{\mu}^x_{t_0}}\big[\bar{\Lambda}^{t-t_0}_{t_0}(x)\big]\leqslant \mathcal{R}^{D\Phi}_{\textsc{kl}}(x,t_0,t) +\mathfrak{C}_+^{t-t_0}\quad \forall \;t\in \Ic,
\end{align}

\smallskip
\noindent  where $\mathcal{R}^{D\Phi}_{\textsc{kl}}$ is defined in (\ref{xlocrt}),
and $(\mathcal{P}^{\Phi^x*}_{t_0,t})_{t\in \Ic}$ are the duals of transition evolutions $(\mathcal{P}^{\Phi^x}_{t_0,t})_{t\in \Ic}$ induced by the centred two-point motion  $\{\Phi^{x}_{t_0,t}\}_{t\in \Ic}$. The required assumptions are  automatically satisfied for flows induced by SDEs with sufficiently regular coefficients (such as those in~(\ref{reg_growth_cond})).  
\end{prop}

\smallskip
Furthermore, we have the following bound which should be important in applications when one deals with empirical uncertainties and ensemble averages over multiple experiments performed on the underlying flows (e.g., \cite{beron19,duran18}):
\begin{cor}\label{stoch_ftle} 
Assumming that the conditions of Proposition~\ref{KL_max_Lyap_stoch} hold,  the bound (\ref{ftle_bnd_1}) implies
 \begin{align}\label{sdsdsds}
-\mathfrak{C}_-^{t_0,t}-\E\!\left[\DK^{t_0,t}\big({\Pi}^{x,\varepsilon}_{t_0,t}\|\tilde{\mu}^{x,\varepsilon}_{t_0}\big)\right]\leqslant  \E^{\tilde{\mu}^{x,\varepsilon}_{t_0}}\big[\bar{\Lambda}^{t-t_0}_{t_0}(x)\big]\leqslant \E\!\left[\DK^{t_0,t}\big({\Pi}^{x,\varepsilon}_{t_0,t}\|\tilde{\mu}^{x,\varepsilon}_{t_0}\big)\right]+\mathfrak{C}_+^{t_0,t},
\end{align}
where ${\Pi}^x_{t_0,t} = \tilde{\mu}^{x,\varepsilon}_{t_0}\circ (D\Phi^x_{t_0,t})^{-1}$ is  the forward random measure (\ref{Random measure1}) induced by $\{D\Phi^x_{t_0,t}\}_{t\in \Ic}$.
 \end{cor}
 \begin{rem}\rm
 Recall from Remark \ref{lapsum} that if the dynamics has a stationary ergodic measure $\bar\mu$, then under relatively mild conditions $\E\big[\DK^{t_0,t}\big({\Pi}^{x,\varepsilon}_{t_0,t}\|\bar\mu\big)\big]$ is bounded by the sum of (infinite-time) Lyapunov exponents.  
 \end{rem}
 
 Finally, we derive tight bounds on the average Lyapunov exponents using a more general  information inequality established in \cite{branuda_luq}. 
 
 \begin{prop}\label{tight_KL}
Consider  
a stochastic flow  of $\mathcal{C}^1$-diffeomorphisms $\{\Solm\}_{t\in \Ic}$  on $\Xt$. Given the measurable  family of time-marginal probability measures $(\tilde \mu^{x,\varepsilon}_t)_{t\in \Ic}$, $\tilde \mu^{x,\varepsilon}_t = \mathcal{P}^{ D\Phi^x*}_{t_0,t}\tilde \mu^{x,\varepsilon}_{t_0}$ such that $\D_\textsc{kl}^{t-t_0}\big(\tilde{\mu}^{x,\varepsilon}_t\|\tilde{\mu}^{x,\varepsilon}_{t_0}\big)<\infty$, the following hold: 
\vspace{.15cm}
\begin{align}\label{KK_1}
{\mathcal{K}}^{x,\varepsilon}_{\scaleobj{1}{-}}\big(-\D_\textsc{kl}^{t_0,t}\big(\tilde{\mu}^{x,\varepsilon}_t\|\tilde{\mu}^{x,\varepsilon}_{t_0}\big)\big)\leqslant \E^{\tilde{\mu}^{x,\varepsilon}_{t_0}}\big[\bar{\Lambda}^{t-t_0}_{t_0}(x)\big] \leqslant \mathcal{K}^{x,\varepsilon}_+\big(\D_\textsc{kl}^{t_0,t}\big(\tilde{\mu}^{x,\varepsilon}_t\|\tilde{\mu}^{x,\varepsilon}_{t_0}\big)\big), \hspace{.8cm} \forall \,t\in \Ic,
\end{align}
and 
\begin{align}\label{KK_2}
{\mathcal{K}}^{x,\varepsilon}_{\scaleobj{.8}{-}}\big(-\mathcal{R}^{D\Phi}_{\textsc{kl}}(x,t_0,t)\big)\leqslant \E^{\tilde{\mu}^x_{t_0}}\big[\bar{\Lambda}^{t-t_0}_{t_0}(x)\big] \leqslant \mathcal{K}^{x,\varepsilon}_+\big(\mathcal{R}^{D\Phi}_{\textsc{kl}}(x,t_0,t)\big), \hspace{.8cm} \forall \,t\in \Ic,
\end{align}
where  $\mathcal{K}^{x,\varepsilon}_{\pm}(\pm s)\rightarrow 0$,  
as $s\,{\downarrow}\, 0$ with 
\begin{align}
{\mathcal{K}}^{x,\varepsilon}_{\pm} (s)&= \frac{1}{\mathscr{F}^{-1}_{x,\varepsilon,\pm}(\pm s)}\bigg(\mathbb{E}^{\tilde\mu^{x,\varepsilon}_{t_0}}\Big[\exp\big(\mathscr{F}^{-1}_{x,\varepsilon,\pm}(\pm s)\big) F\big)\Big] -1+s\bigg),\qquad s\geqslant 0,
\end{align}
where  $\mathscr{F}^{-1}_{x,\varepsilon,\pm}$ are the respective  inverses of $\mathscr{F}_{x,\varepsilon,\pm}(\lambda) = \pm\int_{\Xt} (\pm\lambda F-1)e^{\pm\lambda F}d\tilde\mu^{x,\varepsilon}_{t_0}\pm1$, $\lambda>0$, and  $F(y) = \log\big(\vert y\vert+\gamma\big) - \E^{\tilde{\mu}^{x,\varepsilon}_{t_0}}\big[\log\big(\vert y\vert+\gamma\big) \,\big], \;0<\gamma\ll 1$.
The bounds (\ref{KK_1}) and (\ref{KK_2}) are tight and the required assumptions hold  automatically for flows induced by SDEs with sufficiently regular coefficients; e.g.,~(\ref{reg_growth_cond}). 

 \end{prop}

\subsubsection{Proofs of Proposition \ref{KL_max_Lyap_stoch}, Corollary \ref{stoch_ftle} and }\mbox{}

{\it Proof of Proposition \ref{KL_max_Lyap_stoch}.} 
For any strictly convex, locally bounded function $\varphi$ satisfying the normality conditions (\ref{Normality}), the Fenchel-Young inequality $  \eta f \leqslant \varphi^*(f)+\varphi\big(\eta\big) $, $f,\eta\in \mathbb{M}_\infty(\Xt)$, implies that 
\begin{align}\label{fench}
\int_{\Xt} \eta f d\nu -\int_{\Xt} \varphi^*(f)d\nu \leqslant \int_\Xt \varphi(\eta)d\nu \quad \forall\; \nu\in \PP(\Xt).
\end{align}
In particular, given the assumptions of the proposition, we can set $\nu = \tilde\mu^{x,\varepsilon}_{t_0}$ localised around some $x\in \Xt$ (see (\ref{Mollify})) and  $\eta = d\tilde\mu^{x,\varepsilon}_{t}/d\tilde\mu^{x,\varepsilon}_{t_0}$  so that  (\ref{fench}) becomes 
\begin{align}\label{fench2}
\int_{\Xt}  f d\tilde\mu^{\varepsilon}_{t} -\int_{\Xt} \varphi^*(f)d\tilde\mu^{\varepsilon}_{t_0} \leqslant \Df\big(\tilde\mu^{\varepsilon}_{t}\|\tilde\mu^{\varepsilon}_{t_0}\big). 
\end{align}
To simplify notation, we skip the superscripts until they become relevant. 

Next, set $\varphi(u) = u\log u-u+1$ so that (\ref{fench2}) yields a lower bound on the KL-divergence 
\begin{align}\label{fench2}
\int_{\Xt}  f d\mu_{t} -\int_{\Xt} \big(\exp(f)-1\big)d\mu_{t_0} \leqslant \D_\textsc{kl}\big(\mu_{t}\|\mu_{t_0}\big). 
\end{align}
Take $f(y) = \log\big(\vert y\vert+\gamma\big) - \E^{\mu_{t_0}}\big[\log\big(\vert y\vert+\gamma\big) \,\big]$, $0<\gamma\ll 1$. Then, (\ref{fench2}) can be written as 
\begin{align}\label{fench3}
 \E^{\mu_{t}}\big[\log(|y|+\gamma)\,\big] -  \E^{\mu_{t_0}}\big[\log(|y|+\gamma)\,\big]& 
\leqslant \D_\textsc{kl}\big(\mu_{t}\|\mu_{t_0}\big) +\mathfrak{C}_+\big(\mu_{t_0}\big),  
\end{align}
where $\mathfrak{C}_+ (\mu_{t_0})=  \E^{\mu_{t_0}}\big[|y|+\gamma\big]/\exp\big(\,\E^{\mu_{t_0}}\big[\log|y|+\gamma\,\big]\big)-1\geqslant 0$; $\mathfrak{C}_+\geqslant 0$ can be asserted from the Jensen's inequality.

On the other hand,  take $\hat f(y) = -f(y) =  \E^{{\mu}_{t_0}}\big[\log\big(\vert y\vert+\gamma\big)\big] -\log\big(\vert y\vert+\gamma \big)$, $0<\gamma\ll 1$. Then, the bound~(\ref{fench2}) becomes 
\begin{align}\label{fench4}
-\int_{\Xt}  f d\mu_{t} -\int_{\Xt} \big(\exp(-f)-1\big)d\mu_{t_0} \leqslant \D_\textsc{kl}\big(\mu_{t}\|\mu_{t_0}\big),
\end{align}
and, consequently 
\begin{align}\label{fench5}
 \E^{\mu_{t}}\big[\log\big(|y|+\gamma\big)\,\big] -  \E^{\mu_{t_0}}\big[\log\big(|y|+\gamma\big)\big]& \geqslant -\D_\textsc{kl}\big(\mu_{t}\|\mu_{t_0}\big) -\mathfrak{C}_-\big(\mu_{t_0}\big),  
\end{align}
where 
$\mathfrak{C}_- (\mu_{t_0})=  \exp\big(\E^{\mu_{t_0}}\big[\log\big(|y|+\gamma\big)\big]\big)\E^{\mu_{t_0}}\big[\big(|y|+\gamma\big)^{-1}\big]-1\geqslant 0$. Combining  (\ref{fench3}) and (\ref{fench5}) leads to 
\begin{align}\label{lapbnd1}
& -\D_\textsc{kl}\big(\mu_{t}\|\mu_{t_0}\big) -\mathfrak{C}_-\big(\mu_{t_0}\big)\leqslant  \E^{\mu_{t}}\big[\log\big(|y|+\gamma\big)\big] -  \E^{\mu_{t_0}}\big[\log\big(|y|+\gamma\big)\big]\leqslant  \D_\textsc{kl}\big(\mu_{t}\|\mu_{t_0}\big) +\mathfrak{C}_+\big(\mu_{t_0}\big) . 
\end{align}
Thus,  the first claim is obtained by noting that 
\begin{align*}
 \E^{\mu_{t_0}}[\tilde\Lambda_{t_0}^{t-t_0}(x)] &= \int_\Xt \int_{\Om}\tilde \Lambda_{t_0}^{t-t_0}(x,y,\om)\Ms'_{t_0,y}(d\om)\mu_{t_0}(dy)\notag\\[.2cm]
&= \scaleobj{.8}{\frac{1}{|t-t_0|}}\left( \int_\Xt \int_{\Om}\log\big(|Y_t^{t_0,y}|+\gamma\big)\Ms'_{t_0,y}(d\om)\mu_{t_0}(dy)- \int_\Xt \log\big(|y|+\gamma\big)\mu_{t_0}(dy)\right)\notag\\[.2cm]
&= \scaleobj{.8}{\frac{1}{|t-t_0|}}\left( \int_\Xt \log\big(|y|+\gamma\big)\tilde\mu^{\varepsilon}_{t}(dy)- \int_\Xt \log\big(|y|+\gamma\big)\tilde \mu^{\varepsilon}_{t_0}(dy)\right)\notag\\[.4cm]
&=\scaleobj{.8}{\frac{1}{|t-t_0|}}\Big( \E^{\mu_{t}}\big[\log\big(|y|+\gamma\big)\big] -\E^{\mu_{t_0}}\big[\log\big(|y|+\gamma\big)\big]\Big),
\end{align*} 
which follows from  the fact that 
(see (\ref{k_lift})) 
\begin{align}
\E^{\mu_t}[f] &=  \!\scaleobj{.92}{\int_{\Xt}}f(y) \mu_t(dy) =\!\scaleobj{.92}{\int_{\Xt}}\scaleobj{.92}{\int_{\Om}} f\big(D\Phi_{t_0,t}(\om)y\big)\Ms'_{t_0,y}(d\om)\mu_{t_0}(dy), 
\quad\, \,f\in \mathcal{C}^2_\infty(\Xt),
\end{align}
which follows from the representation of the solutions $\mu_t$ of (\ref{fPDE}) via the solution $\Ms'_{t_0,y}$ of the martingale problem for the linearised SDE (see Theorem \ref{prp_centfl} or \cite[Theorem 3.3.4]{Kunitabook})
\begin{equation}\label{lin_sde}
d Y_t^{t_0,y} = (\nabla_x b)\big(t,X_t^{ t_0,x}\big)Y_t^{t_0,y}\,d t+(\nabla_x \sigma)\big(t,X_t^{t_0,x}\big)Y_t^{t_0,y}d W_{t-t_0},\quad Y_{t_0}^{t_0,y}  \sim \tilde\mu^{x,\varepsilon}_{t_0}\in\PP(\Xt),
\end{equation} 
where $X_t^{ t_0,x}$ solves (\ref{gen_sde}); see, e.g., \cite{Figali,Stroock79}. The bound (\ref{ftle_bnd_2}) follows from (\ref{lapbnd1}) with $\mu_t = \tilde\mu^{x,\varepsilon}_{t}$, $\mu_{t_0} = \tilde\mu^{x,\varepsilon}_{t_0}$, and property~(iii) of  $\varphi$-DRF in Theorem \ref{exp_rate} for $\varepsilon\rightarrow 0^+$.
Finally, the boundedness of  $\mathfrak{C}_\pm^{t_0,t} = \mathfrak{C}_\pm/(t-t_0) \geqslant 0$ is clear from the fact that $|y|+\gamma,(|y|+\gamma)^{-1},\log(|y|+\gamma)\in L^1_{loc}(\Rn,\mu_{t_0})$, and $\mathfrak{C}_\pm^{t_0,t}\rightarrow 0$ as $\mu_{t_0}\overset{\ast}{\rightharpoonup} \delta_0$ follows from standard Taylor expansions of $(|y|+\gamma)^{-1}$ and $\log(|y|+\gamma)$.

 \qed

\noindent {\it Proof of Corollary \ref{stoch_ftle}.} 
Given the assumptions, we have by Theorem \ref{mut_rnd} that 
\begin{align*}
\DK(\tilde{\mu}^{x,\varepsilon}_t\|\tilde{\mu}^{x,\varepsilon}_{t_0})= \!\int_{\Xt}\varphi\!\left(\E\big[{\pi}^{x,\varepsilon}_{t_0,t}(y)\big]\right)\!\tilde{\mu}^{x,\varepsilon}_{t_0}(dy) \leqslant \int_{\Xt}\!\!\E\big[\varphi\big({\pi}_{t_0,t}^{x,\varepsilon}(y)\big)\big]\tilde{\mu}^{x,\varepsilon}_{t_0}(dy)
 =\E\big[\DK({\Pi}^{x,\varepsilon}_{t_0,t}\|\tilde{\mu}^{x,\varepsilon}_{t_0})\big],
\end{align*}
which follows by Jensen's inequality and Fubini's theorem. Combining the above  bound with the results in Proposition~\ref{KL_max_Lyap_stoch} leads to 
 \begin{align}\label{sdsdsds_prf}
-\mathfrak{C}_-^{t_0,t}-\E\!\left[\DK^{t_0,t}\big({\Pi}^{x,\varepsilon}_{t_0,t}\|\tilde{\mu}^{x,\varepsilon}_{t_0}\big)\right]\leqslant  \E^{\tilde{\mu}^{x,\varepsilon}_{t_0}}\big[\tilde{\Lambda}^{t-t_0}_{t_0}(x)\big]\leqslant \E\!\left[\DK^{t_0,t}\big({\Pi}^{x,\varepsilon}_{t_0,t}\|\tilde{\mu}^{x,\varepsilon}_{t_0}\big)\right]+\mathfrak{C}_+^{t-t_0},
\end{align}
where $0\leqslant\mathfrak{C}_\pm^{t_0,t}(\tilde{\mu}^{x,\varepsilon}_{t_0})<\infty$
 and $\mathfrak{C}_\pm^{t_0,t}\rightarrow 0$ as $\tilde\mu^{x,\varepsilon}_{t_0}\overset{\ast}{\rightharpoonup} \delta_0$ are the same as those in Proposition \ref{KL_max_Lyap_stoch}.

\medskip
\noindent {\it Proof of Proposition \ref{tight_KL}.} The proof follows by evaluating the general bounds established in \cite[Theorem 3.1 and Proposition 3.3]{branuda_luq} in the form 
\begin{align}\label{Goal}
\mathfrak{B}_{\varphi,-}(\mu\|\nu;f)\leqslant \E^{\mu}[f] - \E^{\nu}[f]\leqslant \mathfrak{B}_{\varphi,+}(\mu\|\nu;f),\quad f\in L_{\varphi*}(\Xt;\nu),
\end{align}
where 
\begin{equation}
\mathfrak{B}_{\varphi,\pm}(\mu\|\nu;f):= \pm\inf_{\lambda > 0}\bigg\{\frac{1}{\lambda}\int_{\Xt} \varphi^*\big(\pm\lambda \big(f-\E^{\nu}[f]\big)\big)d\nu +\frac{1}{\lambda} \Df(\mu\|\nu)\bigg\},
\end{equation}
\begin{align}\label{Goal2}
\mathfrak{B}_{\varphi,+}(\mu\|\nu;f) = \mathfrak{B}_{\varphi,-}(\mu\|\nu;f) = 0, \;\; \textrm{iff}\; \mu = \nu; \; \textrm{ or if $f$\! is constant} \; \nu\text{-a.s.}, 
\end{align}
and $ L_{\varphi*}(\Xt;\nu)$ is a subspace of an Orlicz space associated with $\varphi, \nu$ and given by 
\begin{align*}
L_{\varphi*}(\Xt;\nu):=\Big\{f\in \mathbb{M}(\Xt): \;\forall \gamma>0,\; \int_{\Xt}\varphi^*(\gamma f)d\nu+ \int_{\Xt}\varphi^*(-\gamma f)d\nu<\infty\Big\},\quad \nu\in \PP(\Xt).
\end{align*}
For $\varphi(u) = u\log u-u+1$, $u>0$ in (\ref{dphi}), and $\varphi^*(v) = \exp(v)-1$ we have 
\begin{align}
\mathfrak{B}_{\varphi,\pm}(\mu\|\nu;f)
& = \pm\inf_{\lambda > 0}\frac{1}{\lambda}\bigg\{\int_{\Xt} \exp\big(\pm\lambda \big(f-\E^{\nu}[f]\big)\big)d\nu -1+\D_\textsc{kl}(\mu\|\nu)\bigg\}.
\end{align}
For $f\ne const$, $\mu\ne\nu$ the respective infima  are attained  at $0<\lambda<\infty$ satisfying  
\begin{align}\label{FF_inv}
\mathscr{F}_{\nu,\pm}(\lambda)= \pm\D_\textsc{kl}(\mu\|\nu),
\end{align}
where 
\begin{equation}
\mathscr{F}_{\nu,\pm}(\lambda):=\pm\int_{\Xt} (\pm\lambda F-1)e^{\pm\lambda F}d\nu \pm1,  \qquad F = f-\E^{\nu}[f]. 
\end{equation}
The respective solutions in (\ref{FF_inv}) are unique since by Jensen's inequality $\mathscr{F}_{\nu,+}>0$, $\mathscr{F}_{\nu,-}<0$, and  $\partial_\lambda \mathscr{F}_{\nu,\pm}(\lambda)  = \pm\lambda\int_{\Xt}  F^2e^{\lambda F}d\nu \ne 0$. 
Finally, it remains to  set $\mu = \tilde\mu^{x,\varepsilon}_t$, $\nu = \tilde\mu^{x,\varepsilon}_{t_0}$ in (\ref{Goal}), and  $F(y) = \log\big(\vert y\vert+\gamma\big) - \E^{\tilde{\mu}^{x,\varepsilon}_{t_0}}\big[\log\big(\vert y\vert+\gamma\big) \,\big], 0<\gamma\ll 1$, in (\ref{Goal}) to obtain the claim. \qed

\medskip
\section{Computational aspects}\label{s_comp}

\medskip
\newcommand{\AAA}{\mathfrak{A}}
In order to provide computational examples of the results obtained in \S\ref{s_exp_rates} and \S\ref{s_lyp_fun_stoch},  we to consider  expansion rates based on the KL-divergence (other choices are clearly possible depending on the divergence suitable for applications).  First, we recall a result  which is useful in numerical approximations which are necessarily carried out at a finite resolution. 

\begin{lem}[\textbf{\textit{Set-oriented KL-divergence}} \cite{Dupuis}]\label{duplem}\rm 
Let $\mu, \nu$ be probability measures on a Polish space $\Xt$ ($\Rd$ or $\bar{\mathbb{T}}^d$).   Let $\AAA$ denote the class of all finite measurable partitions of $\Xt$. Then 
\begin{equation*}
\DK(\mu\|\nu)= \sup_{\mathfrak{a}\in\AAA} \sum_{A\in\mathfrak{a}}\mu(A)\log\frac{\mu(A)}{\nu(A)},
\end{equation*} 
where 
\begin{align*}
\mu(A)\log\frac{\mu(A)}{\nu(A)}= 
\begin{cases} 0, \quad\quad \; \text{if}\; \mu(A) =0,\\
+\infty, \quad \text{if}\; \mu(A)>0 \; \text{and}\; \nu(A) =0.
\end{cases}
\end{align*}

\end{lem}
The above lemma provides a more appropriate starting point for numerical approximations compared direct discretisations of divergences that can be shown to converge in the limit if infinitely fine mesh. However, even in this case we chose a single partition instead of considering the supremum over all $\mathfrak{a}\in \AAA$ and we do not consider the error of the approximation since we are concerned with an illustration of the analytical results established in \S\ref{s_exp_rates} and \S\ref{s_lyp_fun_stoch} and not with the accuracy of numerical approximations.

\newcommand{\nn}{{\scriptscriptstyle n}}
\newcommand{\XX}{\mathcal{X}}

Given the above set-oriented representation of KL-divergence, take $\mu_{t_0}\in \PP(\Xt)$  and a countable measurable partition $\AAA$ of $\Xt$ such that $\mu_{t_0}(A)\ne 0$ for all $A\in \AAA$.
Then,  we have 
\begin{align}\label{dddkl}
\DK(\mathcal{P}^*_{t_0,t}\mu_{t_0}\|\mu_{t_0}) &\approx 
\sum_{A\in \AAA}\big(\mathcal{P}^*_{t_0,t}\mu_{t_0}\big)(A)\log\left(\frac{\big(\mathcal{P}^*_{t_0,t}\mu_{t_0}\big)(A)}{\mu_{t_0}(A)}\right), \quad t\in \Ic = [t_0,\;t_0+T),
\end{align}
where $\mathcal{P}^*_{t_0,t}:\; \PP(\Xt)\rightarrow\PP(\Xt)$, defined in (\ref{P*mu}) and  often referred to as the {\it transfer operator},  is the formal $L^2(\Xt,\mu_{t_0})$ dual of the transition evolution $\mathcal{P}_{t_0,t}: \mathbb{M}(\Xt)\rightarrow\mathbb{M}(\Xt)$ in (\ref{calP}).  
 Recall that in our setup  $\mathcal{P}^*_{t_0,t}\mu_{t_0}\ll m_d \;\,\forall t\in \Ic$ if  $\mu_{t_0}\ll m_d$, since the transition evolutions are induced by a global flow of diffeomorphisms generated by solutions of the SDE (\ref{gen_sde}) with sufficiently regular coefficients that satisfy (\ref{reg_growth_cond}); see also Theorem \ref{FKPt1}. 

In order to derive numerical approximations of the FDTR fields $\mathcal{R}_{\textsc{kl}}^\Phi(x,t_0,t)$ defined in (\ref{xlocrt}) 
 we utilise a discrete,  finite-rank approximation of  $\big\{\mathcal{P}^{\Phi^x*}_{t_0,t}\big\}_{t\in \Ic}$ (see  \S\ref{uni_exp_rates})  induced by the centred motion $\big\{\Phi^x_{t_0,t}\}_{t\in \Ic}$ given by (\ref{cnt_fl}). A standard approach to approximation/coarse-graining the action of any transfer operator $\mathcal{P}^*_{t_0,t}$ is Ulam's method which is obtained via a  projection on a  measurable partition of  $\Xt$ (we will return to the specific operators $\{\mathcal{P}^{\Phi^x*}_{t_0,t}\}_{t\in \Ic}$ later on). 
First,  we partition  $\Xt$  into a finite collection of measurable simply connected, closed  sets $\AAA_N =  \{A_1,A_2,\cdots,A_N\}$  such that  $\cup_{i=1}^N \,A_i = \Xt$,  $\mu_{t_0}(A_i)\ne 0$ for all $i =1, \dots, N$, and $\mu_{t_0}(A_i\cap A_j) = 0$ for $i\ne j$. Then, given the transition evolutions $\{\mathcal{P}_{t_0,t}\}_{t\in \Ic}$  and its duals  $\{\mathcal{P}_{t_0,t}^*\}_{t\in \Ic}$, we have that for $\mu_{t_0}\in \PP(\Xt)$
\begin{equation*}
 \int_\Xt f(x)(\mathcal{P}^*_{t_0,t}\mu_{t_0})(d\xi)=\int_\Xt \big(\mathcal{P}_{t_0,t}f\big)(\xi)\mu_{t_0}(d\xi),  \quad \big(\mathcal{P}_{t_0,t}f\big)(\xi) = \E\big[f\big(\phi_{t_0,t}(\xi)\big)\big], \quad  f\in \mathbb{M}_\infty(\Xt),
\end{equation*}
and, in particular, for $f = \I_\XX$, $\XX\in \mathcal{B}(\Xt)$ one has 
\begin{align}\label{coarse_right}
\mu_t(\XX)= (\mathcal{P}^*_{t_0,t}\mu_{t_0})(\XX)  &= \int_\Xt \E\big[\I_\XX\big(\phi_{t_0,t}(\xi)\big)\big]\mu_{t_0}(d\xi) \notag\\
&=\sum_{i=1}^N\int_\Xt \E\big[\I_\XX\big(\phi_{t_0,t}(\xi)\big)\big]\I_{A_i}(\xi)\mu_{t_0}(d\xi)\notag\\
&=\sum_{i=1}^N \E\big[\mu_{t_0}(A_i\cap\phi^{-1}_{t_0,t}(\XX))\big], 
\end{align}
where $\phi^{-1}_{t_0,t}:=(\phi_{t_0,t})^{-1}$. Furthermore, since for $\mu_{t}\in \PP(\Xt)$, $t\in \Ic$
\begin{equation}\label{coarse_left}
\mu_t(\XX) =   \mu_t\left(\sum_{i=1}^N \XX\cap A_i\right) = \sum_{i=1}^N \mu_t\left( \XX\cap A_i\right), 
\end{equation} 
we have, by combining (\ref{coarse_left}) and (\ref{coarse_right}),  that for $\XX = A_j$
\begin{equation*}
\mu_t\left( A_j\right) = \big(\mathcal{P}^*_{t_0,t}\mu_{t_0}\big)(A_j)=\sum_{i=1}^N \E\big[\mu_{t_0}\big(A_i\cap\phi^{-1}_{t_0,t}(A_j)\big)\big].
\end{equation*}
The above relationship can be written concisely as
\begin{equation}\label{mu_chain}
\widehat\mu_t =\widehat\mu_{t_0}P_{t_0,t},
\end{equation} 
where $\widehat \mu_t:=(\hat\mu_t^1, \hat\mu_t^2, \dots, \hat\mu_t^N)$,  $\hat\mu_t^j :=\mu_t(A_j)= \sum_i \hat\mu^i_{t_0}P_{t_0,t}^{ij}$ and the transition matrix
\begin{equation*}
P_{t_0,t}^{ij} = \frac{\E\big[\mu_{t_0}\big(A_i\cap\phi_{t_0,t}^{-1}(A_j)\big)\big]}{\mu_{t_0}(A_i)}. 
\end{equation*}
Discrete-time, iterative updates of $\widehat \mu_t$ can be obtained by  setting  the sequence  $t_0<t_1<\dots<t_M$, $\{t_k\}_{k=0}^M\subset \Ic$ where, for simplicity,  $t_{k+1}-t_k = \Delta t$ so that  
\begin{equation}\label{mu_iter}
\textstyle \widehat\mu_{t_0+k\Delta t} =  \widehat \mu_{t_0}P_{t_0,t_0+k\Delta t}= \widehat \mu_{t_0}\prod_{\ell = 1}^k P_{t_0+(\ell-1) \Delta t ,t_0+\ell\Delta t}. 
\end{equation}
The transition matrix $P_{s,t}, \,s,t\in \Ic$, can be numerically/empirically approximated using standard sampling techniques, especially since in our setting $\mu_{t_0}\ll m_d$ and $\rho_{t_0}=d\mu_{t_0}/dm_d>0$. However, the use of the iterative formula (\ref{mu_iter}) requires estimation of the whole family of operators $\big\{P_{t_0+(\ell-1)\Delta t,t_0+\ell\Delta t}\big\}\!\!\phantom{|}_{\ell =1}^M$ for any given $t_0$, unless the dynamics is autonomous (in which case  $\textstyle \widehat\mu_{t_0+k\Delta t} =  \widehat \mu_{t_0}P_{t_0,t_0+\Delta t}^{(k)}$). Therefore, we directly estimate $\{P_{t_0,t_0+k\Delta t}\}_{k=1}^M$ in (\ref{mu_iter}). 

\medskip  
In order to derive a set-oriented approximation of the KL-divergence rate field in (\ref{xlocrt}) we consider the following steps:
\begin{itemize}[leftmargin=0.4cm]
\item[-] Choose a finite (Lebesgue) measurable partition  $\AAA_N = \{A_i\}_{i=1}^N$ of $\Xt$ (which is necessarily bounded in all computations). 
\item[-] Chose a collection of points $\{x_i\}_{i=1}^N\subset \Xt$ such that $x_i\in A_i$ is the barycentre of $A_i$. 
\item[-] Set the family of  probability measures $\{\tilde \mu^n_{t_0}\}_{n=1}^N$ on $\Xt$, where  $\tilde \mu^n_{t_0} := \tilde\mu^{x_n,\varepsilon_n}_{t_0}$  are the regularised uniform measures (\ref{Mollify})  with $\varepsilon_n\leqslant \textrm{diam}(A_n)$.

\item[-] Define a family of transition matrices 
\begin{equation*}
\tilde P_{t_0,t}^{n,ij} = \frac{\E\big[\tilde\mu^{n}_{t_0}\big(A_i\cap\Phi_{t_0,t}^{x_i,-1}(A_j)\big)\big]}{\tilde\mu^{n}_{t_0}(A_i)}, \quad n=1,\dots, N.
\end{equation*}
\end{itemize}

Refinements of the above procedure are possible and, in fact, necessary in high-dimensional applications (e.g., utilising adaptive techniques akin to those \cite{froy09,froy12,dellnitz97} in a different setup, and for approximating the probability measures as in \cite{chen18,chen18_2}). However, this is beyond the scope of what is an  illustration of the analytical results. More advanced computational approaches to estimating the $\varphi$-DR fields and utilisation of the bound (\ref{bnd2}) as a loss function in statistical and machine/deep learning of coarse-grained models are postponed to a subsequent publication devoted to applications.

\medskip
Given the above setting, we define the set-oriented counterpart of the KL-divergence rate field~(\ref{xlocrt}) as follows:
\begin{definition}[\textbf{\textit{KL-divergence rate field on partition $\AAA$}}]\label{AKL} \rm
The KL-divergence rate on a finite $\mu_{t_0}$-measurable partition  $\AAA_N = \{A_i\}_{i=1}^N$ of $\Xt$ is given by (cf. (\ref{dddkl}))
\begin{equation}\label{RBi}
\RL^{\Phi}_{\textsc{kl},\varepsilon}(A_n,t_0,t) = \scaleobj{.9}{\frac{1}{\vert t-t_0\vert }}\sum_{j=1}^N \tilde\mu^{n}_t(A_j)\log\frac{\tilde\mu^{n}_t(A_j)}{\tilde\mu^{n}_{t_0}(A_j)}=\scaleobj{.9}{\frac{1}{\vert t-t_0\vert }}\sum_{i,j=1}^N\hat\mu^{n,i}_{t_0}\tilde{P}_{t,t_0}^{n,ij}\log \frac{\sum_{i=1}^N\hat\mu^{n,i}_{t_0}\tilde{P}_{t,t_0}^{n,ij}}{\hat\mu^{n,j}_{t_0}}, 
\end{equation}
where $\hat\mu^{n,j}_{t_0}:= \tilde\mu^{n}_{t_0}(A_j)$.
\end{definition}
In what follows we consider a few concrete examples illustrating the finite-time KL-divergence fields on finite partitions of $\Xt$ in the sense of Definition \ref{AKL}.

\subsection{Case study}\label{pics}

 \begin{figure}[t]
\centering\captionsetup{width=.96\linewidth}
\hspace*{-0.1cm}\includegraphics[width = 16cm]{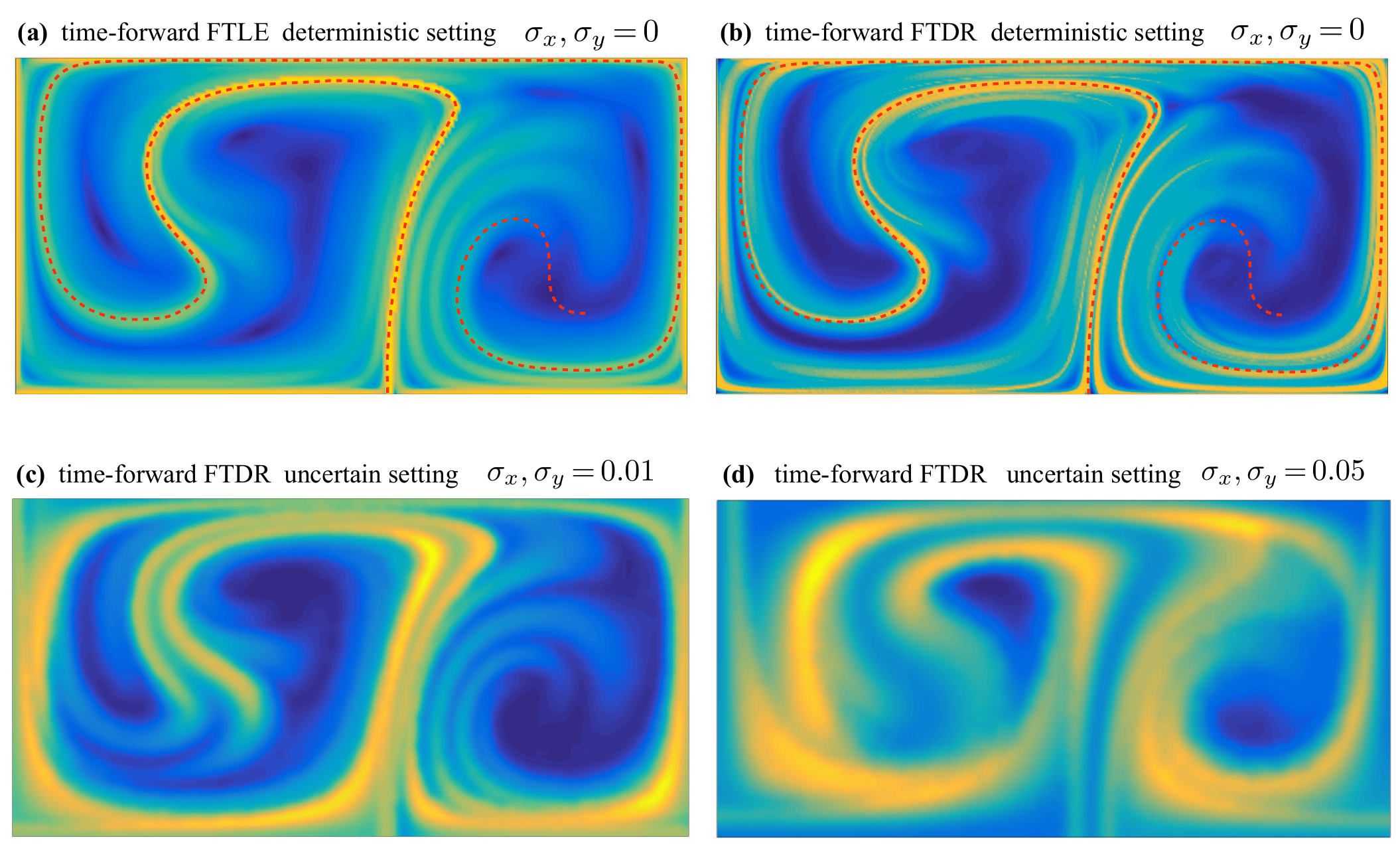}

\vspace{-0.3cm}\caption{\footnotesize  {Time-forward FTLE  (\ref{muE_app_lyap}) and $\textsc{kl}$-DR fields (\ref{xlocrt})  in deterministic and stochastic flows generated  by (\ref{dgyre}).}  The fields are computed over the interval $[t_0,\,T]=[0, \,8]$  for a non-autonomous flows induced by (\ref{dgyre}) with parameters  $\mathfrak{B} \,{=}\,1$, $\delta \,{=}\, 0.25$ and $\Omega \,{=}\, 2$.
 The top row shows the FTLE and $\textsc{kl}$-DR  fields in a deterministic setting; the dashed-red line represents a stable flow-invariant manifold. 
 The good qualitative agreement between agreement the two fields is in line with the bounds (\ref{ftle_bnd_1}) and (\ref{KK_1}). The bottom row shows $\textsc{kl}$-DR fields in uncertain velocity field in  (\ref{dgyre}) for two different values of noise amplitude.  The ridges (yellow) in both types of fields indicate regions of maximal Lagrangian expansion.
  }\label{ftdr_fig}
\end{figure}

 \begin{figure}[t]
\centering\captionsetup{width=.96\linewidth}
\hspace*{-0.1cm}\includegraphics[width = 16cm]{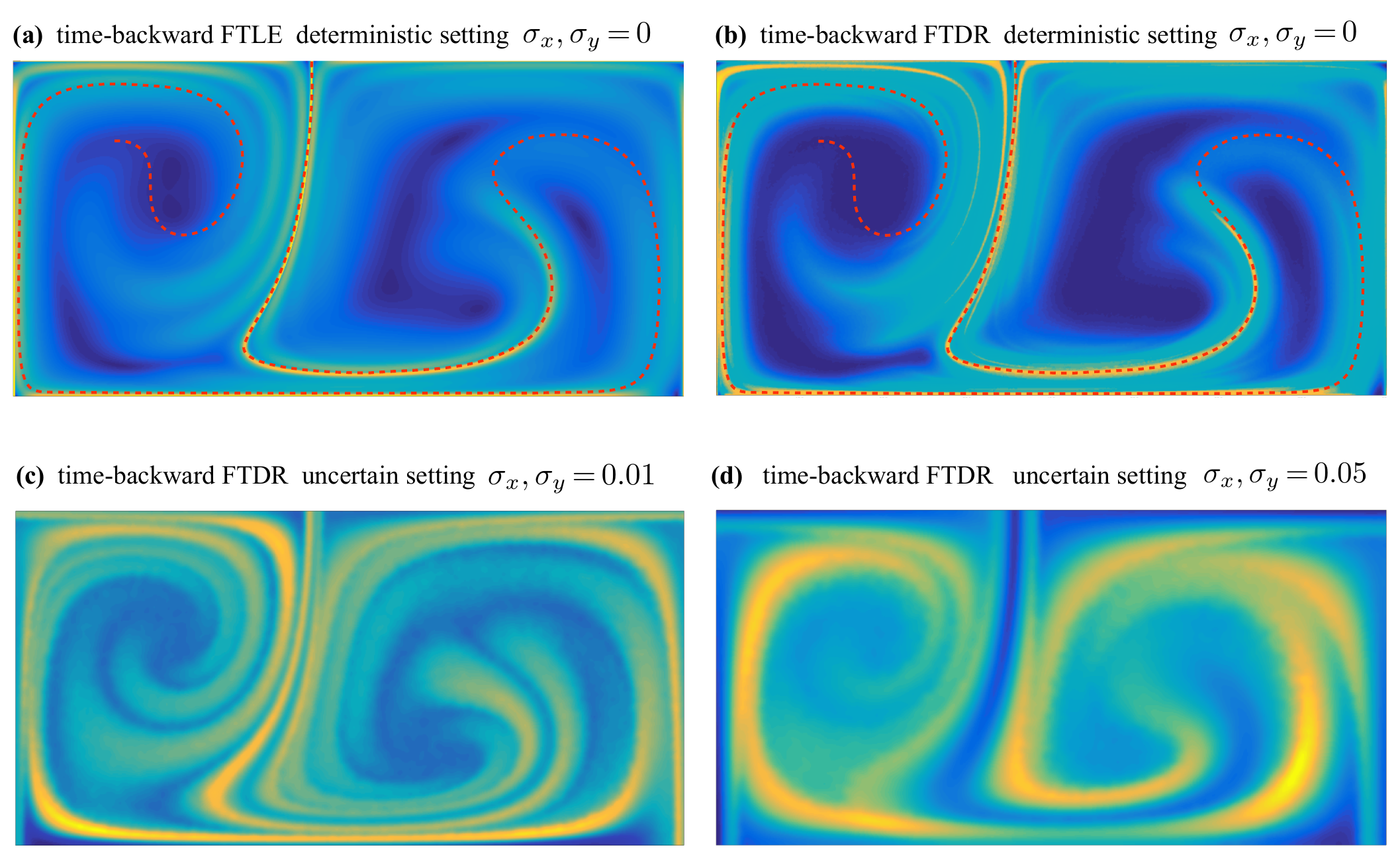}

\vspace{-0.3cm}\caption{\footnotesize {Time-backward FTLE (\ref{muE_app_lyap})  and $\textsc{kl}$-DR fields  (\ref{xlocrt_bck}) in deterministic and stochastic flows generated by  the dynamics (\ref{dgyre}); the fields are computed over the interval $[t_0,\,T]=[-8,\,0]$. The dashed-red line represents the unstable flow-invariant manifold. See figure \ref{ftdr_fig} for more details. } }\label{bck_cell}
\end{figure}

The fields of finite-time divergence rates ($\varphi$-DRF in (\ref{xlocrt})) are needed  in the computational  approach to quantification and mitigation of uncertainty  in Lagrangian predictions through the bounds (\ref{bnd1})-(\ref{bnd2}). The analytical properties of $\varphi$-DRF's were established in Theorem \ref{exp_rate} of \S\ref{uni_exp_rates} and the set-oriented reformulation, which is amenable to numerical approximations, was outlined in earlier in this section (Lemma \ref{duplem} and Definition \ref{AKL}), while the links to Lyapunov exponent fields were established in \S\ref{Almost_sure}. It is worth reiterating that $\varphi$-DRF's are well defined for a general nonlinear stochastic flow and no linearisation in the computations is~required.

 Here, we illustrate our  probabilistic approach for determining the fields of path-based expansion rates via $\textsc{kl}$-DRF (i.e., divergence rate field based on the KL-divergence) on a number of low-dimensional examples, and we compare them with the corresponding FTLE fields.  The relationship between  $\textsc{kl}$-DR and FTLE fields, while interesting, is at best of secondary importance from the perspective of  Lagrangian Uncertainty Quantification (LUQ) in reduced models~\cite{branuda_luq}.  
  However, the abundance and popularity of approaches based on finite-time  Lyapunov exponents in the Lagrangian transport analysis merits a comparison\footnote{\,As remarked at the beginning of \S\ref{s_lyp_fun_stoch}, Lagrangian transport analysis concerns identification of (approximately) flow-invariant structures which represent barriers to Lagrangian transport but are largely irrelevant for our purposes. Ridges of the FTLE field tend to indicate the location of the most prominent barriers  (see e.g., \cite{branwig10,peacock10,hadji17}) subject to a number of caveats (e.g., \cite{branwig10,haller11-b,hadji17}). We do not look for transport barriers.
 }. 
 At the same time, it is  worth stressing that the bounds (\ref{bnd1})-(\ref{bnd2}) utilised in \cite{branuda_luq} generally rely on the global structure of  the expansion/divergence rates (expressed via $\varphi$-DRF's, not FTLE's) and not on a detailed  structure of  their local maximisers and minimisers. Hence,  extensions of Lyapunov exponent-based approaches to the stochastic case  (e.g., via (\ref{muE_app_lyap}))  are not useful in our subsequent applications to LUQ due to the form of the bounds (\ref{bnd1})-(\ref{bnd2}). Moreover, the $\varphi$-DR fields, including  $\textsc{kl}$-DRF  and their set-oriented approximations~(\ref{RBi}), rely on the  evolution of the underlying time-marginal probability measures induced by an arbitrary nonlinear deterministic or stochastic flow.

 \begin{figure}[t]
\centering\captionsetup{width=.96\linewidth}
\hspace*{-0.1cm}\includegraphics[width = 16.3cm]{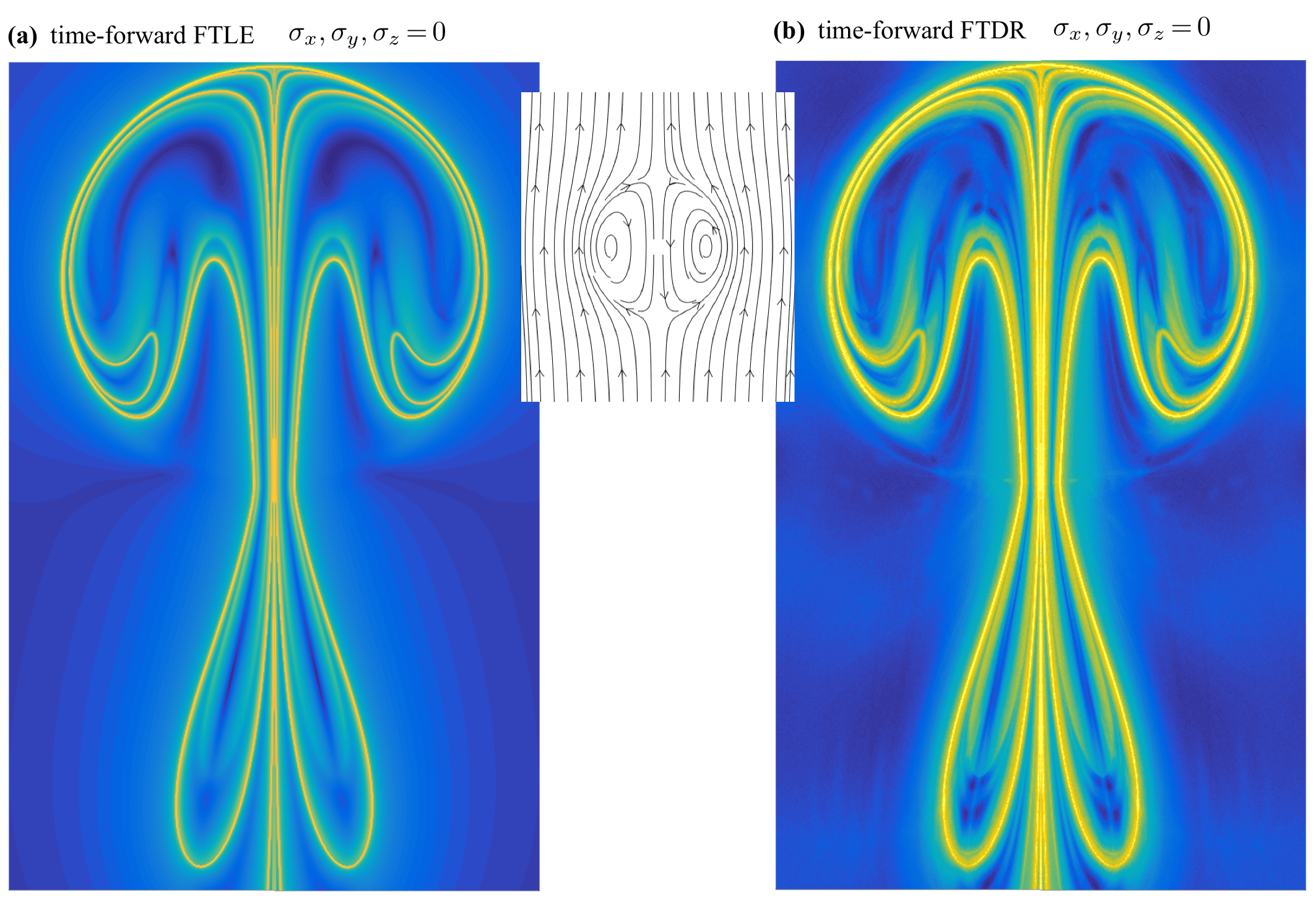}

\vspace{-0.1cm}\caption{\footnotesize Time-forward FTLE  (\ref{muE_app_lyap}) and $\textsc{kl}$-DR fields  (\ref{xlocrt}) in a  deterministic flow generated  by  (\ref{hills}) for $\sigma_x\,{=}\,\sigma_y\,{=}\,\sigma_z\,{=}\,0$.  The fields are computed over the interval  $[t_0,\,T]=[0,\,8]$ for a non-autonomous flow induced by (\ref{hills}). The ridges (yellow) in both types of maps indicate regions of maximal Lagrangian expansion; time-backward fields are shown in figure \ref{bck_hills_det}. See also figure \ref{ftdr_fig}. }\label{frw_hills_det}
\end{figure}

To illustrate our methodology we consider two toy examples. First, we consider the  following 2D system SDEs, which represents a periodically driven `double gyre' flow with additive stochastic noise
\begin{equation}\label{dgyre}
\left.\begin{array}{l}
dx = -\pi \mathfrak{B}\sin(\pi f(t,x))\cos(\pi y)+\sigma_x dW^x_t,\\[.3cm]
dy = \pi \mathfrak{B}\cos(\pi f(t,x))\sin(\pi y)\partial_x f(t,x)+\sigma_y dW^y_t,
\end{array}\right\}
\end{equation}
where  $f(x,t) = \delta \sin(\Omega t)x^2+(1-2\delta\sin(\Omega t))x$; the dynamics is defined on a two-dimensional flat torus, where we take $\Xt:=[0,\,2]\times [0,\,1]$ with doubly-periodic boundary conditions.   

For $\sigma_x\,{=}\,\sigma_y \,{=}\,0$ the dynamics (\ref{dgyre}) reduces to the well-known benchmark for studying Lagrangian transport; note that in the deterministic case the boundary of $\Xt$ is invariant under the induced flow.  When $\delta \,{=}\,0$ in $f(t,x)$ the system (\ref{dgyre}) is autonomous and it has two hyperbolic fixed points on the boundary at $(1,0)$ and $(0,1)$ whose stable and unstable manifolds partition~$\Xt$. For $0<\epsilon\ll 1$ (i.e., the non-autonomous case) these two fixed points morph into  non-trivial hyperbolic trajectories which are confined to the boundary of $\Xt$  which move  along $[1,x]$ and $[0,x]$, $x\in [0,\;2]$ with a period $2\pi/\Omega$. 
For $t=0$ these hyperbolic trajectories are in the same position as the fixed points in the autonomous case, whereas for $t=2\pi/\Omega$ these trajectories are at their extreme locations. The stable and unstable manifolds of these two hyperbolic trajectories form a heteroclinic tangle which has been the focus of many computational studies (e.g., \cite{branwig10,shadden,froy09,lekien07}). 
  In the presence of stochastic noise in (\ref{dgyre}) the above deterministic picture no longer applies but one is still interested in the spatial structure of local expansion rates.

 \begin{figure}[t]
\centering\captionsetup{width=.96\linewidth}
\hspace*{-0.1cm}\includegraphics[width = 16.3cm]{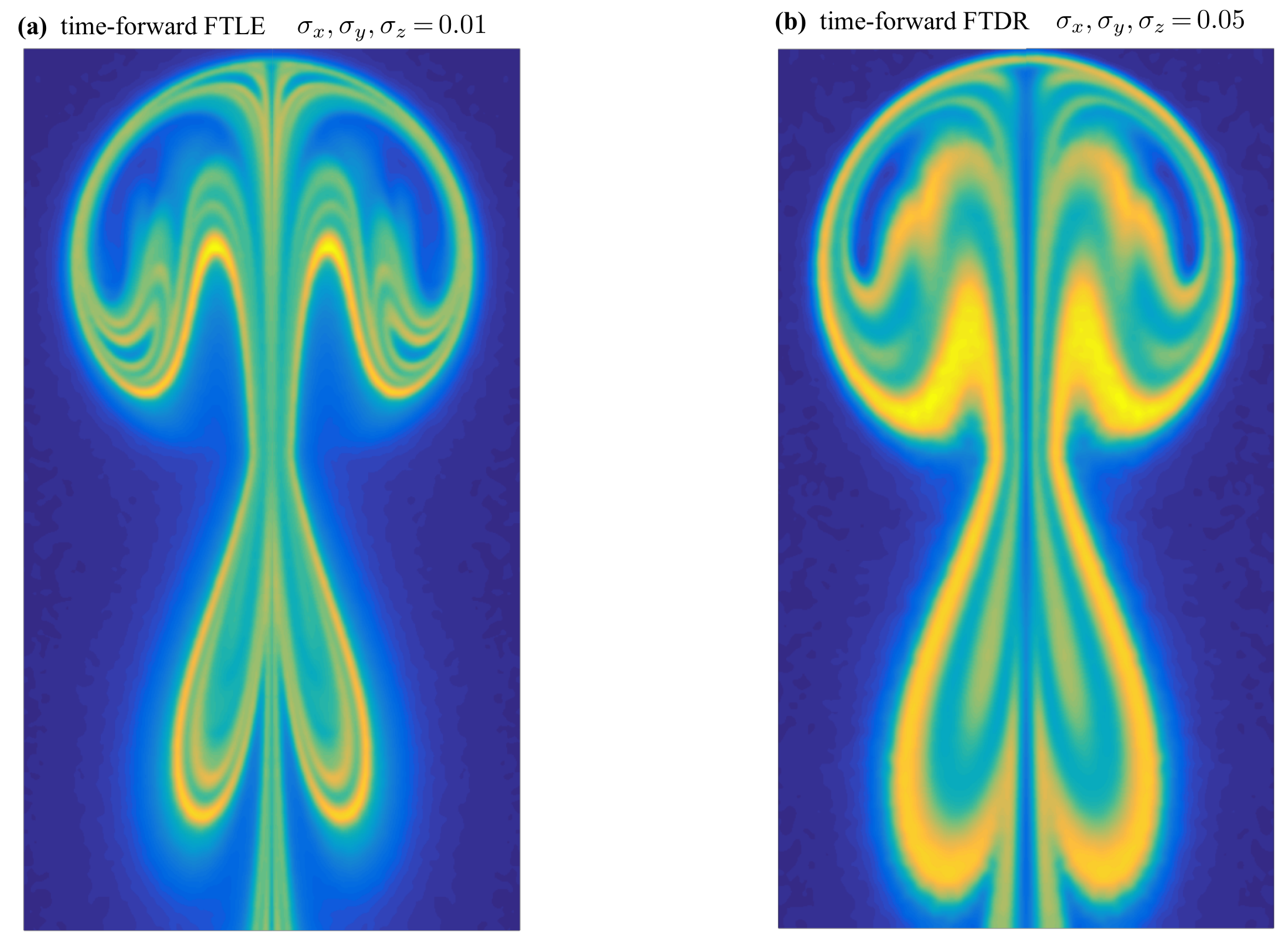}

\vspace*{-0.2cm}\caption{\footnotesize Time-forward FTLE  (\ref{muE_app_lyap}) and $\textsc{kl}$-DR fields  (\ref{xlocrt}) in a stochastic flow generated  by  (\ref{hills}) for two different amplitudes of additive noise.  The fields are computed over the interval  $[t_0,\,T]=[0,\,8]$ for a non-autonomous flow induced by (\ref{hills}). The ridges (yellow) indicate regions of maximal Lagrangian expansion. Compare these fields with those for  the deterministic dynamics of (\ref{hills}) shown in figure \ref{frw_hills_det}.  }\label{frw_hills_stoch}
\end{figure}

\begin{figure}[t]
\centering\captionsetup{width=.96\linewidth}
\hspace*{-0.1cm}\includegraphics[width = 16.3cm]{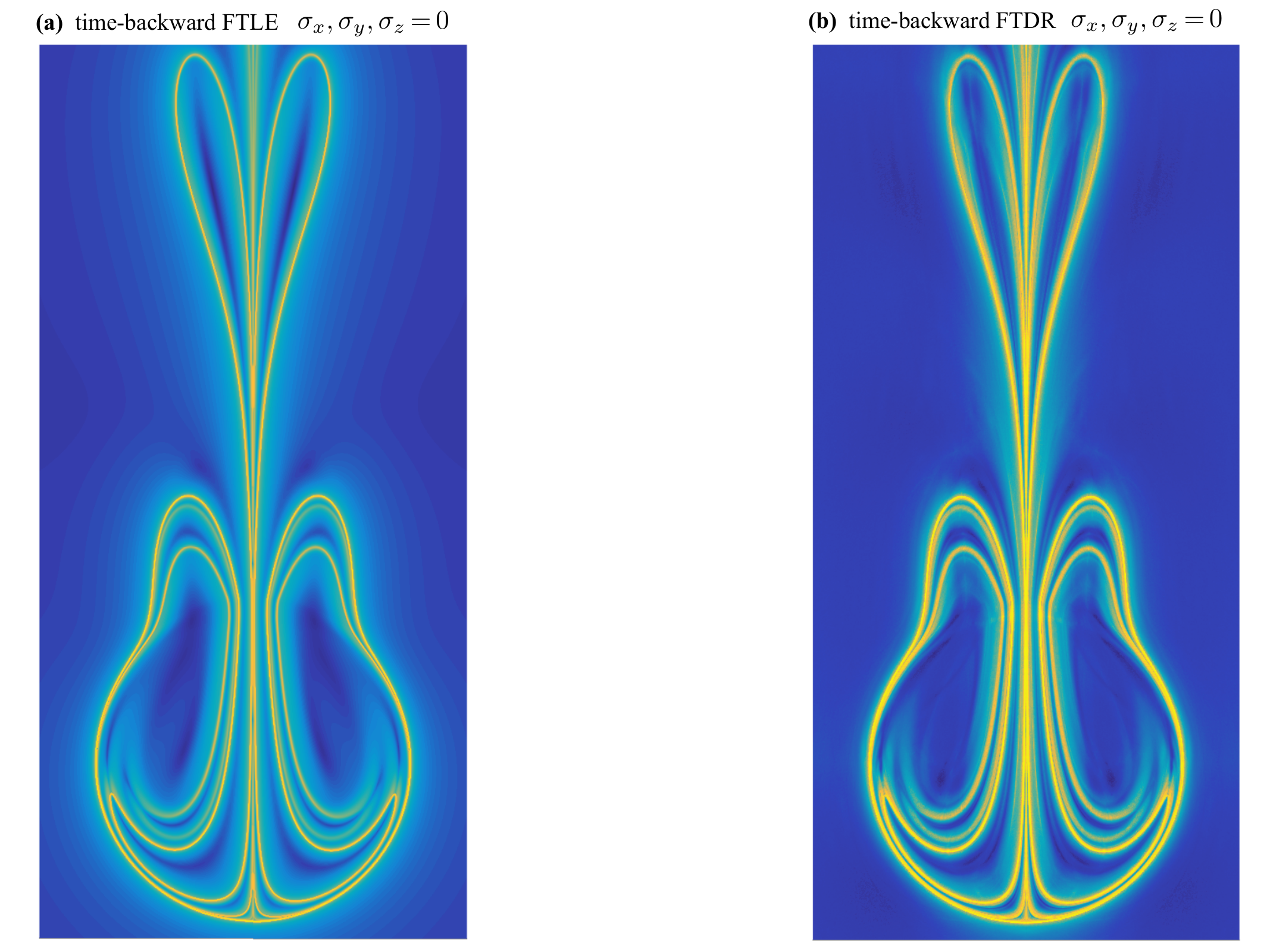}

\vspace{0.0cm}\caption{\footnotesize   Time-backward FTLE  (\ref{muE_app_lyap}) and $\textsc{kl}$-DR fields  (\ref{xlocrt_bck}) and (\ref{dkl}) in a deterministic  flow generated  by (\ref{hills}) for $\sigma_x\,{=}\,\sigma_y\,{=}\,\sigma_z\,{=}\,0$. The fields are computed over the interval  $[t_0,\,T]=[-8,\,0]$ for a non-autonomous flow induced by (\ref{hills}). The ridges (yellow) in both types of maps indicate regions of maximal Lagrangian expansion (backwards in time). See also figure \ref{frw_hills_det} for the associated time-forward FDTR fields. }\label{bck_hills_det}

\vspace*{-.5cm}
\end{figure}

In Figures \ref{ftdr_fig} and \ref{bck_cell} we compare the numerical approximations of $\textsc{kl}$-DRF ((\ref{xlocrt}) and (\ref{RBi}))  and the average FTLE fields (\ref{muE_app_lyap});  where appropriate, we also determine the geometry of stable and unstable manifolds of the two hyperbolic trajectories confined to the top and bottom boundary of $\Xt$ based on methods described in \cite{mswi,msw2,branwig10}. We fix parameter values in (\ref{dgyre}) as $\mathfrak{B} \,{=}\,1$, $\delta \,{=}\, 0.25$ and $\Omega \,{=}\, 2$ and the results show  FTLE and FTDR fields computed over the time interval $[t_0,\;T] = [0,\, 8]$ for the forward fields, and over $[t_0,\;T] =[-8,\, 0]$ for the backward fields. To create a numerical approximation of the transfer operators $(\mathcal{P}^{\Phi^x*}_{t,t_0})_{t\in \Ic}$, the domain $\Xt$ is partitioned into  $N = 1000^2$ boxes.  The transition matrix $\tilde{\mathcal{P}}_{t,t_0}^n$ in (\ref{RBi}) is estimated  by numerically integrating the system with respect to $N^2 = 1000^2$ inner grid points per box using the 4th-order Runge-Kutta scheme in the deterministic case or the Euler-Maruyama scheme in the stochastic case (with time step chosen to guarantee numerical convergence and desired  accuracy). 

As a final example, consider the dynamical system associated with the well-known solution of the Euler's equation of an inviscid incompressible fluid flow given by the Hill's spherical vortex~(e.g.,~\cite{batchelor}). The SDE associated with stochastically perturbed Hill's vortex in Cartesian coordinates are 
\begin{equation}\label{hills}
\left.\begin{array}{l}
dx =  \big(u_r(t) \sin\Theta+u_\Theta(t) \cos\Theta\big)\cos\Phi dx+\sigma_xdW_t^x,\\[.2cm]
dy = \big(u_r(t)\sin\Theta+u_\Theta(t)\cos\Theta\big)\sin\Phi dt+\sigma_ydW_t^y,\\[.2cm]
dz  =\big(u_r(t)\cos\Theta-u_\Theta(t)\sin\Theta\big)dt+\sigma_zdW_t^z,
\end{array}\right\}
\end{equation}   
where $r = \sqrt{x^2+y^2+z^2}$, $\Theta = \textrm{acos(z/r)}$, $\Phi = \textrm{acos}(x/\sqrt{x^2+y^2})$. Assuming that $a>0$ denotes the radius of the vortex, the velocity components in the spherical coordinates are 
\begin{equation*}
u_r(t) = \begin{cases} \hspace{.5cm}U\big(1-a^3(t)/r^3\big)\cos\Theta & \textrm{if} \quad r \geqslant a,\\[.2cm] -\textstyle{\frac{3}{2}} U\big(1-r^2/a^2(t)\big)\cos\Theta & \textrm{if} \quad r< a,\end{cases}
\end{equation*}
\begin{equation*}
u_\Theta(t)= \begin{cases}  -U\big(1+a(t)^3/(2r^3)\big)\sin\Theta & \textrm{if} \quad r \geqslant a,\\[.2cm] \textstyle{\frac{3}{2}}U\big(1-2r^2/a^2(t)\big)\sin\Theta & \textrm{if} \quad r< a, \end{cases}
\end{equation*}
where we set $U=2$.
In the deterministic, autonomous case the steady Hill's vortex flow has two hyperbolic fixed points 
\begin{equation*}
h_1 = (0,0,-a)^T, \quad h_2 = (0,0,a)^T, 
\end{equation*}
which are located on the (flow-invariant) axis of symmetry $\pmb{e}_z$ of the vortex. The fixed point $h_1$ has a two-dimensional unstable manifold in $\R^3$ and the fixed point $h_2$ has a two-dimensional stable manifold in $\R^3$ (\cite{branwig10}).
We consider a non-autonomous version of the dynamics (\ref{hills}) and~set 
\begin{equation*}
a(t) = 2+0.12 \sin(2.2 t).
\end{equation*}
This time-dependence does not break the axial symmetry of the dynamics and, consequently, any plane containing $\pmb{e}_z$ is invariant under the (deterministic) flow. This, in turn, implies that the stable and unstable manifolds of hyperbolic trajectories, $\gamma_{h_1}(t), \gamma_{h_2}(t)$, confined to the axis of symmetry are foliated by planes containing $\pmb{e}_z$ (i.e., at any fixed $t$, the stable and unstable manifolds of $\gamma_{h_1}(t), \gamma_{h_2}(t)$ intersect any invariant plane along a 1D curve). These trajectories can be computed using the algorithms of \cite{idw,jsw}.  Their stable and unstable manifolds are computed as in the previous examples using techniques described in \cite{mswi,msw2,branwig10}. In Figures \ref{frw_hills_det}--\ref{bck_hills_det} we compare numerical approximations of $\textsc{kl}$-DRF fields and FTLE fields associated with the dynamics~(\ref{hills}); figures \ref{frw_hills_det} and \ref{bck_hills_det} correspond to the deterministic dynamics in (\ref{hills}) with $\sigma_x\,{=}\,\sigma_y\,{=}\,\sigma_z\,{=}\,0$, and figure~\ref{frw_hills_stoch} illustrates the structure of FTDR fields in the stochastic  case.

For deterministic version of the dynamics in (\ref{dgyre}) and (\ref{hills})  both $\textsc{kl}$-DRF and FTLE fields  match well the location of the stable and unstable manifolds of the two hyperbolic trajectories (see figures \ref{ftdr_fig}(a,b), \ref{bck_cell}(a,b), \ref{frw_hills_det}, and \ref{bck_hills_det}) which is in line with the bounds (\ref{ftle_bnd_1}) and (\ref{KK_1}).
Remnants of the deterministic Lagrangian structures persist for small values of the noise amplitude but they gradually degrade with increasing amplitude of the stochastic noise  (Figures~\ref{ftdr_fig}(c,d), \ref{bck_cell}(c,d),~\ref{frw_hills_stoch}),  and/or  when the fields are computed over increasingly long time intervals.  This behaviour is to be expected and it does invalidate the bounds (\ref{bnd2}) not their utility for tuning reduced models at finite times.  For example,  for ergodic dynamics the support of the measures $\mathcal{P}_{t,t_0}^{\Phi^x*}\tilde\mu^{x,\varepsilon}_{t_0}$ increases all of $\Xt$ as $t\,{\rightarrow}\,\infty$. Consequently, the KL-divergence $\D_\textsc{kl}(\mathcal{P}_{t,t_0}^{\Phi^x*}\tilde\mu^{x,\varepsilon}_{t_0}\|\tilde\mu^{x,\varepsilon}_{t_0})$ approaches its asymptotic value $\D_\textsc{kl}(\bar\mu\|\tilde\mu^{x,\varepsilon}_{t_0})$ for any $x\in \Xt$ so that the corresponding expansion rate field  $\RL^\Phi_{\textsc{kl},\varepsilon}(x;t_0,t)$ in (\ref{xlocrt}) decays to zero as $t\rightarrow \infty$ after reaching some intermediate maximum value.   

More advanced  approaches to estimating the $\varphi$-DR fields (utilising adaptive techniques akin to those \cite{froy09,froy12,dellnitz97} and for approximating the probability measures in higher dimension as in \cite{chen18,chen18_2}) and utilisation of the bound (\ref{bnd2}) as a loss function in statistical and machine/deep learning of coarse-grained models are postponed to a subsequent publication devoted to applications.

\section{Conclusions and future work}\label{s_concl}
We developed a new probabilistic framework for characterising nonlinear  trajectorial expansion rates in non-autonomous stochastic dynamical systems that can be  defined over a finite time interval and used for the subsequent  uncertainty quantification in Lagrangian (trajectory-based) predictions. These nonlinear expansion rates are quantified via the family of  $\varphi$-divergences  between probability measures induced by the  law of the stochastic flow associated with the underlying dynamics.  The stochastic flow formulation   elucidates connections between the path space measures, the evolution of their time-marginal probability measures, and  path-based local stretching rates. We constructed scalar fields of local finite-time divergence/expansion rates termed  ($\varphi$-DRF), showed their existence for general stochastic flows, continuity in time and the spatial parameter, and we derived a hierarchy of bounds  that are important in practical uncertainty quantification tasks using $\varphi$-DRF. These fields can be subsequently  combined with information inequalities (\ref{bnd1})-(\ref{bnd2}) derived in~\cite{branuda_luq} to mitigate the uncertainty in path-based observables estimated from approximations of the true dynamics in a way that is amenable to algorithmic implementations, and it can be  utilised in information-geometric analysis of statistical estimation and inference on families of models, as well as in machine/deep learning  approaches to the identification of accurate models for Lagrangian predictions.

Moreover, for the particular case of the Kullback-Leibler divergence the corresponding expansion rates were linked to the Lyapunov exponents for probability measures, as well as the finite-time Lyapunov exponents for path-based observables; the latter  are commonly used  for estimating expansion rates in Lagrangian transport considerations and our results should be relevant when considering some of its aspects in a more general probabilistic setting. 

\smallskip
This work was motivated by the desire to quantify the evolution of  uncertainty and improving path-based/Lagrangian predictions in complex dynamical systems based on simplified, data-driven models, especially in geophysical applications.
A follow-up study combining the above results with the framework for Lagrangian uncertainty quantification \cite{branuda_luq}, as outlined in \S\ref{out_main}, will focus on applications of these tools for deriving optimised models for accurate Lagrangian predictions over different time-horizons, and on understanding the role of both the Eulerian characteristics of the underlying dynamics and the dominant Lagrangian structures on the skill of path-based predictions in applications.  In particular, the use of the bound (\ref{bnd2}) as a loss function in machine/deep learning of coarse-grained models is of special interest and will be discussed in a separate publication.

\bigskip
\noindent{\bf Acknowledgements.} M.B. acknowledges the support of Office of Naval Research grant ONR N00014-15-1-2351. K.U. is funded as a postdoctoral researcher on the above grant.

\appendix

\appendix

\newpage
\section{Further proofs}\label{f_prf}
\addtocontents{toc}{\protect\setcounter{tocdepth}{1}}

\subsection{Proof of Proposition \ref{al_pi}}\label{app_al_pi}
{The procedure  is similar to derivations in  \cite[\S4.3]{Kunitabook}. Densities of  $\Pi_{t,s}$ and $\Check{\Pi}_{t,s}$ w.r.t. $\mu_s\in \PP(\Xt)$ are related by the stochastic flow of diffeomoprphisms $\big\{\phi_{s,t}\big\}_{s,t\in \Ic}$.  In order to see this,  recall the change-of-variable formula
\begin{equation}\label{Cha}
\int_{\phi_{t,s}(A,\om)}f(y)m_d(dy) = \int_{A} f(\phi_{t,s}(x,\om))\big|\!\det D \phi_{t,s}(x,\om)\big| m_d(dx) \quad \mathbb{P}\textrm{-a.s.}, \quad f_\infty\in \mathcal{C}(\Xt).
\end{equation}
 If $\mu_s$ has a strictly positive density $\rho_s(x) = d\mu_s/dm_d>0$,  then  by (\ref{Cha})  the measure $\Check{\Pi}_{t,s}$ satisfies 
 \begin{align*}
 \Check{\Pi}_{t,s}(A,\omega)= \mu_s(\phi_{t,s}(A,\omega))&=
 \int_{A}\rho_s(\phi_{t,s}(x,\omega))\rho_s^{-1}(x)\big|\!\det D\phi_{t,s}(x,\omega)\big|\mu_s(dx), \quad \forall A\in \Bb(\Xt), 
 \end{align*}
 where  $D\phi_{t,s}(x,\om)$ denotes  the Jacobian matrix of $\phi_{t,s}(x,\omega)$ and $|\det D\phi_{t,s}(x,\om)|\neq 0 \;\; \p\textrm{-a.s.,}$  and  it follows that
 \begin{equation}\label{formu1}
 \alpha_{t,s}(x,\omega):=\frac{d\Check{\Pi}_{t,s}}{d\mu_s}= \rho^{-1}_s(x)\rho_s(\phi_{t,s}(x,\omega))\big|\!\det D\phi_{t,s}(x,\omega)\big|,
 \end{equation}
 is the Radon-Nikodym derivative of the probability measure $\check{\Pi}_{t,s}$ with respect to $\mu_s$. In fact, if $\rho_s\in \mathcal{C}^{3}(\Xt;\Rp)\cap L^1(\Xt,\R^+)$, then $\alpha_{t,s}(x)$ admits the  integral representation \cite[Lemma 4.3.1]{Kunitabook}
 \begin{align}\label{Dens_pr}
 \nonumber \alpha_{t,s}(x,\omega)&= \exp\bigg\{\int_{s}^t\rho_s^{-1}(\phi_{s,\xi}(x,\omega))\, \text{div}(\rho_s b)(\xi,\phi_{\xi,s}(x,\omega))d\xi\\ &\hspace{2cm}+\int_{s}^t\sum_{k=1}^d \rho_s^{-1}(\phi_{s,\xi}(x,\omega))\,\text{div}(\rho_s\sigma_k)(\xi,\phi_{s,\xi}(x,\omega))\circ dW_\xi^k(\omega)\bigg\},
 \end{align}}
where $\sigma_{k}$ denote the columns of $\sigma$ in (\ref{gen_sde}). Similarly, the random measure $\Pi_{t,s}$ is given by   
 \begin{align*}
 \Pi_{t,s}(A,\omega) = \mu_s(\phi_{t,s}^{-1}(A,\omega)) 
 &= \int_A\frac{\rho_s(\phi_{t,s}^{-1}(x,\omega))}{ \vert \!\det D\phi_{t,s}(\phi_{t,s}^{-1}(x,\omega),\omega)\vert \rho_s(x)}\mu_s(dx),
 \end{align*}
 so that the Radon-Nikodym derivative of $\Pi_{t,s}$ with respect to $\mu_s$ is given by 
 \begin{align}\label{Forw1}
 \pi_{t,s}(x,\omega) := \frac{d\Pi_{t,s}}{d\mu} = \rho_s^{-1}(x)\rho_s(\phi_{t,s}^{-1}(x,\omega)) \big|\! \det D\phi_{t,s}(\phi_{t,s}^{-1}(x,\omega),\omega)\big|^{-1}.
 \end{align}
The two densities in  (\ref{formu1}) and (\ref{Forw1}) are related  by 
 \begin{align}\label{rel1}
 \pi_{t,s}^{-1}(x,\omega) = \alpha_{t,s}(\phi_{t,s}^{-1}(x,\omega),\omega); 
 \end{align}
 the above relationship is crucial as it is generally not possible to write the stochastic integral governing $\pi_{t,s}$ which we need in the derivations of \S\ref{s_exp_rates}.

\subsection{Proof of Theorem \ref{mut_rnd}}\label{app_mutrnd}
 We present the proof of  the  second part of the theorem to make this article self-contained  (see \cite{Kunitabook, Lu09} for more details). Under the imposed regularity assumptions on the coefficients of~(\ref{gen_sde}),  the maps $\Solm(\ccdot,\omega)$  and  $\Solm^{-1}(\ccdot,\om)$ are  $\mathcal{C}^1$-\,diffeomorphisms (e.g., \cite{Kunitanote, Kunitabook}) and there exist $K_\om\in \cap_{p>1}L^p(\Om)$ and $\alpha>1$ such that 
 \begin{align}\label{Spat_est}
 \vert \Solm^{-1}(x,\om)\vert\leqslant K_\om (1+\vert x\vert^{\alpha}), \quad \text{for all $(t_0,t,x)$}\in\Ic\times \Ic\times\Xt.
 \end{align}
 Let $B_r\subset \Xt\subseteq \Rd$ be a ball of radius $r>0$. The spatial estimate (\ref{Spat_est}) implies  that 
 \begin{align*}
 \Solm^{-1}(B_r,\om)\subset B_{K_\om(1+r^{\alpha})},
 \end{align*}
 and, due to the  monotone property of the finite positive measure $\mu_{t_0},$ we have
 \begin{align*}
 \mu_{t_0}(\Solm^{-1}(B_r,\om))\leqslant \mu_{t_0}(B_{K_\om|1+r^{\alpha}|}) = \!\!\int_{B_{K_\om|1+r^{\alpha}|}}\hspace{-.6cm}\rho_{t_0}(x)m_d(dx)\leqslant \!\!\sup_{x\in B_{K_\om|1+r^{\alpha}|}}\!\!\!\!\big\{\rho_{t_0}(x)\big\}V_d(1)K_\om^d|1+r^{\alpha}|^d,
 \end{align*}
 where $V_d(1)$ is the volume of the unit ball in $\Rd$.  Therefore, 
 \begin{align}\label{sp2}
 \E\Big(\sup_{t\in \Ic}\,\big[\mu_{t_0}(\Solm^{-1}(B_r,\om))\big]^p\Big)\leqslant V_d(1)^p(1+r^{\alpha})^{pd} \,\E\big[ M^pK_\om^{pd}\big]<\infty, \quad \forall \;p> 1,
 \end{align}
 where $M:= \sup\{\rho_{t_0}(x): x\in B_{K_\om|1+r^{\alpha}|}\}<\infty$.  
This procedure can be extended to a set of the form $\bigcup_{x\in S}B_{r(x)}(x)$ for some countable subset of $S\subset\Xt$ provided that the density of the initial measure is sufficiently smooth and strictly positive (see, e.g., \cite{Lu09} for details).
Thus,  by the covering lemma, the above construction can be extended to any Borel measurable subset $A\in\Bb(\Xt)$ and we have that $\mu_t\ll m_d, \; \p\,\textrm{-a.s.}$ \,for \;$t\in\Ic.$
    
\smallskip
   Finally, in order to derive (\ref{rho_pi}), note that for any $f\in \mathcal{C}_{c}^{\infty}(\Xt),$ 
\begin{align}\label{sp13}
 \int_{\Xt}f(\Solm(x,\om))\mu_{t_0}(dx) = \int_{\Xt}f(y)\frac{d(\mu_{t_0}\circ \Solm^{-1})(y)}{d\mu_{t_0}}\mu_{t_0}(dy)= \int_{\Xt}f(y)\pi_{t_0,t}(y,\om)\mu_{t_0}(dy).
\end{align}
On the other hand,  we have  
\begin{equation}\label{sp12}
\int_{\Xt}f(x)\mu_t(dx) = \int_{\Xt}\E\big[f\big(\Solm(x,\ccdot)\big)\big]\mu_{t_0}(dx). 
\end{equation}
Combining (\ref{sp13}) and (\ref{sp12}) leads to   
\begin{align*}
\int_{\Xt}f(x)\mu_t(dx) = \int_{\Xt}f(x)\E\big[\pi_{t_0,t}(x,\ccdot)\big]\rho_{t_0}(x)m_d(dx), \quad f\in\mathcal{C}_c^{\infty}(\Xt),
\end{align*}
which gives 
\begin{align*}
\hspace{4cm}\mu_t(dx) = \rho_{t_0}(x)\E\big[\pi_{t_0,t}(x,\ccdot)\big]m_d(dx),\quad t\in\Ic. \hspace{3.5cm}\qed
\end{align*}

\subsection{Comments on Theorem \ref{FKPt1}}\label{app_FKPt1}
The imposed regularity of the coefficients $(\brr,\sigma)$ guarantees  the existence of global strong solutions of the SDE (\ref{gen_sde}) and the existence of the extremal solution\footnote{\,Here, an element  $P_0\in \Delta\subset \mathcal{P}(\mathcal{C}(\Ic; \Xt))$ is said to be {\it extremal}  if $P_0 = \lambda P_1+ (1-\lambda)P_2,$ for some $P_0, \; P_1\in \Delta$ and $0<\lambda<1,$ implies that $P_0 = P_1=P_2$ i.e., a  point in the  convex set $\Delta$ which is not an interior point of any line segment lying entirely in $\Delta$.}  of the martingale problem associated with $\LG_t$ for $m_d$\,-\,a.e.~initial condition $x\in \Xt$, leading to  the existence and uniqueness of the  weak$-^*$ solution of the forward Kolmogorov equation (\ref{fPDE});   e.g.,  \cite{ Figali,Roc-Kry, Stroock79}. Then,  the existence of the Lebesgue density $(\rho_t)_{t\in \Ic}$ solving (\ref{rho_FP}) when $\mu_{t_0}\ll m_d$, and the representation~(\ref{rho_pi}) follows from the existence of a stochastic flow of solutions of (\ref{gen_sde}) and Theorem \ref{mut_rnd}. If $\rho_{t_0} = d\mu_{t_0}/dm_d\in \mathcal{C}^2_\infty(\Xt; \Rp)\cap L^1(\Xt;\R^+)$, the existence of unique solutions $\rho_{t}\in \mathcal{C}^2_\infty(\Xt; \Rp)\cap L^1(\Xt;\R^+)$, $t\in \Ic$, follows from the Feller property of transition evolutions $(\mathcal{P}_{s,t})_{t\geqslant s}$ induced by the stochastic flow and It\^o formula. Finally,  the last statement follows from the smoothing property of transition evolutions in such a setup.

\subsection{Proof of Proposition \ref{prp_centfl}}\label{app_prp_centfl}

Given the assumptions of the proposition, the derivative flow $\big\{D\phi_{t_0,t}(\ccdot,\ccdot)\big\}_{t\in \Ic}$ exists and is bounded from below due to Theorem \ref{der_fl_thm}, and the map $y\mapsto \Phi^x_{t,t_0}(y,\om)$ in (\ref{cnt_fl}) is a $\mathcal{C}^1$\,-diffeomprphism.  Differentiating  $\Phi^x_{t_0,t}(y,\om)$  with respect to $y\in \Xt$  at $y=0$, (see, e.g.,  \cite[Theorem 3.3.4]{Kunitabook}) yields for all $t\in \Ic$
\begin{align}\label{Var_Centred}
D\Phi^{x}_{t_0,t}(\om) = D\Solm(x,\om) \quad  \p\,\textrm{-\,a.s.}
\end{align}
The relationship in (\ref{Var_Centred})  implies, in particular, that the solutions $y_t = D\Solm(x,\om)y$ associated with (\ref{var_dphi}) coincide $ \p\,\textrm{-\,a.s.}$ with  $y_t = D\Phi^x_{t,t_0}(\om)y$ for\footnote{\,Note that $y\in TM_{x+v}(\Xt)$ and $y_t\in TM_{\phi_{t_0,t}(x+v,\om)}(\Xt)$ but we utilise the isomprphism between  $TM_{x}(\Xt)$ and $\Xt$ given the assumed `flat' geometry of $\Xt$.} $x,y\in \Xt$. Thus, if $D\phi_{t_0,t}(x,\om)$  is bounded away from zero for $t\in \Ic$, then $D\Phi_{t_0,t}^x(\om)^{-1} = D\Solm(x,\om)^{-1}$ $\p$\,-\,a.s. which implies that $D\Phi_{t_0,t}^{x}(\om)$ is nonsingular.

The rest of the proof exploits the It\^o formula for continuous semimartingales obtained from the two-point motion $\big\{(\phi_{t,t_0}(x_1,\omega),\phi_{t,t_0}(x_2,\omega))\,{:}\; (x_1,x_2)\,{\in}\, \Xt\otimes\Xt,\; t,t_0\,{\in}\, \Ic\big\}$.

First, consider a semi-martingale
 $V_t^{\scriptscriptstyle (N)}\!: \R^{Nd}\rightarrow\R^{Nd}$ with $V_{t,t_0}^{\scriptscriptstyle(N)} \!= (X^{x_1}_t, \dots,X^{x_N}_t)$, $(x_1,\cdots x_N)\in \otimes^{N}\Xt$, s.t. $X^{t_0,x_i}_t$ solves (\ref{gen_sde})  with the initial conditoin $x_i\in \Xt$.
The continuous dependence the solutions $X^{t_0,x_i}_t$ on the initial condition  leads to the $N$-point motion generated by  
$\big\{(\Solm(x_1,\om),\cdots , \Solm( x_N,\om))\!\!: (x_1,\cdots x_N)\in \otimes^{N}\Xt,\; t_0,t\in \Ic\big\}$ which may be interpreted as a stochastic flow generated  by `running' (\ref{gen_sde}) simultaneously from $N$ different initial conditions (see, e.g., \cite{Arnold1, Kunitabook, Baxendale2} for details). The $N$-point motion of stochastic flows is a process on $\R^{Nd}$ with the generator $\mathcal{L}^{\scriptscriptstyle(N)}$ induced by (\ref{Ngen}) and transition evolutions $(\mathcal{P}^{\scriptscriptstyle(N)}_{t_0,t})_{t \in \Ic}$ in the form~(\ref{calP}) and generated by the transition kernel 
\begin{align*}
P^{\scriptscriptstyle(N)}(t_0,(x_1,\dots,x_N);t ,E) &:= \p\big\{\om: (\Solm(x_1,\om), \dots,\Solm(x_N,\om)){\in} E\big\}, \quad t_0,t\in \Ic, \; E\in \Bb({\otimes}^{\scriptscriptstyle N}\Xt).
\end{align*}

\vspace{-0.2cm}
 \begin{theorem}{\rm [It\^o formula for continuous semimartingales \cite{Kunitabook}]. }\label{gIT}\rm
 Let $V_{t,t_0}^{\scriptscriptstyle (N)} = (X^{x_1}_t,\cdots, X_t^{x_N})$ be a continuous semi-martinagale with $X^{t_0,x_i}_t = \phi_{t,t_0}(x_i,\om)$. If $g\in\mathcal{C}^2(\R^{Nd})$, then $g(V_t^{\scriptscriptstyle(N)})$ is a continuous semi-martingale which satisfies
 \begin{align}\label{Ngen}
 g(V_{t,t_0}^{\scriptscriptstyle(N)})-g(V_{t_0}^{\scriptscriptstyle(N)}) = \sum_{i=1}^{\scriptscriptstyle N}\int_{t_0}^t\partial_{x_i}g(V_{s,t_0}^{\scriptscriptstyle(N)})dX_{s}^i+\frac{1}{2}\sum_{i,j=1}^{\scriptscriptstyle N}\int_{t_0}^t\partial^2_{x_ix_j}g(V_{s,t_0}^{\scriptscriptstyle(N)})d\big\langle X_s^i, X_s^j\big\rangle_s.
 \end{align}
 If $g\in \mathcal{C}^3(\R^{Nd})$, then we have 
 \begin{align}\label{SIt}
  g(V_{t,t_0}^{\scriptscriptstyle(N)})-g(V_{t_0}^{\scriptscriptstyle(N)}) = { \sum_{i=1}^{\scriptscriptstyle N}}\int_{t_0}^t\partial_{x_i}g(V_{s,t_0}^{\scriptscriptstyle(N)})\circ dX_s^i.
\end{align}  
 \end{theorem}

\smallskip
\noindent {\it Proof of Proposition \ref{prp_centfl}:} Consider the transition kernel for  the two-point motion \\$\big\{(\Solm(x_1,\om), \Solm(x_2,\om))\!: (x_1,x_2)\in \Xt\times\Xt,\; t,t_0\in \Ic\big\}$ given by 
\begin{align*}
P^{\scriptscriptstyle(2)}(t_0,(x_1,x_2);t ,E) &:= \p\{\om: (\Solm(x_1,\om), \Solm(x_2,\om))\in E\}, \quad t_0,t\in \Ic, \quad E\in \Bb(\Xt{\times}\Xt),
\end{align*}  
and the corresponding transition evolution $\mathcal{P}_{t_0,t}^{\scriptscriptstyle(2)}$  defined by 
\begin{align*}
\mathcal{P}^{\scriptscriptstyle(2)}_{t,t_0}g(x_1,x_2) &:= \int_{\R^{2d}}g(u,z)P^{\scriptscriptstyle(2)}(t_0,(x_1,x_2);t, du\otimes dz)=\E\big[g(\Solm(x_1,\om), \Solm(x_2,\om))\big],
\end{align*}
for $g\in \mathbb{M}(\Xt\times\Xt)$ and  $\mu\otimes\nu \in \PP(\Xt\times\Xt)$  with the dual 
\begin{align*}
(\mathcal{P}^{\scriptscriptstyle(2)*}_{t,t_0}\mu\otimes\nu)(E) = \int_{\Xt\times\Xt}\!\!\!P^{\scriptscriptstyle(2)}(t_0,(x_1,x_2); t, E)\mu\otimes\nu(dx_1,dx_2), \quad t\geqslant t_0, \quad E\in \Bb(\Xt{\times}\Xt), 
\end{align*}
where the extension of $g\in \mathbb{M}_\infty$ to $g\in \mathbb{M}$ is done in the standard way.

Now, let $X_t^{t_0,x_1}(\om)=\Solm(x_1,\om)$ and $X_t^{t_0,x_2}(\om)=\Solm(x_2,\om)$ be solutions of (\ref{gen_sde}) starting from two initial conditions $x_1, x_2\in \mathcal{M}$. Then,  according to the It\^o's formula (\ref{SIt}) with $g({x}_1,x_2) = x_1-x_2$, the process  $V^{\scriptscriptstyle(2)}_{t,t_0} = \Solm(x_1,\cdot)-\Solm(x_2,\cdot)$ is a semimartingale  and it can be represented as 
 \begin{align}\label{App_gIT}
 V_{t,t_0}^{\scriptscriptstyle(2)}= {x}_1{-}x_2+ \!\!\int_{t_0}^t\!\Big[b(s,\phi_{s,t_0}(x_1)) - b(s,\phi_{s,t_0}(x_2))\Big]ds+\!\!\int_{t_0}^t\!\Big[\sigma(s,\phi_{s,t_0}(x_1))-\sigma(s,\phi_{s,t_0}(x_2))\Big]{\circ} \,dW_s .
 \end{align}
Next, for $g\in \mathcal{C}^2(\Xt\,{\times}\, \Xt)$, define 
\begin{align*}
\mathcal{P}^{\scriptscriptstyle(2)}_{t,t_0}\,g(x_1 -x_2):=\E\Big [g\big(V_{t,t_0}^{\scriptscriptstyle(2)}\big)\Big |V_{t,t_0}^{\scriptscriptstyle(2)} = x_1 -x_2\Big] = \int_{\Xt\times \Xt} g(y-u)P^{\scriptscriptstyle(2)}(t_0, (x_1, x_2)); t, dy\otimes du).
\end{align*}
Furthermore, the generator of the two-point motion $\{(\Solm(x,\om), \Solm(y, \om)): t_0,t\in \Ic, (x,y)\in\Xt\times\Xt\}$ is given by 
\begin{align*}
\mathcal{L}^{\scriptscriptstyle(2)}_t:=\sum_{i=1}^d {b}_i(t,x)\partial_{x_i}+\sum_{i=1}^d {b}_i(t,y)\partial_{y_i} +{\textstyle \frac{1}{2}}\sum_{i,j=1}^{2d} \sum_{k=1}^d\sigma_{k,i}(t,x,y)\sigma_{k,j}(t,x,y)\partial^2_{x_i x_k},
\end{align*}
where $\sigma_{k,i}(t,x,y) := \sigma_{ki}(t,x)$ and $\sigma_{k,i+d}(t,x,y) := \sigma_{kj}(t,y), $ $i=1,2,\cdots,d$. Finally, application of It\^o's formula (Theorem \ref{gIT}) to  the process $ g\left(V^{\scriptscriptstyle(2)}_{t,t_0}(\om)\right)=g\left(\phi_{t,t_0}(x_1,\om)-\phi_{t,t_0}(x_2,\om)\right), $ $t_0,t\in\Ic$, for $g\in \mathcal{C}_\infty^2(\Xt\times\Xt)$, gives a simplified form of the generator $ \mathcal{L}_t^{\scriptscriptstyle(2)}$ as 
\begin{align*}
\mathcal{L}^{\scriptscriptstyle(2)}_tg(x_1-x_2) &= \left({b}(t,x_1)-{b}(t,x_2)\right)\partial_x \,g(x_1-x_2)\notag\\[.2cm]
&\hspace{2cm}+{\textstyle \frac{1}{2}}\text{tr}\left( \big[\sigma(t, x_1)-\sigma(t,x_2)\big]^T\mathcal{H}_g(x_1-x_2)\big[\sigma(t,x_1)-\sigma(t,x_2)\big]\right).
\end{align*}
Setting $x_1 = x+y$ and $x_2 = x$ in the above two expressions for  $\mathcal{P}^{\scriptscriptstyle (2)}_{t,t_0}$ and $\mathcal{L}^{\scriptscriptstyle(2)}_t$ yields, $\mathcal{P}^{\Phi^x}_{t,t_0}$ in~(\ref{2pme}) and $\mathcal{L}^x_t$ in~(\ref{2pp1}).
\qed

\subsection{Proof of Lemma~\ref{alpha_bnd1}}\label{app_alpha_bnd1}
Given the same assumptions as in Theorem~\ref{Dphi_posit} we have for $t,t_0\in \Ic$
\begin{align*}
\D_{\varphi}(\mu_t\|\mu_{t_0}) =& \int_{\Xt}\varphi\left(\rho_t(x)/\rho_{t_0}(x)\right)\rho_{t_0}(x)m_d(dx)\\
&\hspace{-1.15cm}\underset{\textrm{Thm \ref{mut_rnd},\,\ref{FKPt1}}}{=} \int_{\Xt}\varphi\Big(\E[\pi_{t_0,t}(x)]\Big)\rho_{t_0}(x)m_d(dx)\\[.2cm]
&\hspace{-1.15cm}\underset{\text{Jensen's ineq.}}{\leqslant} \int_{\Xt}\E[\varphi\big(\pi_{t_0,t}(x)\big)]\rho_{t_0}(x)m_d(dx)\\
&\hspace{-.8cm}\underset{\textrm{Thm \ref{al_pi}}}{=} \int_{\Xt}\E\Big[\varphi\Big(\alpha_{t_0,t}^{-1}\big(\Solm^{-1}(x)\big)\Big)\Big]\rho_{t_0}(x)m_d(dx)\\[.2cm]
&\hspace{-1.1cm}\underset{\text{Fubini's thm.}}=\E\bigg[\int_{\Xt}\varphi\Big(\alpha_{t_0,t}^{-1}\big(\Solm^{-1}(x)\big)\Big)\rho_{t_0}(x)m_d(dx)\bigg]\\[.2cm]
&\hspace{-.4cm}=\E\bigg[ \int_{\Xt}\varphi\big(\alpha_{t_0,t}^{-1}(\xi)\big)\rho_{t_0}(\Solm(\xi))\big|\!\det D\Solm(\xi)\big|\rho_{t_0}^{-1}(\xi)\rho_{t_0}(\xi)m_d(d\xi)\bigg] \\[.2cm]
&\hspace{-.4cm}= \E\bigg[\int_{\Xt}\varphi\big(\alpha_{t,t_0}^{-1}(\xi)\big)\alpha_{t_0,t}(\xi)\rho_{t_0}(\xi)m_d(d\xi)\bigg]=\int_{\Xt}\E\big[\varphi_{\ddagger}\big(\alpha_{t_0,t}(\xi)\big)\big]\mu_{t_0}(d\xi),
\end{align*}
where $\E[\pi_{t_0,t}(x)]:=\int_\Om \pi_{t_0,t}(x,\om)\p(d\om)$, $\E[\alpha_{t_0,t}(x)]:=\int_\Om \alpha_{t_0,t}(x,\om)\p(d\om)$, and similarly \\$\E\big[\varphi\big(\alpha_{t_0,t}^{-1}\big(\Solm^{-1}(x)\big)\big)\big] := \int_\Om\big[\varphi\big(\alpha_{t_0,t}^{-1}\big(\Solm^{-1}(x,\om),\om\big)\big)\big]$ to simplify notation.

\medskip

\subsection{Proof of Theorem \ref{local-bound}}\label{app_local-bound}
First, note that  for $\mu_{t_0}\in \PP(\Xt)$ with strictly positive Lebesgue density $\rho_{t_0}>0$ the (random) density process  $\{\alpha_{t_0,t}(x,\om)\!: \,t_0,t\in \Ic\}$, $\alpha_{t_0,t} = d\Check\Pi_{t_0,t}/d\mu_{t_0}$, given in (\ref{alph_Ga})  is a  continuous semi-martingale w.r.t.~the natural (complete, right-continuous) filtration $(\F_{t_0}^{\,t})_{t\in \Ic}$, since by It\^o's formula it satisfies 
 (see also (\cite[Corollary 4.3.5]{Kunitabook}))
\begin{align}\label{alph_dw}
\alpha_{t_0,t}(x,\om) = 1 + \int_{t_0}^t\alpha_{t_0,s}(x,\om) G(s,\phi_{t_0,s}(x,\om))dW_s +\int_{t_0}^t\alpha_{t_0,s}(x,\om)\frac{\check{\LG}_s^*\rho_{t_0}(\phi_{t_0,s}(x,\om))}{\rho_{t_0}(\phi_{t_0,s}(x,\om))}ds.
\end{align}  
Moreover,  note that for any  convex $\varphi(u)$, $u>0$, the function $\varphi_\ddagger(u) = u\varphi(u^{-1})$ is also convex and locally bounded. In what follows we consider functions $\varphi$ satisfying (\ref{Normality}) so that $\varphi_\ddagger(1)=0$. Thus, by It\^o-Tanaka-Meyer formula (e.g., \cite{Yor}), we have for any $t\in \Ic$
\begin{align}\label{inter11}
\varphi_{\ddagger}(\alpha_{t_0,t}) &= \varphi_{\ddagger}(\alpha_{t_0,t_0})+\int_{t_0}^t D^{-}\varphi_{\ddagger}(\alpha_{t_0,s})d\alpha_{t_0,s}
\notag\\&\hspace{2cm}+\frac{1}{2}\int_{(0,\infty)}L^{\ell}_t(\alpha)\varphi^{\prime\prime}_{\ddagger}(d\ell)+\sum_{t_0\leqslant s\leqslant t}\big(\Delta\varphi_{\ddagger}(\alpha_{t_0,s})-\varphi'_{\ddagger}(\alpha_{t_0,s^-})\Delta\alpha_{t_0,s}\big)\notag\\
&= \int_{t_0}^t D^{-}\varphi_{\ddagger}(\alpha_{t_0,s})d\alpha_{t_0,s}
\notag\\&\hspace{2cm}+\frac{1}{2}\int_{(0,\infty)}L^{\ell}_t(\alpha)\varphi^{\prime\prime}_{\ddagger}(d\ell)+\sum_{t_0\leqslant s\leqslant t}\Delta\varphi_{\ddagger}(\alpha_{t_0,s}), 
\end{align}
where $L^{\ell}_t(\alpha)$ is the local time of $\alpha_{t_0,t}(x,\om)$ at level $\ell\geqslant 0$ (e.g., \cite{Yor}) 
given by 
\begin{align*}
L_t^\ell(\alpha) = 2\Big[(\alpha_{t,t_0}-\ell)_{+}-(\alpha_{t_0,t_0}-\ell)_{+} -\int_{t_0}^t\I_{\{\alpha_{t_0,s}>\ell\}}d\alpha_{t_0,s}-\sum_{t_0\leqslant s\leqslant t}\Delta(\alpha_{t_0,s}-\ell)_+\Big],
\end{align*}
$\Delta(\ccdot)$ denotes a pure jump of the argument and the left-sided derivative $D^{-}\varphi_{\ddagger}$ is bounded on $\Ic$ due to the convexity of $\varphi_{\ddagger}$ on $(0, \infty)$.

Taking the expectation of (\ref{inter11}) conditioned on $\F_{t_0}^t,$  we have that 
\begin{align*}
\E\big[\varphi_{\ddagger}(\alpha_{t_0,t})\big] &= \sum_{t_0\leqslant s\leqslant t}\E\big[\Delta\varphi_{\ddagger}(\alpha_{t_0,s})\big]+\E\left[\int_{t_0}^tD^{-}\varphi_{\ddagger}(\alpha_{t_0,s})\alpha_{t_0,s}\frac{\check{\LG}_s^*\rho_{t_0}(\phi_{t_0,s})}{\rho_{t_0}(\phi_{t_0,s})}ds\right]\\[.2cm]
&\qquad +\frac{1}{2}\E\bigg[\int_{(0,\infty)}L^\ell_t(\alpha)\varphi^{\prime\prime}_{\ddagger}(d\ell)\bigg],
\end{align*}
which, combined with the inequality (\ref{N121}), leads to the first assertion
\begin{align*}
\Df(\mu_t||\mu_{t_0})&\leqslant  \sum_{t_0\leqslant s\leqslant t}\E\big[\Delta\varphi_{\ddagger}(\alpha_{t_0,s})\big]+\int_{\Xt}\E\left[\int_{t_0}^tD^{-}\varphi_{\ddagger}(\alpha_{t_0,s})\alpha_{t_0,s}\frac{\check{\LG}_s^*\rho_{t_0}(\phi_{t_0,s})}{\rho_{t_0}(\phi_{t_0,s})}ds\right]d\mu_{t_0}\\[.2cm]
&\qquad +\frac{1}{2}\int_{\Xt}\E\bigg[\int_{(0,\infty)}L^\ell_t(\alpha)\varphi^{\prime\prime}_{\ddagger}(d\ell)\bigg] d\mu_{t_0}.
\end{align*}
\smallskip
Next, if $\varphi$ is twice continuously differentiable on $(0,\infty)$, we obtain
\begin{align*}
\hspace{2cm}\Df(\mu_t\|\mu_{t_0})&\leq \int_{\Xt}\E\left[\int_{t_0}^t\varphi^\prime_{\ddagger}(\alpha_{t_0,s})\alpha_{t_0,s}\frac{\check{\LG}_s^*\rho_{t_0}(\phi_{t_0,s})}{\rho_{t_0}(\phi_{t_0,s})}ds\right]\mu_{t_0}(dx)\\[.2cm]
&\qquad +\frac{1}{2}\int_{\Xt}\E\bigg[\int_{t_0}^t \varphi_{\ddagger}^{\prime\prime}(\alpha_{t_0,s})d\langle \alpha_{t_0,s}\rangle\bigg]\mu_{t_0}(dx)\\[.2cm]
&=\int_{\Xt}\E\left[\int_{t_0}^t\varphi^\prime_{\ddagger}(\alpha_{t_0,s})\alpha_{t_0,s}\frac{\check{\LG}_s^*\rho_{t_0}(\phi_{t_0,s})}{\rho_{t_0}(\phi_{t_0,s})}ds\right]\mu_{t_0}(dx)\\[.2cm]
&\qquad +\frac{1}{2}\int_{\Xt}\E\bigg[\int_{t_0}^t \varphi_{\ddagger}^{\prime\prime}(\alpha_{t_0,s})\alpha^2_{t_0,s}G^2(s,\phi_{t_0,s})ds\bigg]\mu_{t_0}(dx). \hspace{2cm}\qed
\end{align*}


\newpage
\section*{Glossary}\label{glossary}

\vspace{-.2cm}\noindent Here, we list further definitions and notation which recurs throughout the paper.

\smallskip
\noindent{\small \bf (2) Probability spaces and function spaces}
{\small
\vspace*{-.0cm}
\begin{itemize}[leftmargin=0.4cm]
\item[\tiny$\bullet$]{\it Wiener space.} We  fix the probability space $\PS$ as the Wiener space, i.e.,
$\Om\simeq\mathcal{C}_{0}(\Ic;\R^m)$, $m\in \N$, $\Ic:= \tint\subset\R$, $T>0$, is a subspace of continuous functions $\mathcal{C}(\R;\R^m)$ which are zero at $t_0\in \Ic$. $\F$ is the Borel  $\mathfrak{S}$-algebra generated by open subsets  in the compact-open topology on $\Om$ defined~via
\begin{align*}
\varrho(\om,\hat{\om})= \sum_{\ell=0}^\infty\frac{1}{2^\ell}\frac{\Vert \om -\hat{\om}\Vert_\ell}{1+\Vert \om -\hat{\om}\Vert_\ell}, \qquad \Vert \om -\hat{\om}\Vert_{\ell} := \sup_{t\in [-\ell, \,\ell\,]}\vert \om(t) -\hat{\om}(t)\vert, \quad \om,\hat\om\in \Omega, 
\end{align*}
with $|{\ccdot}|$ the Euclidean norm on $\R^m$. Finally, $\p$ is the Wiener measure on $\mathcal{F}$.

\vspace{.05cm}
\item[\tiny$\bullet$]  $\mathcal{W}_d:=\mathcal{C}(\Ic,\Xt)$ is the path space defined over $\Xt$. The Borel $\mathfrak{S}$-algebra  $\mathcal{B}(\mathcal{W}_d)$ on $\mathcal{W}_d$  are defined analogously to those in the Wiener space.

\vspace{.05cm}
\item[\tiny$\bullet$] For $f\,{:}\; \Xt\rightarrow\R$, where $(\Xt, \mathcal{B}(\Xt))$ is a Polish  space equipped with a Borel $\mathfrak{S}$\,-\,algebra, 
the following function spaces are relevant:
\begin{itemize}[leftmargin=0.4cm]
 \item $\mathbb{M}_\infty(\Xt)$ space of bounded Borel measurable functions  $\mathbb{M}(\Xt)$ on $\Xt$. 
 \item $\mathbb{M}^{+}(\Xt)$ space of non-negative Borel measurable functions on $\Xt$.
 \item $\mathcal{C}_\infty(\Xt)$ space of bounded continuous functions on $\Xt$.
\item  $\mathcal{C}^l(\Xt)$, $l\geqslant 1$, space of  $l$-times continuously differentiable functions on $\Xt$. 
\item  $\mathcal{C}^k_\infty(\Xt)$,  $l\geqslant 1$, functions in $\mathcal{C}^l(\Xt)$ which are bounded with bounded derivatives up to order $l$ on $\Xt$. 
 \item $\mathcal{C}_c^{+}(\Xt)$ space continuous non-negative functions on $\Xt$ with compact supports.
\item  $\mathcal{C}_c^{\infty}(\Xt)$ space of smooth functions on $\Xt$ with compact support.
 \end{itemize}

\item[\tiny$\bullet$] $\tilde{\mathcal{C}}^{l,\delta}(\Xt)$,  is the space of functions $f\,{:} \;\Xt\,{\rightarrow}\, \Xt$ with the countable family of semi-norms
\begin{align*}
& \tilde{\Vert} f\tilde{\Vert}_{l,\delta;N}:= \Vert f\Vert_{l;N}+\sum_{\vert \alpha\vert =l}\sup_{x,y\in \textsf{B}_{N}\!,\, x\neq y}\frac{\vert D^{\alpha}f(x) -D^{\alpha}f(y)\vert}{\vert x-y\vert^{\delta}}<\infty, \quad 0<\delta\leqslant 1,\;N_1\in\mathbb{N},\hspace{.6cm}\\[.1cm]
&\tilde{\Vert} f\tilde{\Vert}_{l;N}:= \sup_{x\in\Xt}\frac{\vert \langle f(x),x\rangle\vert}{1+\vert x\vert^2}+\sum_{1\leqslant \vert \alpha\vert \leqslant l}\sup_{x\in \textsf{B}_{N}}\vert D^{\alpha}f(x)\vert,
\end{align*}

\vspace{-.3cm}\noindent where $\textsf{B}_N:=\{x\in\Xt: |x|\leqslant N\}$, $\displaystyle D^{\alpha}f(x): = \frac{\partial^{\vert \alpha\vert}f}{(\partial x_1)^{\alpha_1}\cdots(\partial x_d)^{\alpha_n}}, \; \vert \alpha\vert := \sum_{i=1}^n\alpha_i$, $\alpha_i\in \mathbb{N}_0$, and \;$D^0\equiv \textrm{Id}$.

\smallskip
\item[\tiny$\bullet$] $\bar{\mathcal{C}}^{l,\delta}(\Xt)$,  is the space of functions $f\,{:} \;\Xt\,{\rightarrow}\, \Xt$ with the countable family of semi-norms
\begin{align*}
& \tilde{\Vert} f\tilde{\Vert}_{l,\delta;N}:= \Vert f\Vert_{l;N}+\sum_{\vert \alpha\vert =l}\sup_{x,y\in \textsf{B}_{N}\!,\, x\neq y}\frac{\vert D^{\alpha}f(x) -D^{\alpha}f(y)\vert}{\vert x-y\vert^{\delta}}<\infty, \quad 0<\delta\leqslant 1,\;N_1\in\mathbb{N},\hspace{.6cm}\\[.1cm]
&\tilde{\Vert} f\tilde{\Vert}_{l;N}:= \sup_{x\in\Xt}\frac{\vert  f(x)\vert}{1+\vert x\vert}+\sum_{1\leqslant \vert \alpha\vert \leqslant l}\sup_{x\in \Xt}\vert D^{\alpha}f(x)\vert,
\end{align*}

\item[\tiny$\bullet$] $\mathcal{C}\big(\Ic; \bar{\mathcal{C}}^{l,\delta}(\Xt)\big)$, $\Ic\subseteq \R$,  is the set of all continuous fields $f\,{:}\; \R\times\Xt\rightarrow\Xt$ such that $f(t,\ccdot)\in \bar{\mathcal{C}}^{l,\delta}(\Xt)$. 

\end{itemize}

\noindent{\small \bf (3) Frequently used notation}

\vspace*{-0.1cm}
\begin{itemize}[leftmargin=0.4cm]

\vspace{.1cm}
\item[\tiny$\bullet$]  $\PP(\Xt)$ is a set of all Borel probability measures on $\Xt$.

\vspace{.1cm}
\item[\tiny$\bullet$] $\Df\big(\mu\|\nu\big)$ is a $\varphi$-divergence between measures in $\PP(\Xt)$; $\varphi$  is a strictly convex function (see~(\ref{dphi})).

\vspace{.1cm}
\item[\tiny$\bullet$] $\Df^{t_0,t}\big(\mu\|\nu\big) = \scaleobj{.9}{|t-t_0|^{-1}}\hspace{.004cm}\Df^{t_0,t}\big(\mu\|\nu\big)$ is a $\varphi$-divergence rate between measures in $\PP(\Xt)$.

\vspace{.1cm}
\item[\tiny$\bullet$] $\big\{\phi_{s,t}(\ccdot,\om)\,{:}\; s,t\in\Ic\big\}$, $\om\in \Om$, denotes a stochastic flow on $\Xt$  (usually, but not exclusively) generated by the  SDE (\ref{gen_sde}); see~\S\ref{s_setup}. $\{\Phi^x_{t_0,t}\}_{t\in \Ic}$,  $\Phi^x_{t_0,t}(y,\om) = \phi_{t_0,t}(x+y,\om)-\phi_{t_0,t}(x,\om)$, is the two-point motion associated with $\{\phi^x_{t_0,t}\}_{t\in \Ic}$.

\vspace{.1cm}
\item[\tiny$\bullet$] $(\brr,\sigma)$ are the drift and diffusion coefficients of the SDE (\ref{gen_sde}). $\{\sigma_{k}\}_{k=1}^m$ stand for  columns of the matrix field $\sigma$ with coefficients $\sigma_{ik}$.

\vspace{.1cm}
\item[\tiny$\bullet$] $ {b}_i(t,x) \,{:=}\, \brr_i(t,x)\,{+}\,c_i(t,x)$,  $c_i(t,x) := \frac{1}{2}\sum_{k,j=1}^{m,d} \sigma_{jk}(t,x)\partial_{x_j} \sigma_{ik}(t,x)$, $a_{ij} :=  \sum_{k=1}^m \sigma_{ik}\sigma_{jk}$, $i=1,\dots,\ell$, is the Stratonovich-corrected drift in the It\^o SDE (\ref{gen_sde_ito}). 

\vspace{.1cm}
\item[\tiny$\bullet$] $\Vert \sigma\Vert^2_{\textsc{hs}} :=\sum_{i=1}^n \sum_{k=1}^m\vert \sigma_{ik}\vert^2$ is the Hilbert-Schmidt (or Frobenius) norm of the matrix field $\sigma$.

\vspace{.1cm}
\item[\tiny$\bullet$] $t\mapsto \phi_{t_0,t}(x,\om)$, $t\in \Ic$, $\phi_{t_0,t_0}(x,\om)=x\in \Xt$ is a random path of the original dynamical system on $\Xt$.

\vspace{.1cm}
\item[\tiny$\bullet$] $\mu_t\in \PP(\Xt)$ is the time-marginal probability measure associated with the dynamics on $\Xt$.  For dynamics induced by the SDE (\ref{gen_sde}), $\mu_t$ solves (weakly) the forward Kolmogorov equation (\ref{fPDE}). 

\vspace{.1cm}
\item[\tiny$\bullet$] $\rho_t$ is the density of $\mu_t$ w.r.t.~Lebesgue measure $m_\ell$ on $\Xt$ (whenever $\mu_t\ll m_\ell$).

\vspace{.1cm}
\item[\tiny$\bullet$] $\tilde\mu^{x,\varepsilon}_{t_0}$ is a regularised uniform measure in $\PP(\Xt)$ localised on the ball $B_\varepsilon(x)$; see Definition~\ref{gmolf}.

\vspace{.1cm}
\item[\tiny$\bullet$] $(\mathcal{P}_{t_0,t})_{t\geqslant s}$ is a family of transition evolutions  induced by $\{\phi_{t_0,t}\}_{t\in \Ic}$, and acting on $f\in \mathbb{M}(\Xt)$; see (\ref{calP}).

\vspace{.1cm}
\item[\tiny$\bullet$] $(\mathcal{P}^{\Phi^x}_{t_0,t})_{t\geqslant s}$ is a family of transition evolutions  induced by $\{\Phi_{t_0,t}\}_{t\in \Ic}$.

\vspace{.1cm}
\item[\tiny$\bullet$]
 $\mathcal{P}^*_{t_0,t}$ and $\mathcal{P}^{\Phi^x*}_{t_0,t}$ are  $L^1$ duals of  $\mathcal{P}_{t_0,t}$ and $\mathcal{P}^{\Phi^x,}_{t_0,t}$ acting on probability measures in $\PP(\Xt)$;~see (\ref{P*mu}).

\vspace{.1cm}
\item[\tiny$\bullet$] $\Ms_{t_0,x}$ solves  the martingale problem for the operator $\mathcal{L}_t$ starting at $(t_0,x)\in \Ic\times\Xt$.
   $\Ms_{t_0,x}$ is identified with a path space probability measure on $(\mathcal{W}_d, \mathcal{B}(\mathcal{W}_d))$ s.t.~$\Ms_{t_0,x}((t_0,x))=1$.
 
\vspace{.1cm}
\item[\tiny$\bullet$] $\Ms_{t_0}(d\om):=\int_\Xt \Ms_{t_0,x}(d\om)\mu_{t_0}(dx)$ is a  path space probability measure on $(\mathcal{W}_d, \mathcal{B}(\mathcal{W}_d))$ which such that  $\Ms_{t_0}\simeq \mu_{t_0}\otimes\p$ and (formally) $\Ms_{t_0}\circ\phi_{t_0,t}^{-1} = \mu_t\in \PP(\Xt)$; see (\ref{k_lift}).

\vspace{.1cm}
\item[\tiny$\bullet$] $\E[f(\phi_{t_0,t}(x))]:=\int_{\Om} f(\phi_{t_0,t}(x,\om))\Ms_{t_0,x}(d\om)=\int_{\Om} f(\phi_{t_0,t}(x,\om))\p(d\om)$  denotes an `observable' based on $f\in \mathbb{M}(\Xt)$, and defined on the paths $t\mapsto\phi_{t_0,t}(x,\om)$.

\vspace{.1cm}
\item[\tiny$\bullet$]
$\E^{\mu_t}[f]:= \int_\Xt f(x)\mu_t(dx)= \int_\Xt\int_{\Om} f(\phi_{t_0,t}(x,\om))\Ms_{t_0,x}(d\om)\mu_{t_0}(dx)$ is an observable based on $f\in \mathbb{M}(\Xt)$, and evaluated on  paths $t\mapsto\phi_{t_0,t}(x,\om)$.

\vspace{.1cm}
\item  $\mathcal{R}^{\Phi}_{\varphi,\varepsilon}(x,t_0,t)$ is a $\varphi$-divergence rate field at time $t\in \Ic$ given by $\D_\varphi^{t_0,t}\big(\mathcal{P}^{\Phi^x*}_{t_0,t}\tilde\mu_{t_0}^{x,\varepsilon}\|\tilde\mu_{t_0}^{x,\varepsilon}\big)$ (see (\ref{xlocrt})).

\end{itemize}
}

\medskip
{\small

\end{document}